\documentclass[pdflatex,sn-mathphys-num]{sn-jnl}% Math and Physical Sciences Numbered Reference Style
\usepackage{graphicx}
\usepackage{placeins}   % \FloatBarrier, used in 6-comp-results.tex

\usepackage{amsmath,amssymb,amsfonts}

\hypersetup{hypertexnames=false}
\usepackage{amsthm}
\usepackage{mathtools}
\usepackage{xcolor}
\usepackage{multirow}
\usepackage{array}      % column width in table
\usepackage{longtable}  % the per-instance reference-point table spans pages
\usepackage{subcaption} % subfigures (a)/(b) in the results section

\DeclareCaptionFont{sncaptionfont}{\fontsize{8bp}{9.5bp}\selectfont\rmfamily}
\DeclareCaptionLabelSeparator{snfigsep}{\hskip.7em}
\DeclareCaptionLabelSeparator{sntabsep}{\hskip2mm}
\newlength{\sncapwidth}
\AddToHook{env/table/begin}{\setlength{\sncapwidth}{\textwidth}}
\AddToHook{env/sidewaystable/begin}{\setlength{\sncapwidth}{0pt}}
\usepackage{algorithm}
\usepackage{algorithmicx}
\usepackage{algpseudocode}

\usepackage{tikz}
\usetikzlibrary{backgrounds,patterns,arrows.meta,calc,positioning}
\usetikzlibrary{svg.path}

\theoremstyle{thmstyleone}
\newtheorem{theorem}{Theorem}[section]
\newtheorem{corollary}[theorem]{Corollary}
\newtheorem{proposition}[theorem]{Proposition}
\newtheorem{lemma}[theorem]{Lemma}
\newtheorem{assumption}{Assumption}

\theoremstyle{thmstylethree}
\newtheorem{definition}[theorem]{Definition}

\theoremstyle{thmstyletwo}

\newcommand{\refpt}{\ensuremath{p^\mathrm{ref}}}

\DeclareMathOperator{\MAS}{MAS}

\DeclareMathOperator{\EM}{EM}

\let\AA\relax
\DeclareMathOperator{\AA}{AA}
\DeclareMathOperator{\HV}{HV}
\DeclareMathOperator{\Top}{Top}

\DeclareMathOperator{\argmin}{argmin}
\DeclareMathOperator{\lex}{lex}

\newcommand{\R}{\ensuremath{\mathbb{R}}}
\newcommand{\N}{\ensuremath{\mathbb{N}}}
\newcommand{\X}{\ensuremath{\mathcal{X}}}
\newcommand{\Y}{\ensuremath{\mathcal{Y}}}

\newcommand{\Init}{\texttt{InitSubset}}

\newcommand{\UB}{\mathrm{UB}}
\newcommand{\LB}{\mathrm{LB}}
\newcommand{\EMmax}{\ensuremath{\mathrm{EM}_\mathrm{max} }}
\newcommand{\CCG}{C\&CG}
\newcommand{\Rs}[1][s]{\ensuremath{R_{#1}}}
\newcommand{\OPT}{\ensuremath{\mathrm{OPT}}}
\newcommand{\point}{p}

\newcommand{\PND}{P_\mathrm{ND}}
\newcommand{\PSOL}{P_{\mathrm{E}}}

\newcommand{\addsn}{\ensuremath{\mathrm{add}^{s,n}}}
\newcommand{\XY}{\ensuremath{\X \times \Y^s}}
\newcommand{\TG}{\ensuremath{\mathrm{TG}}}

\newcommand{\round}[1]{\ensuremath{\lfloor#1\rceil}}

\newcommand{\zenodolink}{\href{https://doi.org/10.5281/zenodo.23127921}{Zenodo}}

\definecolor{fillblue}{RGB}{173,216,230}
\definecolor{fillred}{RGB}{255,127,127}
\definecolor{pointred}{RGB}{153,0,0}
\definecolor{pointblue}{RGB}{0,0,255}
\definecolor{fillpurple}{RGB}{127,221,127}
\definecolor{pointpurple}{RGB}{0,110,0}
\definecolor{fillorange}{RGB}{255,190,120}
\definecolor{pointorange}{RGB}{204,102,0}
\definecolor{fillgreen}{RGB}{205,160,235}
\definecolor{pointgreen}{RGB}{120,40,160}
\definecolor{mygray}{RGB}{150,150,150}

\pgfdeclarelayer{foreground}
\pgfsetlayers{main,foreground}

\newif\ifconnectfront
\connectfronttrue
\colorlet{connectcolor}{black!45}        % line color (light gray)
\newcommand{\connectfrontwidth}{0.9pt}   % line thickness
\newcommand{\connectfrontopacity}{0.6}   % line alpha
\tikzset{connectline/.style={draw=connectcolor,
		line width=\connectfrontwidth, opacity=\connectfrontopacity}}

\newcommand{\drawfrontline}[1]{%
	\ifconnectfront
		\xdef\cfpath{}%
		\foreach \cfa/\cfb [count=\cfc] in #1 {%
			\ifnum\cfc=1\relax
				\xdef\cfpath{(\cfa,\cfb)}%
			\else
				\xdef\cfpath{\cfpath -- (\cfa,\cfb)}%
			\fi
		}%
		\draw[connectline] \cfpath;%
	\fi
}

\newcommand{\drawfrontlineheights}[2]{% #1 = heights list, #2 = x-coordinate csname prefix
	\ifconnectfront
		\xdef\cfpath{}%
		\foreach \cfh [count=\cfi] in #1 {%
			\pgfmathtruncatemacro{\cfshow}{\cfh>0 ? 1 : 0}%
			\ifnum\cfshow=1\relax
				\pgfmathsetmacro{\cfy}{7-\cfh}%
				\edef\cfx{\csname #2\cfi\endcsname}%
				\ifx\cfpath\empty
					\xdef\cfpath{(\cfx,\cfy)}%
				\else
					\xdef\cfpath{\cfpath -- (\cfx,\cfy)}%
				\fi
			\fi
		}%
		\draw[connectline] \cfpath;%
	\fi
}

\makeatletter
\newcommand{\algmargin}{\the\ALG@thistlm}
\makeatother

\begin{document}

\title[Hypervolume-based \CCG for bi-objective two-stage robust MIPs]{Hypervolume-based column-and-constraint generation for bi-objective two-stage robust mixed-integer programs}

%%=============================================================%%
%% GivenName	-> \fnm{Joergen W.}
%% Particle	-> \spfx{van der} -> surname prefix
%% FamilyName	-> \sur{Ploeg}
%% Suffix	-> \sfx{IV}
%%=============================================================%%

\author*[1]{\fnm{Johannes} \sur{Kager}}\email{johannes.kager@tum.de}

\author[1]{\fnm{Clemens} \sur{Thielen}}\email{clemens.thielen@tum.de}

\affil[1]{\orgdiv{Professorship of Optimization and Sustainable Decision Making}, \orgname{TUM~Campus Straubing for Biotechnology and Sustainability},
	\orgaddress{\street{Am~Essigberg 3}, \city{Straubing},
		\postcode{94315},
        \country{Germany}}}

\abstract{Extending column-and-constraint generation (\CCG{}) for two-stage robust mixed-integer programs from the single-objective to the bi-objective setting is challenging: for a fixed first-stage decision and scenario, the second stage yields an entire nondominated set rather than a single objective value. We develop a \CCG{} algorithm that uses the hypervolume indicator---one of the most widely studied quality indicators in multi-objective optimization---to evaluate these sets and seeks a first-stage decision maximizing the worst-case hypervolume over a finite set of scenarios. To avoid the nonlinearities of an exact hypervolume formulation, we introduce a slab decomposition of the objective space that yields a linear mixed-integer master problem with an explicit bound on the discretization error. This bound enables a computable optimality-gap certificate and supports an adaptive algorithm that alternates between adding scenarios and refining the decomposition. We prove the correctness of the algorithm and establish sufficient conditions for its finite termination. Computational experiments on two-stage robust selection and assignment problems show that only a small number of scenarios need to be included in the master problem, even when the full scenario set contains 128 scenarios. Moreover, we identify the number of initial slabs as a key factor influencing the running time and develop a simple, instance-dependent rule for choosing the initial slab width that outperforms fixed choices on average for the selection problem and remains competitive for the assignment problem.}


\maketitle

\section{Introduction}\label{sec:intro}

Robust optimization~\citep{ben2009robust} addresses decision-making under uncertainty by seeking solutions that remain feasible and perform well under uncertainty represented by a prescribed uncertainty set. 
In this paper, we consider two-stage robust optimization~\citep{bib27ben2004adjustable} for mixed-integer programs (MIPs). This setting distinguishes between \emph{first-stage decisions} that must be fixed before the uncertainty is revealed and \emph{second-stage decisions} that can be postponed until the realization of uncertain parameters is known. 
We restrict ourselves to finite uncertainty sets, whose elements we call \emph{scenarios}.

A standard solution paradigm for single-objective two-stage robust problems with finite uncertainty sets is column-and-constraint generation (\CCG{})~\citep{zeng,goerigk2025asbp}. \CCG{} iteratively restricts attention to a subset of scenarios, solves a \emph{master problem} to find an optimal first-stage decision for this restricted scenario set, and identifies a worst scenario by solving the \emph{second-stage problems}---the problems of optimizing the second-stage decisions given the current first-stage decision. Here, a scenario is called \emph{worst} if it yields the worst objective value among all scenarios. A worst scenario is then added to the considered subset, and the process repeats until a termination criterion is met. A central benefit of \CCG{} is its favorable scaling with the number of scenarios: in many settings, only a small subset of critical scenarios needs to be included before convergence, even if the full uncertainty set is large.

In this paper, we extend this framework to \emph{bi-objective} two-stage robust MIPs, which involve two usually conflicting objective functions that are optimized simultaneously. In the single-objective case, both the master problem and the second-stage problems are scalar-valued optimization problems. In a bi-objective problem, however, there is typically no single solution that optimizes both objectives simultaneously. In this situation, the solutions of interest are the so-called \emph{efficient} solutions~\citep{ehrgott2005multicriteria}, each of which is characterized by the property that no other feasible solution exists that is not worse in any of the two objectives and strictly better in at least one. The image of an efficient solution (the pair of its objective values) in the two-dimensional objective space is called a \emph{nondominated image} or a \emph{nondominated point}, and the set of all nondominated points is called the \emph{nondominated set}~\citep{ehrgott2005multicriteria}. The goal in a bi-objective problem then typically consists of computing all nondominated points together with one corresponding efficient solution for each of these points. As a consequence, neither component of the \CCG{} framework can be directly formulated: the master problem has no well-defined scalar objective over first-stage decisions, and scenarios cannot be ranked without a way to aggregate the nondominated set into a single number. 
A \emph{quality indicator}~\citep{zitzler2003performance,li2019quality} provides exactly such an aggregation. Following~\citep[Definition~2]{zitzler2003performance}, we let~$\Omega$ denote the collection of all sets of points in the objective space that are under consideration---in our setting, sets of nondominated points arising from the bi-objective problem---and, for an \emph{arity}~$m\in\N$, we call a function~$I\colon \Omega^m \to \R$
an \emph{$m$-ary quality indicator}. An indicator thus assigns a real value to each $m$-tuple of point sets. Since~$\R$ is totally ordered, these values can be compared and optimized.\footnote{Being scalar-valued is part of the formal definition of a quality indicator~\citep[Definition~2]{zitzler2003performance}. We nevertheless emphasize this property, as it is precisely what reduces the set-valued outcome of a bi-objective problem to a single number. Note, however, that an $m$-ary indicator with~$m\ge2$ is sometimes instead interpreted as assigning~$m$ relative quality values, one per compared set~\citep{li2019quality}.} 
In the multi-objective optimization literature, quality indicators are commonly unary or binary, corresponding to arities~$m=1$ and~$m=2$, respectively: a \emph{unary} indicator~$I\colon\Omega\to\R$ evaluates a single set in isolation, whereas a \emph{binary} indicator~$I\colon\Omega^2\to\R$ assesses the relative quality of an ordered pair of sets~\citep{zitzler2003performance,li2019quality}. Both kinds are widely employed---binary indicators, for instance, to guide indicator-based evolutionary search~\citep{zitzler2004}---and comprehensive surveys of the indicators proposed in the literature are given in~\cite{zitzler2003performance,riquelme2015performance,li2019quality}. Note that fixing the second argument of a binary indicator~$I$ to a reference set~$R\in\Omega$ yields the unary indicator~$I(\cdot,R)$.

Only unary indicators are directly usable in our setting: a unary indicator~$I$ induces a total preorder on~$\Omega$ via $A\preceq B :\Longleftrightarrow I(A)\le I(B)$, thereby providing a scalar objective for the master problem and a basis for comparing scenarios. Which unary indicator to use is a modeling choice, and several alternatives have been used within the framework we follow (see below).
We use the \emph{hypervolume indicator}~\citep{zitzler1998multiobjective,zitzler2003performance}, one of the most widely used quality indicators in multi-objective optimization. It measures the hypervolume---an area in our bi-objective setting---of the region in the objective space weakly dominated\footnote{A point~$b$ is weakly dominated by a point~$a$ if each entry of~$a$ is better than or equal to the corresponding entry of~$b$; it is weakly dominated by a point set if it is weakly dominated by at least one point in that set.} by the point set being evaluated and weakly dominating a given reference point. The indicator enjoys strong theoretical properties and offers a natural geometric interpretation (see Section~\ref{sec:literature}).

Computing the hypervolume inside the master problem leads to a nonlinear formulation. We address this by a \emph{slab decomposition} of the objective space, which yields a tractable linear MIP that approximately solves the master problem. The resulting algorithm iteratively either adds a worst-case scenario or refines the slab decomposition until a user-specified optimality gap is reached. Section~\ref{sec:fund} introduces the formal problem setting and notation; Section~\ref{sec:alg} develops the slab decomposition model and algorithm; Section~\ref{sec:experiments} describes the application problems and presents the computational results; and Section~\ref{sec:concl} concludes the paper.

\paragraph{Contributions}
We propose a \CCG{}-based algorithm for bi-objective two-stage robust MIPs using the hypervolume indicator. Our main contributions are:
\begin{itemize}
	\item A column-and-constraint generation approach for bi-objective two-stage robust MIPs using the hypervolume indicator to evaluate the scenario-dependent nondominated sets in the second stage.
	\item A slab decomposition of the objective space that approximates the hypervolume and enables the master problem to be formulated as a linear MIP.
	\item An optimality gap certificate based on bounding the error of the slab decomposition, together with a strategy that adaptively decides between adding new scenarios and refining the decomposition.
    \item A computational evaluation on two-stage robust selection and assignment problems, showing that the algorithm typically adds only a few scenarios to the master problem and developing a simple instance-dependent rule for choosing the initial slab width, while also analyzing the effects of instance size, scenario-set size, target gap, reference point, and slab orientation.
\end{itemize}

\subsection{Related literature}\label{sec:literature}

Two-stage robust optimization was introduced in~\cite{bib27ben2004adjustable} as an extension of classical robust optimization~\citep{ben2009robust} that allows recourse decisions to be made after some or all uncertain parameters have been revealed. Problems of this kind are formulated as min-max-min problems, where a first-stage decision is chosen to minimize the worst second-stage cost over all scenarios. In the case of finite uncertainty sets, which have been studied extensively (see, e.g.,~\cite{bib29bienstock2008computing,toenissen11,RODRIGUES2021499}), the min-max-min problem can be reformulated as a single-stage model using an epigraph formulation~\citep{goerigk2024buch}. This model, however, quickly becomes impractical as the number of scenarios or the problem size grows.

To overcome the scalability limitations of such reformulations, \CCG{} methods iteratively add new scenarios and corresponding second-stage variables to the master problem~\citep{zeng,toenissen11,RODRIGUES2021499}. A generic \CCG{} method proceeds as follows:
\begin{enumerate}
	\item Initialize a typically small subset of scenarios, often consisting of a single scenario.
	\item Solve the master problem using only the subset of scenarios and obtain a first-stage decision.
	\item Find a worst scenario among the remaining scenarios by evaluating the objective values of the second-stage problems given the fixed first-stage decision.
	\item Unless a termination condition is met, add a worst scenario to the subset and go to Step~2.
\end{enumerate}

\noindent
Several improvements have been proposed for settings where the second-stage problems are particularly hard to solve~\citep{toenissen11,RODRIGUES2021499,goerigk2025asbp}. {\NoHyper\citeauthor{TSANG202392}\endNoHyper}~\cite{TSANG202392} consider a variant in which the master problem is solved only up to a certain optimality gap, and a new scenario is added only if its objective improvement is sufficiently large. Otherwise, the to-be-achieved gap of the master problem is reduced. Our method follows an analogous adaptive principle, with the inaccuracy arising from the slab-decomposition approximation rather than an inexact master solve. Accordingly, we refine the decomposition by reducing the slab width instead of tightening the master-problem gap.

The systematic development of robust approaches to multi-objective optimization is comparatively recent. Early contributions addressed uncertainty in multi-objective problems through sequential or scenario-evaluation approaches rather than integrated robust formulations~\citep{mattson2003concept,tock2015decision,lemos2018new,gabrielli2019robust,groetzner2022multiobjective}. Theoretical foundations for robust multi-objective optimization---in particular, multi-objective generalizations of min-max problems---were established in~\cite{kuroiwa2012robust,ehrgott2014minmax,bokrantz2017necessary}. Good overviews are provided in~\cite{kuhn2016bi,ide2016robustness,wiecek2016robust}. 

A natural paradigm for two-stage robust multi-objective optimization is to first evaluate each nondominated set by a scalar-valued quality indicator and then apply standard two-stage robust methods---an approach termed \emph{first-scalarize-then-robustify} in~\cite{halser2025multicriteria}.  
In the terminology introduced in Section~\ref{sec:intro}, the set evaluated by the indicator is here the nondominated set for a given scenario and a fixed first-stage decision.
Several quality indicators have been applied within this framework. The authors of~\cite{halser2025multicriteria} use a weighted-sum indicator and a constraint indicator, both of which are unary. The first scalarizes the vector-valued objective by taking a weighted sum of the individual objectives, and the second optimizes only one objective while constraining the others to certain bounds. The latter was applied to a building energy supply design problem in~\cite{halser2024multicriteria}. In~\cite{goerigk2021flexible}, the authors use the binary $\varepsilon$-indicator of~\cite{zitzler2003performance}, which measures the distance between two nondominated sets---here, between the set for a fixed first-stage decision and the set obtained when the first-stage decision is not fixed. Since its second argument thus serves as a fixed reference set, the binary indicator is effectively employed as a unary one. They apply this approach to an energy system design selection problem. 

The hypervolume indicator~\citep{zitzler1998multiobjective,zitzler2003performance} is a unary quality indicator that measures the hypervolume of the region in the objective space weakly dominated by the point set being evaluated and weakly dominating a given reference point. 
It is one of the most widely studied quality indicators in multi-objective optimization~\citep{li2019quality}. Among the commonly used unary quality indicators, the hypervolume indicator, for a suitably chosen reference point, is currently the only one known to be fully consistent with the \emph{better} relation on solution sets (a property called \emph{strict Pareto compliance}): whenever a set~$A$ is better than a set~$B$---i.e., each element of~$B$ is weakly dominated by some element of~$A$ and at least one element of~$A$ is not weakly dominated by any element of~$B$---the hypervolume of~$A$ is strictly larger than that of~$B$~\citep{li2019quality}. Furthermore, empirical and theoretical observations suggest that the hypervolume indicator captures aspects of both coverage and uniformity~\citep{li2019quality}. Here, given a distance measure in the objective space (usually defined by a metric), the \emph{coverage error} of a set of points is the smallest value~$\varepsilon\geq0$ such that, for any nondominated point~$y$, there exists a point~$y'$ in the set within a distance of at most~$\varepsilon$ from~$y$. The set is \emph{$\delta$-uniform} for~$\delta\geq0$ if the distance between any two distinct points in the set is at least~$\delta$. Increasing the hypervolume typically reduces the coverage error, and, for a fixed number of points, hypervolume improvements are usually associated with a more uniform distribution, both properties that are considered desirable~\citep{li2019quality}. The reference point is often chosen as the so-called \emph{Nadir point} (i.e., the point obtained by taking the worst value of each objective among all nondominated points),
or a slight relaxation thereof (e.g., $1.1$ times the Nadir point in the case of a minimization problem). Since determining the Nadir point exactly is difficult in our setting, we discuss strategies for selecting an appropriate reference point in detail in Appendix~\ref{sec-app:refpoint}.

\section{Fundamentals}\label{sec:fund}

In this section, we define some notation and basic concepts used throughout the paper. To this end, we let~$S$ denote a finite set of scenarios. The set of feasible vectors of values for the first-stage decision variables~$x$ is denoted by~$\X$. The second-stage variables are called~$y^s$ for each scenario~$s\in S$, and their feasible set~$\Y^s(x)$ depends on the scenario~$s$ and the first-stage decision~$x$.
We assume~$\X$ to be nonempty and that~$\X$ and~$\Y^s(x)$ can be expressed as feasible sets of mixed-integer programs~(MIPs), i.e., using finitely many linear equality and inequality constraints and integrality constraints on some or all of the variables, where the linear constraints defining~$\Y^s(x)$ are jointly linear in~$x$ and~$y$. Further, for each~$s \in S$, $f_1^s(x, y^s)$ and $f_2^s(x, y^s)$ describe two nonnegative, affine, and scalar-valued objective functions (in~$x$ and~$y^s$). We use~$f^s(x,y^s)$ as a shorthand for the vector~$(f_1^s(x, y^s), f_2^s(x, y^s))\in\R^2$.
The restriction to nonnegative objective functions is without loss of generality for objective functions that are bounded below over all feasible first- and second-stage decisions: we can translate the value of an objective function and the corresponding entry of the reference point by the same constant without affecting the resulting hypervolume. Moreover, we restrict our presentation to minimization problems. Maximization problems whose objective functions are bounded above over all feasible first- and second-stage decisions can be transformed into this setting without changing the hypervolume values.\footnote{For each objective~$\ell$, choose a constant~$c_\ell$ that is an upper bound on~$f_\ell^s(x,y)$ for every~$s\in S$ and every feasible pair~$(x,y)$, and that also satisfies $c_\ell\ge\refpt_\ell$. Replacing~$f_\ell^s$ by~$c_\ell-f_\ell^s$ and the corresponding reference-point coordinate~$\refpt_\ell$ by~$c_\ell-\refpt_\ell$ yields a nonnegative minimization objective and a nonnegative reference-point coordinate. Applying this transformation to both objectives maps the dominated regions of the maximization problem to those of the transformed minimization problem and preserves their areas. Hence, all scenario-wise hypervolume values and the optimal first-stage decisions remain unchanged.}

Whenever a scenario~$s\in S$ is fixed and clear from the context, we may omit the superscript~$s$ from~$y^s$, writing~$y$.
If, instead, a used method relies on multiple vectors of second-stage variables for a scenario~$s\in S$, we additionally equip these vectors with a superscript~$i$ and write $\left(y^{s,i}\in\Y^s(x)\right)_{i\in I}$ for some index set~$I$. In general, we denote the set of positive integers by~$\N$ and, for $n\in\N$, the set~$\{1,\dots,n\}$ by~$[n]$.

\subsection{Bi-objective optimization}\label{sec:fund-biobj}

We now introduce some notation and definitions regarding bi-objective optimization. To this end, for each scenario~$s\in S$, we write
\begin{align*}
    &\XY \coloneqq \{(x, y): x\in\X, y\in\Y^s(x) \}\\
	&f^s(\X, \Y^s) \coloneqq \{(f_1^s(x, y), f_2^s(x, y)): (x,y) \in \XY \}.%
\end{align*}
Here, each pair~$(x,y)\in\XY$ is called a \emph{(feasible) solution} for scenario~$s$, and $\XY$ is the set of all such solutions in the decision space. Moreover, $f^s(\X, \Y^s)$ is the corresponding set of 
\emph{points}
in the objective space~$\R^2$. For a subset~$\X'\subseteq\X$ of first-stage decisions and a scenario~$s\in S$, we write~$\X'\times\Y^s \coloneqq \{(x,y): x\in\X', y\in\Y^s(x)\}$ and~$f^s(\X',\Y^s) \coloneqq \{f^s(x,y): (x,y)\in\X'\times\Y^s\}$ for the corresponding restricted sets; in particular, $f^s(\{x\},\Y^s)$ denotes the set of points in the objective space attainable for a fixed first-stage decision~$x$. We say that a solution~$(x,y)\in\XY$ corresponds to a point~$p\in f^s(\X,\Y^s)$ if~$p=f^s(x,y)$ and vice versa. If~$x$ is fixed and clear from the context, we also say that a second-stage decision~$y\in\Y^s(x)$ corresponds to the point~$p= f^s(x,y)$ and vice versa. Note that, in both cases, the relation is not 1:1, meaning that there can be multiple distinct solutions corresponding to one point in the objective space.

We now define the concept of (weak) dominance for minimization problems (see, e.g., \cite{ehrgott2005multicriteria}):

\begin{definition}
    We say that~$\point\in\mathbb{R}^2$ \emph{weakly dominates}~$\point'\in\mathbb{R}^2$ (or, equivalently, that~$\point'$ \emph{is weakly dominated} by~$\point$) if $\point_\ell\le\point'_\ell$ for $\ell=1,2$. If additionally $\point \ne \point'$, we say that~$\point$ \emph{dominates}~$\point'$ (or, equivalently, that~$\point'$ \emph{is dominated} by~$\point$).
	If, for a scenario~$s\in S$, a subset~$\X'\subseteq\X$ of first-stage decisions, and a solution~$(x,y)\in\X'\times\Y^s$, the point~$\point=f^s(x,y)\in f^s(\X',\Y^s)$ is not dominated by any other point~$\point'\in f^s(\X',\Y^s)$, we call the point~$\point$ \emph{nondominated for~$s$ and~$\X'$} and the solution~$(x, y)$ \emph{efficient for~$s$ and~$\X'$}. The set~$\PND^s(\X')$ of all nondominated points for~$s$ and~$\X'$ is called the \emph{nondominated set for~$s$ and~$\X'$}.

	In the case~$\X'=\X$, we omit the argument~$\X'$ and simply speak of \emph{nondominated points} and \emph{efficient solutions for~$s$}, and we write~$\PND^s\coloneqq\PND^s(\X)$ for the \emph{nondominated set for~$s$}. 
    In the case~$\X'=\{x\}$ for a fixed first-stage decision~$x\in\X$, we replace the singleton argument~$\{x\}$ by~$x$ in our terminology and notation: we speak of \emph{nondominated points} and \emph{efficient solutions for~$s$ and~$x$}, and we write~$\PND^s(x)\coloneqq\PND^s(\{x\})$ for the \emph{nondominated set for~$s$ and~$x$}.

\end{definition}

\noindent
Besides dominance, we also use the lexicographic order on the objective space:

\begin{definition}\label{def:lex}%
	We say that~$\point\in\R^2$ is \emph{lexicographically smaller} than~$\point'\in\R^2$, written~$\point <_{\lex} \point'$, if either $\point_1 < \point'_1$, or $\point_1 = \point'_1$ and $\point_2 < \point'_2$. We write $\point \le_{\lex} \point'$ if $\point <_{\lex} \point'$ or $\point = \point'$, and we call a finite sequence $\point^1,\dots,\point^k$ of points in~$\R^2$ \emph{lexicographically non-decreasing} if $\point^i \le_{\lex} \point^{i+1}$ for all~$i\in[k-1]$.
	Let $s\in S$ and $\mathcal{F}\subseteq \XY$. A solution $(x^*,y^*) \in \mathcal{F}$ is a \emph{lexicographic minimum} of $(f^s_1, f^s_2)$ over $\mathcal{F}$ if
	\begin{equation*}
		\left(f_1^s(x^*,y^*),f_2^s(x^*,y^*)\right) \le_{\lex} \left(f_1^s(x,y),f_2^s(x,y)\right) \quad \text{for all } (x,y) \in \mathcal{F}.
	\end{equation*}
	We denote the problem of finding a lexicographic minimum of~$(f_1^s,f_2^s)$ over $\mathcal{F}$ by $\lex\;\min_{(x,y) \in \mathcal{F}} \left(f^s_1(x,y),\, f^s_2(x,y)\right)$.
\end{definition}

\noindent
Note that the order of the two objectives matters here: a lexicographic minimum of~$(f^s_1, f^s_2)$ first minimizes~$f^s_1$ and, among all solutions attaining this minimum, then minimizes~$f^s_2$. Accordingly, we write $\lex\min_{(x,y)\in\mathcal{F}}\left(f^s_2(x,y),\, f^s_1(x,y)\right)$ when the roles of the two objectives are interchanged, in which case the corresponding points are compared with respect to the lexicographic order on their reversed components.

Throughout this paper, we impose the following assumption on the nondominated sets, known as \emph{external stability} and commonly used in multi-objective optimization (see, e.g., \cite{ehrgott2005multicriteria}).

\begin{assumption}\label{ass:external-stability}
For every scenario~$s\in S$ and every first-stage decision~$x\in\X$, every point in~$f^s(\X,\Y^s)$ is weakly dominated by a nondominated point for~$s$, and every point in~$f^s(\{x\},\Y^s)$ is weakly dominated by a nondominated point for~$s$ and~$x$.
\end{assumption}

Assumption~\ref{ass:external-stability} is satisfied, in particular, whenever the sets~$\XY$ of feasible solutions are finite, as is the case for both applications considered in Section~\ref{sec:applications}. Moreover, it is automatically satisfied whenever the constraints describing~$\X$ and~$\Y^s(x)$ and the objective functions are given by rational data: for any point~$p\in f^s(\{x\},\Y^s)$, the problem $\lex\min_{y\in\Y^s(x)}\left(f_1^s(x,y), f_2^s(x,y)\right)$ with the additional constraints $f_\ell^s(x,y)\le p_\ell$ for~$\ell=1,2$ has an optimal solution~\citep{meyer1974existence}, and its point is nondominated for~$s$ and~$x$ and weakly dominates~$p$ (analogously for~$f^s(\X,\Y^s)$). 

\begin{algorithm}
	\caption{(Lexicographic) Epsilon-Constraint Method}\label{alg:eps-constraint}
	\begin{algorithmic}[1]
		\Statex \textbf{Input:} A step size~$\delta > 0$, a scenario~$s$, and a subset~$\X'\subseteq\X$ of first-stage decisions such that~$\X'=\X$ or that~$\X'$ is a finite set
		\Statex \textbf{Output:} A subset~$P \subseteq \PND^s(\X')$ of the nondominated set for~$s$ and~$\X'$, and a set~$\PSOL\subseteq \X'\times\Y^s$ of one corresponding efficient solution for each~$p\in P$
		\State $P \gets \varnothing$
		\State $\PSOL \gets \varnothing$
		\State Solve the problem
		\begin{equation}\label{alg:epsconstr-leftanchor}
			\begin{aligned}
				\underset{(x,y)\in\X'\times\Y^s}{\lex\;\min} & \quad \left(f_1^s(x,y),f_2^s(x,y)\right)
			\end{aligned}
		\end{equation}
        and call an optimal solution~$(\hat x,\hat y)$ if one exists.\label{alg-line:epsconstr-leftanchor}
		\If{Problem~\eqref{alg:epsconstr-leftanchor} has an optimal solution}
			\State $\point_1^{s,\LB} \gets f_1^s(\hat x,\hat y)$
			\State $\point_1^{s,\UB} \gets \max_{(x,y)\in\X'\times\Y^s} f_1^s(x,y)$ (set to $+\infty$ if unbounded)
			\State $\varepsilon \gets \point_1^{s,\UB}$
			\While{$\varepsilon \ge \point_1^{s,\LB}$}
            \State Set $\mathcal{F}\gets\{(x,y)\in\X'\times\Y^s:f_1^s(x,y)\le \varepsilon\}$, solve the problem
			\begin{equation}\label{alg:epsconstr-lexmin}
				\begin{aligned}
					\underset{(x,y)\in\mathcal F}{\lex\;\min} & \quad \left(f_2^s(x,y), f_1^s(x,y)\right)
				\end{aligned}
			\end{equation}
			\hspace{\algorithmicindent}\hspace{\algorithmicindent}and call an optimal solution~$(x^\ast, y^\ast)$.\label{alg-line:epsconstr-lexmin}
			\State $p \gets (f_1^s(x^\ast,y^\ast),f_2^s(x^\ast,y^\ast))$
			\State $P \gets P \cup \{p\}$
			\State $\PSOL \gets \PSOL\cup\{(x^\ast, y^\ast)\}$
			\State $\varepsilon \gets f_1^s(x^\ast,y^\ast) - \delta$
			\EndWhile
			\State $p \gets (f_1^s(\hat x,\hat y),f_2^s(\hat x,\hat y))$ % left anchor point
			\If{$p \notin P$}\label{alg-line:epsconstr-anchor-guard}
				\State $P \gets P \cup \{p\}$
				\State $\PSOL \gets \PSOL\cup\{(\hat x, \hat y)\}$
			\EndIf
		\EndIf\\
        \Return $P$ and $\PSOL$
	\end{algorithmic}
\end{algorithm}

One method to algorithmically construct a finite subset~$P$ of the nondominated set for a scenario~$s\in S$ and a subset~$\X'\subseteq\X$ of first-stage decisions, together with a set of one corresponding efficient solution for each~$p\in P$, is the \emph{lexicographic epsilon-constraint method} \citep{HAMACHER2007336,prinz2025parallel,haimes1971bicriterion,chankong2008multiobjective}. 
Algorithm~\ref{alg:eps-constraint} states the pseudocode for our setting. For brevity, we refer to the method simply as the \emph{epsilon-constraint method}.
To simplify the exposition, we restrict attention to the cases $\X'=\X$ and finite~$\X'\subseteq\X$. We further assume that every lexicographic optimization problem encountered in Algorithm~\ref{alg:eps-constraint} either is infeasible or admits an optimal solution. This mild attainment assumption is satisfied, in particular, whenever the relevant feasible sets are finite, as in both applications considered in Section~\ref{sec:applications}, and more generally under standard rational-data assumptions for the encountered
MIPs~\citep{meyer1974existence}.

The method bounds one objective (without loss of generality, the first objective~$f_1^s$) from above by a parameter~$\varepsilon$ and minimizes the second objective~$f_2^s$ subject to this bound.
To exclude weakly nondominated points, it then minimizes~$f_1^s$ in a second step, with~$f_2^s$ fixed at the optimal value obtained in the first step.
Let~$(x^\ast,y^\ast)$ denote the resulting solution.
The bound is then tightened to~$\varepsilon \leftarrow f_1^s(x^\ast,y^\ast) - \delta$ for a suitable~$\delta>0$, and the procedure is repeated.
The search interval is given by the smallest value that~$f_1^s$ attains over~$\X'\times\Y^s$ and by the supremum~$\point_1^{s,\UB}$ of~$f_1^s$ over this set: the upper end initializes~$\varepsilon$, and the search terminates once~$\varepsilon$ falls below the lower end, since no feasible solution satisfying the bound can exist beyond that point.

By the assumption on the existence of optimal solutions stated above, Problem~\eqref{alg:epsconstr-leftanchor} (see Line~\ref{alg-line:epsconstr-leftanchor} of Algorithm~\ref{alg:eps-constraint}) has an optimal solution whenever $\X'\times\Y^s\ne\emptyset$; otherwise, the problem is infeasible and the method returns empty sets of nondominated points and efficient solutions. Moreover, Problem~\eqref{alg:epsconstr-lexmin} (see Line~\ref{alg-line:epsconstr-lexmin}) has an optimal solution in every iteration of the while loop. Indeed, the loop condition implies $\varepsilon\ge\point_1^{s,\LB}=f_1^s(\hat x,\hat y)$, so that $(\hat x,\hat y)$ is feasible for Problem~\eqref{alg:epsconstr-lexmin} in every iteration. By the above assumptions, this problem therefore admits an optimal solution in every iteration.

Note that~$f_1^s$ need not be bounded from above on $\X'\times\Y^s$, in which case~$\point_1^{s,\UB}=+\infty$. In this case, we interpret the constraint $f_1^s(x,y)\le\varepsilon$ as vacuous in the first iteration of the while loop, i.e., Problem~\eqref{alg:epsconstr-lexmin} is solved over the entire set~$\X'\times\Y^s$. The point returned in this iteration nevertheless has finite coordinates, since it corresponds to an optimal solution of Problem~\eqref{alg:epsconstr-lexmin}. Consequently, all subsequent values of~$\varepsilon$ are finite.
Unless stated otherwise, we apply the method with~$\X'=\{x\}$ for a fixed first-stage decision~$x\in\X$, in which case it constructs (a subset of) the nondominated set~$\PND^s(x)$ for~$s$ and~$x$; this is the case relevant for Lemma~\ref{lemma:hv-eps-constr-error} below and for the evaluation step of our algorithm in Section~\ref{sec:outline}. The only place where we use a different subset is Section~\ref{sec:comp}, where the method is applied with~$\X'=\X$ in order to construct (a subset of) the nondominated set~$\PND^s$ for~$s$ over all first-stage decisions.
Note that, due to the strictly positive step size~$\delta$, the algorithm may not return all nondominated points, but only a subset of them.

To evaluate the quality of a nondominated set, we use the \emph{hypervolume indicator}~\citep{zitzler2003performance}. More generally, for any set~$P$ of points, it measures the hypervolume---an area in our bi-objective setting---of the region in the objective space weakly dominated by~$P$ and weakly dominating a given reference point:

\begin{definition}\label{def:HV}
	Given a \emph{reference point}~$\refpt=(\refpt_1,\refpt_2)\in\mathbb{R}_{\geq 0}^2$ and a set~$P\subseteq\mathbb{R}_{\geq 0}^2$ of points in the objective space, the \emph{(dominated) hypervolume of~$P$ (with respect to~$\refpt$)}, denoted by~$\HV(P)$, is defined as the area of the \emph{dominated region}
	\begin{align*}
		R(P)\coloneqq\bigcup_{\point'\in P}\left\{\point\in\mathbb{R}^2: \point'_\ell\le\point_\ell\le\refpt_\ell \text{ for } \ell=1,2 \right\}\subseteq\mathbb{R}^2.
	\end{align*}
	For a finite set~$P=\{\point^1,\dots,\point^n\}$, the region~$R(P)$ is a finite union of axis-parallel rectangles, and~$\HV(P)$ is the area of this union. For a general set~$P$, the region~$R(P)$ is still measurable: writing~$\mu_P(\tau)\coloneqq\inf\{\point'_2: \point'\in P,\, \point'_1\le\tau\}$ (with $\inf\varnothing\coloneqq+\infty$) for the non-increasing extended-real-valued function whose graph bounds~$R(P)$ from below, the area of~$R(P)$ equals $\int_0^{\refpt_1}\max\{\refpt_2-\mu_P(\tau),\,0\}\,\mathrm{d}\tau$.

	Let~$\PND^s(x)$ be the nondominated set for a scenario~$s\in S$ and a first-stage decision~$x\in\X$. We call $\HV^s(x) \coloneqq \HV(\PND^s(x))$ the \emph{(dominated) hypervolume for~$s$ and~$x$}.
\end{definition}

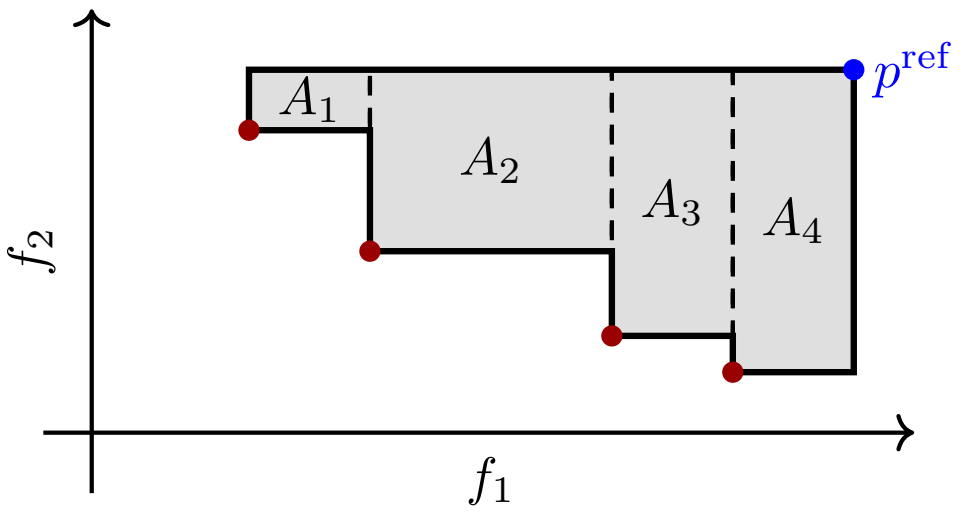
\begin{figure}
	\centering
	\begin{tikzpicture}[scale=0.8, thick]

% points
\coordinate (P1) at (1,2.5);
\coordinate (P2) at (2,1.5);
\coordinate (P3) at (4,0.8);
\coordinate (P4) at (5,0.5);

% reference point
\coordinate (R)  at (6,3);

% rectangles + labels
\path[fill=mygray, opacity=0.3] (P4) rectangle (R);
\path[fill=mygray, opacity=0.3] (P3) rectangle (R -| P4);
\path[fill=mygray, opacity=0.3] (P2) rectangle (R -| P3);
\path[fill=mygray, opacity=0.3] (P1) rectangle (R -| P2);
\draw[black,dashed] (P4) rectangle (R) node[pos=.5] {$A_4$};
\draw[black,dashed] (P3) rectangle (R -| P4) node[pos=.5] {$A_3$};
\draw[black,dashed] (P2) rectangle (R -| P3) node[pos=.5] {$A_2$};
\draw[black,dashed] (P1) rectangle (R -| P2) node[pos=.5] {$A_1$};

\draw[black, very thick]
(R) -- (R -| P1) -- (P1) -| (P2) -| (P3) -| (P4) -- (P4 -| R) -- (R);

% points
\fill[pointred] (P1) circle (2.5pt);
\fill[pointred] (P2) circle (2.5pt);
\fill[pointred] (P3) circle (2.5pt);
\fill[pointred] (P4) circle (2.5pt);

\filldraw[pointblue] (R) circle (2pt);
 \node[pointblue, right] at (R) {$\refpt$};

% Axes
\draw[->, thick] (-0.7,0) -- (6.5,0);      % x-axis
\draw[->, thick] (-0.3,-0.5) -- (-0.3,3.5); % y-axis
% Axes labels
\node at (-0.8,1.5) [rotate=90] {$f_2$};
\node at (3.0,-0.4) {$f_1$};

\end{tikzpicture}
    \caption{Illustration of the hypervolume as a union of rectangles. For the nondominated point set~$P$ shown in red and the reference point~$\refpt$, the dominated region~$R(P)$ is outlined in black, shaded in gray, and partitioned into the rectangles $A_1,\dots,A_n$.}
	\label{fig:hv-partition}
\end{figure}

Note that a point~$\point\in P$ that does not weakly dominate the reference point
(i.e., for which $\point_\ell>\refpt_\ell$ for some~$\ell$) contributes an empty set
to the union in Definition~\ref{def:HV} and therefore does not affect~$R(P)$ or
$\HV(P)$. Moreover, the hypervolume is monotone with respect to set inclusion:
if~$P'\subseteq P$, then~$R(P')\subseteq R(P)$ and, consequently,
$\HV(P')\le\HV(P)$. In particular, any subset of the nondominated set, such as
one returned by the epsilon-constraint method, has hypervolume at most that of
the complete nondominated set.

Let~$\refpt$ be a reference point and~$P\subseteq\R_{\ge0}^2$ be a set of~$n$ nondominated points, each of which weakly dominates~$\refpt$. Assuming the points in~$P$ are ordered by non-decreasing first component (and, thus, non-increasing second component), the hypervolume~$\HV(P)$ can be computed in $\mathcal{O}(n)$ time~\citep{fonseca2006improved,beume2009complexity} by partitioning the dominated region into $n$ rectangles $A_1, \dots, A_n$ as shown in Figure~\ref{fig:hv-partition}:
\begin{align} \label{eq:hv-calculation}
	\HV(P)=\sum_{i=1}^n \operatorname{area}(A_i)=\sum_{i=1}^n (\point^{i+1}_1-\point^i_1)\cdot(\refpt_2-\point^{i}_2), \qquad \textrm{with } \point^{n+1}_1\coloneqq\refpt_1.
\end{align}

When the set~$\PND^s(x)$ of nondominated points is not explicitly available, the hypervolume for~$s$ and~$x$ can still be modeled exactly, namely by a quadratic program with linear constraints in which finitely many second-stage decisions describe the nondominated points weakly dominating~$\refpt$. We state this model in Appendix~\ref{sec-app:hv-opt}. The fact that an exact model of the hypervolume is quadratic is the reason why we develop a linear approximation of it in Section~\ref{sec:model}.

Since the epsilon-constraint method may return only a (finite) subset~$P$
of the nondominated set, the hypervolume~$\HV(P)$ that results from using~\eqref{eq:hv-calculation} may be smaller than the hypervolume~$\HV^s(x)$ of the complete nondominated set. The following lemma bounds this approximation error.
\begin{lemma}\label{lemma:hv-eps-constr-error}
Let~$s\in S$, $x\in\X$, and let~$P\subseteq \PND^s(x)$ be a finite subset of the nondominated set for~$s$ and~$x$ returned by the epsilon-constraint method applied to~$s$ and~$\X'=\{x\}$ with step size~$\delta>0$.
If~$\PND^s(x)\neq\emptyset$, then
\begin{align*}
	\HV^s(x)-\HV(P)\;\le\;\delta\cdot\Bigl(\max_{\point\in P}\point_2-\min_{\point\in P}\point_2\Bigr).
\end{align*}
In all cases, $\HV^s(x)-\HV(P)\;\le\;\delta\cdot\refpt_2$.
Finally, if the nondominated set~$\PND^s(x)$ is finite and the step size~$\delta$ is at most the smallest difference $|q_1-q'_1|$ between the first components of two distinct nondominated points $q,q'\in\PND^s(x)$ (with this minimum taken to be~$+\infty$ if $|\PND^s(x)|\le 1$), then the epsilon-constraint method returns the complete nondominated set, i.e., $P=\PND^s(x)$ and
\begin{align*}
	\HV^s(x)-\HV(P)=0.
\end{align*}
This holds, in particular, whenever the distinct values attained by $f_1^s$ have a positive minimum distance and~$\delta$ is at most this distance---for instance, whenever~$f_1^s$ is integer-valued and $\delta\le1$.
\end{lemma}
\begin{proof}[Proof idea.]
If~$\PND^s(x)=\emptyset$, then also~$P=\emptyset$, since $P\subseteq\PND^s(x)$, and therefore $\HV^s(x)=\HV(P)=0$. Thus, all assertions applicable to this case follow immediately. Hence, suppose that~$\PND^s(x)\neq\emptyset$. Then the epsilon-constraint method returns at least one point, so $P\neq\emptyset$. Since~$P\subseteq\PND^s(x)$, we have~$R(P)\subseteq R(\PND^s(x))$, and the difference~$\HV^s(x)-\HV(P)$ is the area of
$R(\PND^s(x))\setminus R(P)$. Writing~$P=\{\point^1,\dots,\point^k\}$ with strictly increasing first components (and, thus, strictly decreasing second components), the construction of the epsilon-constraint method implies that every nondominated point missing from~$P$ lies strictly between two consecutive returned points~$\point^j$ and~$\point^{j+1}$ and, in addition, less than~$\delta$ to the left of~$\point^{j+1}$: otherwise, the iteration with bound~$\varepsilon=\point^{j+1}_1-\delta$ would have returned it (or a point with an even smaller second component) instead of~$\point^j$. Hence, the missed region is contained in the union of the~$k-1$ rectangles $\left(\point^{j+1}_1-\delta,\,\point^{j+1}_1\right)\times\left(\point^{j+1}_2,\,\point^j_2\right)$, $j\in[k-1]$, each of which has width~$\delta$. Summing their areas telescopes to $\delta\cdot(\point^1_2-\point^k_2)$, which proves the bound $\HV^s(x)-\HV(P) \le \delta\cdot \Bigl(\max_{\point\in P}\point_2-\min_{\point\in P}\point_2\Bigr)$. 
To obtain the bound~$\HV^s(x)-\HV(P)\le\delta\cdot\refpt_2$, intersect the rectangles with the box~$[0,\refpt_1]\times[0,\refpt_2]$. 
Their vertical projections are pairwise disjoint and contained in~$[0,\refpt_2]$, so their total height is at most~$\refpt_2$ and their total area is at most~$\delta\cdot\refpt_2$.

To prove that~$P=\PND^s(x)$ under the stated minimum-distance condition on the first components of the nondominated points, suppose first that $|\PND^s(x)|\ge2$, and let~$d>0$ denote the smallest difference between the first components of two distinct points in~$\PND^s(x)$. If a nondominated point~$q$ were missing from~$P$, the argument above would yield two consecutive returned points~$\point^j,\point^{j+1}$ such that $\point^j_1<q_1<\point^{j+1}_1$ and $\point^{j+1}_1-q_1<\delta$. Since~$q$ and~$\point^{j+1}$ are distinct nondominated points, the definition of~$d$ gives $d\le\point^{j+1}_1-q_1<\delta\le d$, a contradiction. Hence, $P=\PND^s(x)$; the case $|\PND^s(x)|\le1$ is immediate. Finally, suppose that the distinct values attained by~$f_1^s$ have positive minimum distance~$d$ and that~$\delta\le d$. The point returned in the first iteration of the while loop of the epsilon-constraint method minimizes~$f_2^s$ over all feasible second-stage decisions. Hence, no nondominated point can have a larger first component than this point. Since~$f_1^s$ is nonnegative and its distinct attainable values are separated by at least~$d$, the nondominated set~$\PND^s(x)$ is finite. Moreover, the first components of any two distinct nondominated points differ by at least~$d$. Thus, the preceding condition for $P=\PND^s(x)$ is satisfied. In particular, if~$f_1^s$ is integer-valued, distinct attainable values differ by at least~$1$, so $\delta\le1$ suffices.
\end{proof}

\subsection{Bi-objective two-stage robust optimization}\label{sec:fund-biobjrobust}

In terms of the notation introduced in
Section~\ref{sec:fund-biobj}, the bi-objective two-stage robust
mixed-integer programs considered in this paper can be represented as
\begin{align}\label{eq:bi-obj-model}
	\min_{x \in \X} \max_{s \in S} \min_{y^s\in\Y^s(x)} \bigl( f_1^s(x, y^s), f_2^s(x, y^s) \bigr).
\end{align}
Here, the inner minimization is understood in the Pareto sense and
therefore yields the nondominated set~$\PND^s(x)$ for the scenario~$s$
and the fixed first-stage decision~$x$.
Since smaller objective values weakly increase the hypervolume, minimizing
the two objectives corresponds to maximizing the hypervolume. We therefore evaluate the inner bi-objective problem for scenario~$s$ and first-stage decision~$x$ by the scalar value~$-\HV^s(x)$. This yields the problem
\[
    \min_{x\in\X}\max_{s\in S}\bigl(-\HV^s(x)\bigr).
\]
For every~$x\in\X$, we have $\max_{s\in S}(-\HV^s(x))=-\min_{s\in S}\HV^s(x)$. Consequently, the preceding problem is equivalent to
\begin{equation}\label{eq:hv-maxmin-problem}
    \max_{x\in\X}\min_{s\in S}\HV^s(x).
\end{equation}
We refer to Problem~\eqref{eq:hv-maxmin-problem} as the
\emph{overall problem}. 
Since~$\X\neq\emptyset$ by assumption and~$0\leq\HV^s(x)\leq\refpt_1\refpt_2$ for every~$s\in S$ and~$x\in\X$, Problem~\eqref{eq:hv-maxmin-problem} is feasible and has a finite, nonnegative optimal value.

\section{Algorithm}\label{sec:alg}

We now present our slab decomposition algorithm for solving Problem~\eqref{eq:hv-maxmin-problem}.
Throughout this section, we fix a user-defined reference point~$\refpt=(\refpt_1,\refpt_2)\in\R_{>0}^2$. We exclude reference points with a zero coordinate, since the nonnegativity of the objective functions would then imply that all hypervolume values are zero. Apart from this restriction, the validity of the method is independent of the specific choice of~$\refpt$. 
Appendix~\ref{sec-app:refpoint} discusses several strategies for
choosing the reference point.

\begin{figure}
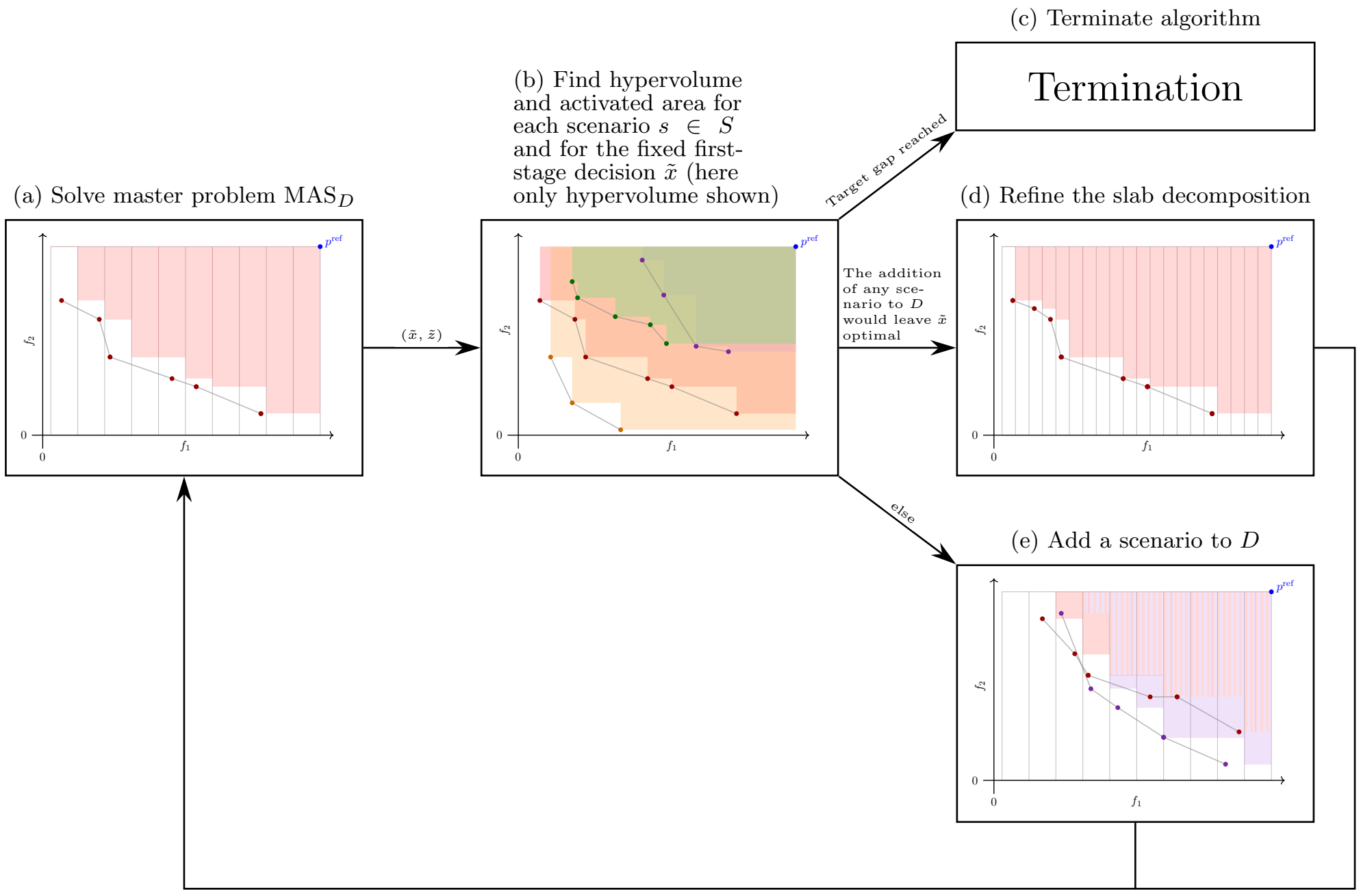

	\centering
	% The whole picture is scaled to \textwidth so it always fits the page width.
\resizebox{\textwidth}{!}{%
\begin{tikzpicture}[
		box/.style={draw, line width=0.7pt, inner sep=4pt},
		term/.style={draw, line width=0.7pt, minimum height=1.2cm, font=\Large},
		arr/.style={-{Stealth[length=10pt,width=6pt]}, line width=0.7pt},
		toplabel/.style={font=\scriptsize, anchor=south, inner sep=4pt},
		arrlabel/.style={font=\tiny, fill=white, inner sep=2pt},
	]
	% each sub-plot is its own tikzpicture; scale it down to fit as a node
	\newcommand{\subplot}[1]{\scalebox{0.46}{\input{figures/alg-overview-#1}}}

	%            T
	%          /
	%  A -- B -- D
	%          \
	%		     C

	\node[box]                       (A) {\scalebox{0.46}{\input{figures/alg-overview-A}}};
	\node[box, right=1.6cm of A]     (B) {\scalebox{0.46}{\input{figures/alg-overview-B}}};
	\node[box, right=1.6cm of B]     (D) {\scalebox{0.46}{\input{figures/alg-overview-D}}};
	% capture D's width so the Terminate box can be made exactly as wide
	\path let \p1=($(D.east)-(D.west)$) in \pgfextra{\xdef\boxwidth{\x1}};
	\node[term, minimum width=\boxwidth, above=1.2cm of D] (T) {Termination};
	\node[box, below=1.2cm of D]     (C) {\scalebox{0.46}{\input{figures/alg-overview-C}}};

	% one-line labels directly above the top border of each box
	\node[toplabel] at (A.north) {(a) Solve master problem~$\MAS_D$};
	\node[toplabel, text width=4cm, align=left] at (B.north) {(b) Find hypervolume and activated area for each scenario $s\in S$ and for the fixed first-stage decision $\tilde x$ (here only hypervolume shown)};
    \node[toplabel] at (T.north) {(c) Terminate algorithm};
	\node[toplabel] at (C.north) {(e) Add a scenario to~$D$};
	\node[toplabel] at (D.north) {(d) Refine the slab decomposition};

	% arrows (with labels on them)
	\draw[arr] (A) -- (B) node[arrlabel, midway, above] {$(\tilde x, \tilde z)$};
	\draw[arr] (B.south east) -- (C.north west) node[arrlabel, sloped, above, pos=0.5] {else};
	\draw[arr] (B.north east) -- (T.south west) node[arrlabel, sloped, midway, above, fill=none, yshift=2pt] {Target gap reached};
	\draw[arr] (B) -- (D)
		node[arrlabel, pos=0.5, text width=1.5cm, above, fill=none]
		{
			The addition of any scenario to~$D$ would leave~$\tilde x$ optimal
		};

	% loop back: after refining the decomposition (D) or adding a scenario (C),
	% the master problem is solved again
	\coordinate (rety) at ($(C.south)+(0,-0.9)$);
	\draw[arr] (C.south) -- ++(0,-0.9) -| (A.south);
	\draw[line width=0.7pt] (D.east) -- ++(0.55,0) |- (rety);
\end{tikzpicture}%
}
	\caption{Schematic overview of the slab decomposition algorithm (Algorithm~\ref{alg:slabdecomp}). Each box corresponds to a step of the algorithm, and the small plots sketch the associated situation in the objective space. Starting from an initial subset~$D\subseteq S$, shown here as a single scenario, the master problem~$\MAS_D$ is solved, yielding a master solution~$(\tilde x,\tilde z)$. For the first-stage decision~$\tilde x$, the lower bounds~$\HV_{\text{LB}}^s(\tilde x)$ and~$\AA_{\text{LB}}^s(\tilde x)$ are computed for every scenario~$s\in S$, and the optimality gap is evaluated as in Equation~\eqref{eq:gap}. The algorithm then branches: if the computed gap does not exceed the target gap~$\TG$, it terminates and returns~$\tilde x$; otherwise, if~$\AA_{\text{LB}}^s(\tilde x)\ge\tilde z$ holds for all scenarios~$s\in S\setminus D$, the slab decomposition is refined by reducing the slab width (by a factor~$\alpha$); else, a scenario not yet contained in~$D$ with the smallest hypervolume lower bound among those satisfying~$\AA_{\text{LB}}^s(\tilde x)<\tilde z$ is added to~$D$. In the latter two cases, the master problem is re-solved and the procedure repeats.}
	\label{fig:alg-overview-all}
\end{figure}

We first give a high-level description of the algorithm, illustrated schematically in Figure~\ref{fig:alg-overview-all}; the letters in parentheses below refer to the corresponding boxes in the figure.
The algorithm follows a column-and-constraint generation (\CCG) scheme. It maintains a subset~$D\subseteq S$ of scenarios and, in each iteration, solves a linear MIP, called the \emph{master problem}~$\MAS_D$, over this subset to obtain a candidate first-stage decision~$\tilde x$~(a). The master problem approximates the hypervolume by means of a \emph{slab decomposition} of the objective space, i.e., a subdivision into thin rectangular slabs of width~$w$ whose activated portions yield a lower bound on the hypervolume. Because~$\MAS_D$ only approximates the hypervolume and only accounts for the scenarios in~$D$, its solution generally does not yet solve the overall problem~\eqref{eq:hv-maxmin-problem}.
The algorithm therefore evaluates the candidate~$\tilde x$ over all scenarios~$s\in S$, derives a lower and an upper bound on the optimal objective value of the overall problem, and computes an optimality gap from these bounds~(b).
If this gap does not exceed the user-specified target gap, the algorithm terminates and returns~$\tilde x$~(c). Otherwise, it improves the approximation in one of two ways: either it adds a scenario not yet contained in~$D$ that has the
smallest hypervolume lower bound among the eligible scenarios~(e), or---if adding any of the remaining scenarios would leave the optimal objective value of the master problem unchanged and~$\tilde x$ optimal---it refines the slab decomposition by reducing the slab width~(d).

The remainder of this section is organized as follows.
Section~\ref{sec:outline} outlines the algorithm and presents its pseudocode.
Section~\ref{sec:model} develops the slab decomposition and the master problem~$\MAS_D$.
Section~\ref{sec:error-and-gap} derives a bound on the discretization error, the resulting optimality
gap, and the slab-refinement rule.
Finally, Section~\ref{sec:correctness} establishes the correctness and finite termination of the algorithm.

\subsection{Algorithm outline}\label{sec:outline}

\begin{algorithm}
	\caption{Slab decomposition algorithm}\label{alg:slabdecomp}
	\begin{algorithmic}[1]
        \Statex \textbf{Input:} Scenario set~$S$, target gap~$\TG\in[0,1)$, initial slab width~$w\in(0,\refpt_1]$, safety factor~$\sigma\in(0,1)$, lower bound~$\alpha_{\min}\in(0,\sigma)$ on the refinement factor, step size~$\delta>0$ for the epsilon-constraint method.
		\Statex \textbf{Output:} Upon termination, either a first-stage decision with optimality gap at most~$\TG$ for Problem~\eqref{eq:hv-maxmin-problem} or the information that the optimal objective value of Problem~\eqref{eq:hv-maxmin-problem} is zero. Finite termination is guaranteed under Condition~\eqref{eq:termination-condition}.
		\State $D \gets $ \Init().
		\State Initialize $\EMmax$ as in Equation~\eqref{eq:EMmax}.
        \If{$\EMmax = 0$}
        \State \Return ``Optimal objective value of Problem~\eqref{eq:hv-maxmin-problem} is zero.'' \label{algline:emmax-zero}
        \EndIf
		\State $\text{Terminate} \gets \text{False}$.
		\While{$\text{Terminate} = \text{False}$}\label{alg-line:general-while-loop}
		\State Try to optimally solve~$\MAS_D$. \label{algline:solveMAS}
		\If{$\MAS_D$ is infeasible}
        \State \Return ``Optimal objective value of Problem~\eqref{eq:hv-maxmin-problem} is zero.'' \label{algline:mas-infeasible}
		\Else
		\State Let $(\tilde x, \tilde z)$ denote the master solution of $\MAS_D$.
		\EndIf
		\State \parbox[t]{\dimexpr\linewidth-\algmargin\relax}{Find $\AA_{\text{LB}}^s(\tilde x)$ and $\HV_{\text{LB}}^s(\tilde x)$ for each scenario $s\in S$, using Algorithm~\ref{alg:eps-constraint} with step size~$\delta$.}   \label{algline:findHV}
		\State Calculate gap as in Equation~\eqref{eq:gap}. \label{algline:gap}
		\If{$\text{gap} \le \TG$} \label{algline:generalterminationcrit}
		\State $\text{Terminate} \gets \text{True}$.
		\Else
		\If{$\AA_{\text{LB}}^s(\tilde x) \ge \tilde z$ for each scenario $s\in S\setminus D$}\label{algline:branchcondition}
		\State \parbox[t]{\dimexpr\linewidth-\algmargin\relax}{Compute the refinement factor~$\alpha$ as in Equation~\eqref{eq:alpha}, using the safety factor~$\sigma$ and the lower bound~$\alpha_{\min}$.} \label{algline:alpha}
		\State Refine slab width by $\alpha$, i.e., $w \gets \alpha\cdot w$. \label{algline:refine}
		\State Refine $\EMmax$ by $\alpha$, i.e., $\EMmax \gets \alpha\cdot\EMmax$.
		\Else
		\State Add a scenario in~$\argmin_{s\in S\setminus D, \AA_{\text{LB}}^s(\tilde x) < \tilde z} \HV_{\text{LB}}^s(\tilde x)$ to~$D$. \label{algline:addscenario}
		\EndIf
		
		\EndIf
		\EndWhile \label{alg-line:general-while-loop-end}
		\State \Return $\tilde x$.
	\end{algorithmic}
\end{algorithm}

The slab decomposition algorithm is stated in pseudocode in
Algorithm~\ref{alg:slabdecomp}. Its inputs are the scenario set~$S$;
a target gap~$\TG\in[0,1)$, specifying a prescribed upper bound on the
optimality gap of the returned first-stage decision; the initial slab
width~$w\in(0,\refpt_1]$ of the slab decomposition
(Section~\ref{sec:model}); a safety factor~$\sigma\in(0,1)$ and a lower
bound~$\alpha_{\min}\in(0,\sigma)$ on the factor by which the slab
width is reduced in a refinement step
(Section~\ref{sec:gapeval}); and a step size~$\delta>0$ for the
epsilon-constraint method used in the evaluation step. After an
initialization phase, the algorithm enters a loop
(Line~\ref{alg-line:general-while-loop}) that terminates once the computed
optimality gap is at most~$\TG$ or the master problem turns out to be infeasible, in which case the optimal objective value of Problem~\eqref{eq:hv-maxmin-problem} is zero.

During initialization, the algorithm calls the subroutine~\Init{} to obtain a small initial scenario subset~$D\subseteq S$. Possible implementations of~\Init{} are discussed at the end of this subsection, and the problem-specific choices used in our experiments are reported in Section~\ref{sec:applications}. Next, the algorithm initializes the quantity~$\EMmax$, an upper bound on the discretization error of the slab decomposition over all scenarios and all first-stage decisions, which is defined in Equation~\eqref{eq:EMmax} in Section~\ref{sec:gapeval}. This quantity does not depend on the master solution and can be computed up front by solving one small minimization problem per scenario. If~$\EMmax=0$, the optimal objective value of Problem~\eqref{eq:hv-maxmin-problem} is zero (see Section~\ref{sec:gapeval}). A relative optimality gap is then meaningless, so the algorithm terminates and reports that the optimal objective value is zero; in this case, every first-stage decision~$x\in\X$ is optimal.

In each iteration of the loop, the algorithm solves the master
problem~$\MAS_D$ to optimality or establishes its infeasibility in
Line~\ref{algline:solveMAS}. The master problem, defined in Section~\ref{sec:model}, is a linear approximation of~$\max_{x\in\X}\min_{s\in D}\HV^s(x)$ based on the slab decomposition. If it is infeasible, the optimal objective value of Problem~\eqref{eq:hv-maxmin-problem} is zero (see Section~\ref{sec:model}) and, as before, the algorithm terminates and reports this. Otherwise, let~$(\tilde x, \tilde z)$ denote the resulting master solution, where~$\tilde x$ is the first-stage decision and~$\tilde z$ is the optimal objective value of the master problem.

To determine whether~$\tilde x$ meets the target-gap requirement,
the algorithm computes, for every scenario~$s\in S$, lower bounds on
the hypervolume and the activated area for~$s$ and~$\tilde x$ in
Line~\ref{algline:findHV}. For each scenario~$s\in S$, we apply the epsilon-constraint method
(Algorithm~\ref{alg:eps-constraint}) to obtain a subset of the
nondominated set~$\PND^s(\tilde x)$ and one corresponding efficient
solution for each point in this subset.
We first remove nondominated points that do not weakly dominate the reference point together with their corresponding efficient solutions, and then use the resulting sets to determine lower bounds on~$\HV^s(\tilde x)$ and~$\AA^s(\tilde x)$, using Equation~\eqref{eq:hv-calculation} and the construction in Appendix~\ref{sec-app:AAcomb}, respectively. Here, $\AA^s(\tilde x)$ denotes the \emph{activated area} for~$s$ and~$\tilde x$, i.e., the objective value of the single-scenario master problem~$\MAS_{\{s\}}(\tilde x)$, which we formally introduce in Section~\ref{sec:model}. We denote these lower bounds by~$\HV_{\text{LB}}^s(\tilde x)$ and~$\AA_{\text{LB}}^s(\tilde x)$.

Once these quantities are known, the optimality gap is evaluated in Line~\ref{algline:gap} as described in Section~\ref{sec:gapeval} (Equation~\eqref{eq:gap}). If the gap is at most the target gap~$\TG$ (Line~\ref{algline:generalterminationcrit}), the algorithm terminates and returns~$\tilde x$.

Otherwise, the approximation has to be improved, and the algorithm distinguishes two cases. It first checks whether~$\AA_{\text{LB}}^s(\tilde x) \ge \tilde z$ holds for all scenarios~$s\in S\setminus D$. For every scenario~$s\in D$, the exact activated area already satisfies
$\AA^s(\tilde x)\geq\tilde z$; see
Lemma~\ref{lemma:AA-master-feas}. 
If so, adding any further scenario to~$D$ would provably leave the optimal objective value of the master problem unchanged and~$\tilde x$ optimal (we prove this in Section~\ref{sec:correctness}, Lemma~\ref{lemma:scenario-addition}). The algorithm therefore makes the slab decomposition more accurate by reducing the slab width~$w$ by a factor~$\alpha\in(0,1)$ (Line~\ref{algline:refine}). Since~$\EMmax$ scales linearly with~$w$, it is updated by the same factor.
Reducing the slab width entails recomputing the number~$n$ of slabs, the positions~$t^i$, and the index sets~$G^s$ for the new width. The factor~$\alpha$ is recomputed in every refinement step (Line~\ref{algline:alpha}) from the current master solution, the safety factor~$\sigma$, and the lower bound~$\alpha_{\min}$; we show in Section~\ref{sec:gapeval} how this is done. If, on the other hand, there is at least one scenario~$s\in S\setminus D$ with~$\AA_{\text{LB}}^s(\tilde x) < \tilde z$, the algorithm adds, among all such scenarios, one with the smallest hypervolume lower bound~$\HV_{\text{LB}}^s(\tilde x)$ to~$D$ (Line~\ref{algline:addscenario}). In both cases, the master problem is solved again with the refined decomposition or the extended scenario set, and the loop repeats.

\paragraph{Choice of \Init{}.}
The implementation of~\Init{} can be chosen in different ways. Possible options include:
\begin{enumerate}
	\item \Init{} chooses some first-stage decision~$x^0\in\X$ (e.g., heuristically), computes the lower bounds~$\HV_{\text{LB}}^s(x^0)$ for every scenario~$s\in S$ as in the evaluation step in Line~\ref{algline:findHV}, and returns a scenario attaining the smallest of these values as a proxy for a worst scenario for~$x^0$, using the same hypervolume-lower-bound criterion as the scenario-addition rule in Line~\ref{algline:addscenario}.
    In single-objective formulations with a separate scalar first-stage cost, this explicit initialization can sometimes be avoided by using $\Init()=\emptyset$~\citep{goerigk2025asbp}: there, the master objective consists of the first-stage cost plus an epigraph variable for the second-stage cost, so the master problem over the empty scenario set returns a first-stage decision with minimal first-stage cost, for which the first iteration then identifies a worst scenario. In our setting, this shortcut is not available because there is no separate scalar first-stage objective: the first-stage contributions are incorporated into the two objective functions~$f_1^s$ and~$f_2^s$. Consequently, for~$D=\emptyset$, the variable~$z$ in~$\MAS_D$ is unconstrained, so the master problem is unbounded and does not yield a candidate first-stage decision.
	\item \Init{} guesses a presumably worst scenario directly from the structure of the problem at hand, without solving any auxiliary problem; we state such rules for our two applications in Section~\ref{sec:applications}.
	\item \Init{} returns one or several randomly chosen scenarios.
\end{enumerate}
In our computational experiments in Section~\ref{sec:comp}, we use the second option.

\subsection{Slab decomposition and master problem}\label{sec:model}

A \emph{slab} of width~$w\in\R_{>0}$ and \emph{position}~$t\in\R_{\geq 0}$ is the rectangle $[t,t+w]\times[0,\refpt_2]$ in the objective space, where~$t$ is the horizontal coordinate of its left boundary. For each slab, the model chooses an \emph{activation height}
$h\in[0,\refpt_2]$, which defines the sub-rectangle
$[t,t+w]\times[\refpt_2-h,\refpt_2]$. We call this sub-rectangle a
\emph{vertical portion} of the slab and say that the slab is
\emph{activated} if~$h>0$. The area of the vertical portion contributes
to the hypervolume approximation.

For a slab width~$w\in(0,\refpt_1]$, we let $n\coloneqq\left\lfloor\frac{\refpt_1}{w}\right\rfloor$ denote the number of slabs and define the slab-index set $G\coloneqq[n]$. The slab with index~$i\in G$ has position $t^i\coloneqq\refpt_1-w\cdot(n-i+1)$. The resulting~$n$ adjacent vertical slabs partition the rectangle $[t^1,\refpt_1]\times[0,\refpt_2]$.\footnote{Alternatively, the objective space can be decomposed using horizontal slabs of a fixed height that span from~$0$ to~$\refpt_1$ in the horizontal direction. For ease of presentation, we restrict ourselves to vertical slabs. By symmetry, all results presented below apply analogously to horizontal slabs.} For each scenario~$s$, we use a separate subset~$G^s\subseteq G$ of slab indices, whose construction is described later.

Let~$D\subseteq S$ be a subset of scenarios. 
We develop a model that approximately solves Problem~\eqref{eq:hv-maxmin-problem} by constructing, for each scenario~$s\in D$, a lower bound on the hypervolume~$\HV^s(x)$. These lower bounds are obtained from the slab decomposition by activating a vertical portion of a slab whenever there exists a feasible solution whose corresponding point weakly dominates every point in that portion. More precisely, for a given first-stage decision~$x$, scenario~$s\in D$, and slab~$i\in G^s$, the model selects a second-stage decision~$y^{s,i}\in\Y^s(x)$ whose corresponding point lies on or to the left of the slab's left boundary. Its activation height~$h^{s,i}$ is chosen such that the lower boundary $\refpt_2-h^{s,i}$ is at least the vertical coordinate~$f_2^s(x,y^{s,i})$ of this point. Summing the areas of the activated portions over all slabs yields a lower bound on~$\HV^s(x)$ for each scenario~$s\in D$. We call the resulting model the \emph{master problem} and define it formally below as~$\MAS_D$. The variables and parameters we use in~$\MAS_D$ are as follows:

\medskip
\noindent
\begin{tabular}{wc{1.6cm}|l}
	\multicolumn{2}{l}{\emph{Variables:}} \\
	$x \in \X$     & first-stage decision  \\
	$y^{s,i} \in \Y^s(x)$ & second-stage decision for a scenario~$s$ and slab~$i$            \\
	$h^{s,i} \in \R_{\ge 0}$ & height of activation of slab $i$ for scenario $s$ (from the top) \\
    $a^{s,i} \in \{0, 1\}$ & activation indicator for slab~$i$ and scenario~$s$; $a^{s,i}=0$ forces~$h^{s,i}=0$, \\ &whereas $a^{s,i}=1$ permits activation \\
	$z \in \R$     & auxiliary variable (stemming from the epigraph reformulation)    \\[0.5ex]
	\multicolumn{2}{l}{\emph{Parameters:}} \\
    $D \subseteq S$ & subset of scenarios included in the master problem \\
	$G^s \subseteq G$   & set of slab indices for scenario $s$                             \\
	$w \in (0,\refpt_1]$     & slab width                                                       \\
	$t^i \in \R_{\ge 0}$   & coordinate of the left boundary of slab~$i$
\end{tabular}
\medskip

\noindent
The master problem~$\MAS_D$ then reads:
\begin{subequations}
\begin{align}
	\max \quad        & z                                                                                                                    \\
	\text{s.t.} \quad & z \le \sum_{i \in G^s} w \cdot h^{s,i}                            &  & \forall s \in D \label{eq:master-3}            \\[1ex]
	                  & h^{s,i} \le \refpt_2 \cdot a^{s,i}                                   &  & \forall s \in D, i \in G^s \label{eq:master-4} \\
	                  & f_1^s(x, y^{s,i}) \le (1-a^{s,i}) \cdot \refpt_1 + \,a^{s,i} \cdot t^i &  & \forall s \in D, i \in G^s \label{eq:master-5} \\
	                  & f_2^s(x, y^{s,i}) \le \refpt_2 - h^{s,i}                             &  & \forall s \in D, i \in G^s \label{eq:master-6} \\[1ex]
	                  & h^{s,i} \in \R_{\ge 0},\, a^{s,i} \in \{0, 1\}                       &  & \forall s \in D, i \in G^s                     \\[1ex]
	                  & y^{s,i} \in \Y^s(x)                                                &  & \forall s \in D, i \in G^s                     \\[1ex]
	                  & x \in \X, \quad z \in \R
\end{align}
\end{subequations}

If~$\MAS_D$ is infeasible, there exists no first-stage decision~$x\in\X$ such that, for each scenario~$s\in D$, the set~$f^s(\{x\},\Y^s)$ contains a point weakly dominating the chosen reference point~$\refpt$. Consequently, for every first-stage decision~$x\in\X$, there is a scenario~$s\in D$ with $\HV^s(x)=0$, so that the optimal objective value of Problem~\eqref{eq:hv-maxmin-problem} is zero and every first-stage decision is optimal (see Proposition~\ref{prop:finite-termination}~(a)). In this case, the algorithm terminates and reports this (see Line~\ref{algline:mas-infeasible} in Algorithm~\ref{alg:slabdecomp}).

Otherwise, we call an optimal solution of~$\MAS_D$ a \emph{master solution}. We denote a master solution by $(\tilde x,\tilde y,\tilde z,\tilde h,\tilde a)$ or, when only its first-stage decision and objective value are relevant, briefly by $(\tilde x,\tilde z)$. We regard~$\tilde y$, $\tilde h$, and~$\tilde a$ as block vectors indexed by the scenarios~$s\in D$, whose scenario-specific blocks are indexed by the slab indices~$i\in G^s$. For example, $\tilde h=(\tilde h^s)_{s\in D}$, where $\tilde h^s=(\tilde h^{s,i})_{i\in G^s}$.

To exclude slabs that cannot contribute to the hypervolume approximation and thereby reduce the size of the master problem, we restrict the decomposition for each scenario~$s\in S$ to the subset
\[
G^s \coloneqq
\bigl\{i\in G:
\exists (x,y)\in\XY:
f_1^s(x,y)\le t^i,\;
f_2^s(x,y)\le\refpt_2
\bigr\}
\subseteq G.
\]
Thus, for every slab~$i\in G^s$, there exists a feasible solution whose corresponding point lies on or to the left of the slab's left boundary and no higher than the reference point. Note that this construction is performed separately for each scenario and does not, by itself, guarantee feasibility of~$\MAS_D$, which requires a common first-stage
decision for all scenarios in~$D$.
We compute~$G^s$ by solving
\begin{align*}
    \underline f_1^s \coloneqq \min \quad
        & f_1^s(x,y) \\
    \text{s.t.}\quad
        & f_1^s(x,y) \le \refpt_1, \\
        & f_2^s(x,y) \le \refpt_2, \\
        & x\in\X,\quad y\in\Y^s(x).
\end{align*}
If this problem is infeasible or if $\underline f_1^s>t^n$, we set~$G^s\coloneqq\emptyset$. Otherwise, we let~$i'\in G$ be the smallest index satisfying $t^{i'}\geq\underline f_1^s$ and set $G^s\coloneqq\{i',\dots,n\}$.

To strengthen~$\MAS_D$ computationally, we additionally impose the following monotonicity constraints:
\begin{align}\label{constr:monotonicity}
    a^{s,i}\le a^{s,i+1}
    \qquad
    \forall\,s\in D,\;
    i\in G^s\setminus\{n\}.
\end{align}
These constraints do not change the optimal value of the master problem: whenever a slab~$i$ can be activated, the same second-stage decision is feasible for every subsequent slab because $t^{i+1}\geq t^i$. The following analysis is independent of whether these additional constraints are included. 

We now quantify how much of the hypervolume the master model captures for a single scenario and a fixed first-stage decision. To this end,
fix a first-stage decision~$x\in\X$ and a scenario~$s\in S$. For any
scenario subset~$D\subseteq S$, we write $\MAS_D(x)$ for the restriction of~$\MAS_D$ obtained by fixing the first-stage variables to~$x$. We consider the problem $\MAS_{\{s\}}(x)$ and let $\tilde y^s(x)$, $\tilde h^s(x)$, and~$\tilde a^s(x)$ denote the components corresponding to the variables~$y$, $h$, and~$a$, respectively, of an optimal solution.
We then define~$\AA^s(x) \coloneqq \sum_{i \in G^s} w \cdot \tilde h^{s,i}(x)$ as the \emph{activated area} (for~$s$ and~$x$). This equals the optimal objective value of~$\MAS_{\{s\}}(x)$.
The problem~$\MAS_{\{s\}}(x)$ can be infeasible, namely if $G^s\ne\emptyset$ but no $y\in\Y^s(x)$ yields a point in $[0,\refpt_1]\times[0,\refpt_2]$. For example, in the selection
problem from Section~\ref{sec:sel-problem}, this occurs for every
first-stage decision~$x$ that selects more than~$p$ items. In this case, we use the conventions $\AA^s(x)\coloneqq0$, $\AA_{\text{LB}}^s(x)\coloneqq0$, and $\tilde h^{s,i}(x)\coloneqq0$ for all~$i\in G^s$. Since no point for~$s$ and~$x$ then weakly dominates the reference point, we also have $\HV^s(x)=0$ and all statements below remain valid in this case (in particular, the error margin~$\EM^s(x)$ defined in Section~\ref{sec:error-estimation} is also zero).
Because~$\MAS_{\{s\}}(x)$ contains only one scenario, an optimal solution can alternatively be obtained without solving an MIP once the nondominated set~$\PND^s(x)$ is available, by deriving the activation heights directly from it. Applying the same construction to only a subset of the nondominated set yields a lower bound on~$\AA^s(x)$, as shown in Appendix~\ref{sec-app:AAcomb}.
For~$|D|>1$, however, the scenario-specific components of an optimal solution of~$\MAS_D(x)$ need not maximize the activated area for each individual scenario, because the objective maximizes only their common lower bound~$z$. Consequently, $\sum_{i\in G^s}w\cdot\tilde h^{s,i}(x)$ is generally only a lower bound on~$\AA^s(x)$, and the exact value~$\AA^s(x)$ must be determined by solving~$\MAS_{\{s\}}(x)$ separately or by using the combinatorial
construction described above.

The following lemma makes precise how the master problem decomposes across the scenarios for a fixed first-stage decision. It is used repeatedly in the analysis below.

\begin{lemma}\label{lemma:master-decomposition}
	Let~$\emptyset\neq D\subseteq S$ be such that~$\MAS_D$ is feasible, and let~$(\tilde x, \tilde z)$ be a master solution of~$\MAS_D$. Then
	\begin{equation*}
		\tilde z \;=\; \min_{s\in D} \AA^s(\tilde x) \;=\; \max_{x\in\X}\min_{s\in D} \AA^s(x).
	\end{equation*}
\end{lemma}
\begin{proof}
	Let~$(\tilde x, \tilde y, \tilde z, \tilde h, \tilde a)$ be a master solution of~$\MAS_D$. For each~$s\in D$, the scenario-specific vectors~$\tilde y^s$, $\tilde h^s$, and~$\tilde a^s$ together with~$\tilde x$ and the objective value set to~$\sum_{i\in G^s} w\cdot\tilde h^{s,i}$ form a feasible solution of~$\MAS_{\{s\}}(\tilde x)$. Consequently, $\AA^s(\tilde x) \ge \sum_{i\in G^s} w\cdot\tilde h^{s,i} \ge \tilde z$, where the second inequality follows from Constraint~\eqref{eq:master-3}. Hence, $\min_{s\in D}\AA^s(\tilde x) \ge \tilde z$.

	Conversely, let~$x\in\X$ be arbitrary. If~$\MAS_{\{s\}}(x)$ is feasible for every~$s\in D$, we combine the scenario-specific components of optimal solutions to the single-scenario problems~$\MAS_{\{s\}}(x)$ and set $z=\min_{s\in D}\AA^s(x)$. This yields a feasible solution of~$\MAS_D$, and therefore $\tilde z\geq\min_{s\in D}\AA^s(x)$.

    If~$\MAS_{\{s\}}(x)$ is infeasible for some~$s\in D$, then $\min_{s\in D}\AA^s(x)=0$ by our convention. Moreover, $\tilde z\geq0$, since replacing a negative~$z$-component of any feasible solution of~$\MAS_D$ by~$0$ preserves feasibility. Hence, $\tilde z\geq\min_{s\in D}\AA^s(x)$ also in this case.
    Since~$x\in\X$ was arbitrary, we obtain
	\begin{equation*}
		\tilde z \;\ge\; \max_{x\in\X}\min_{s\in D}\AA^s(x) \;\ge\; \min_{s\in D}\AA^s(\tilde x) \;\ge\; \tilde z,
	\end{equation*}
	so equality holds throughout.
\end{proof}

We next show that the activated area is a lower bound on the corresponding hypervolume.

\begin{lemma}\label{lemma:AVleHV}
    For every scenario~$s\in S$ and every first-stage decision
    $x\in\X$, we have
    \[
        \AA^s(x)\le\HV^s(x).
    \]
\end{lemma}

\begin{proof}
    If~$\MAS_{\{s\}}(x)$ is infeasible, the claim follows immediately from our convention~$\AA^s(x)=0$. Hence, assume that~$\MAS_{\{s\}}(x)$ is feasible. For ease of notation, we write~$\tilde y \coloneqq \tilde y^s(x)$, $\tilde h \coloneqq \tilde h^s(x)$, and~$\tilde a \coloneqq \tilde a^s(x)$ for the $y$-, $h$-, and $a$-components of the optimal solution of~$\MAS_{\{s\}}(x)$ introduced above.
    The region activated by this solution is
    \[
        A^s(x)
        \coloneqq
        \bigcup_{\substack{i\in G^s:\\
                           \tilde h^i>0}}
        [t^i,t^i+w]
        \times
        [\refpt_2-\tilde h^i,\refpt_2].
    \]
    Since the slabs have pairwise disjoint interiors, its area is
    \[
        \operatorname{area}(A^s(x))
        =
        \sum_{i\in G^s}w\cdot\tilde h^i
        =
        \AA^s(x).
    \]

    Consider any~$i\in G^s$ with
    $\tilde h^i>0$. Constraint~\eqref{eq:master-4} implies $\tilde a^i=1$, and Constraints~\eqref{eq:master-5} and~\eqref{eq:master-6} yield
    \[
        f_1^s(x,\tilde y^i)\le t^i
        \qquad\text{and}\qquad
        f_2^s(x,\tilde y^i)
        \le
        \refpt_2-\tilde h^i.
    \]
    Since also~$t^i+w\le\refpt_1$, the activated portion of slab~$i$ is contained in the rectangle dominated by the point $f^s(x,\tilde y^i)$ and bounded by~$\refpt$:
    \[
        [t^i,t^i+w]
        \times
        [\refpt_2-\tilde h^i,\refpt_2]
        \subseteq
        [f_1^s(x,\tilde y^i),\refpt_1]
        \times
        [f_2^s(x,\tilde y^i),\refpt_2].
    \]

    By Assumption~\ref{ass:external-stability}, the point $f^s(x,\tilde y^i)$ is weakly dominated by some nondominated point in~$\PND^s(x)$ (note that~$\PND^s(x)\subseteq\R^2_{\ge0}$ by the nonnegativity of the objectives, so that Definition~\ref{def:HV} applies). The rectangle on the right-hand side is therefore contained in~$R(\PND^s(x))$. Taking the union over all activated slabs gives $A^s(x)\subseteq R(\PND^s(x))$. Consequently,
    \[
        \AA^s(x)
        =
        \operatorname{area}(A^s(x))
        \le
        \operatorname{area}(R(\PND^s(x)))
        =
        \HV^s(x).
    \]
\end{proof}

\subsection{Discretization error and optimality gap}\label{sec:error-and-gap}

Lemma~\ref{lemma:AVleHV} shows that the activated area~$\AA^s(x)$ is a lower bound on the hypervolume~$\HV^s(x)$. Since the slab decomposition uses a positive slab width, $\AA^s(x)$ generally underestimates~$\HV^s(x)$. This discretization error can contribute to a nonzero optimality gap for the first-stage decision of a master solution. To compute this gap, we bound the error induced by the slab decomposition. We proceed in two steps: first, for a fixed scenario and first-stage decision, we bound the hypervolume from above by the activated area plus an \emph{error margin} (Section~\ref{sec:error-estimation}); second, we use this bound to derive computable lower and upper bounds on the optimal objective value of the overall problem, the resulting optimality gap, and the slab-refinement rule (Section~\ref{sec:gapeval}).

\subsubsection{Error margin}\label{sec:error-estimation}

\begin{figure}
	\centering
	\begin{tikzpicture}[scale=0.8, thick]

\input{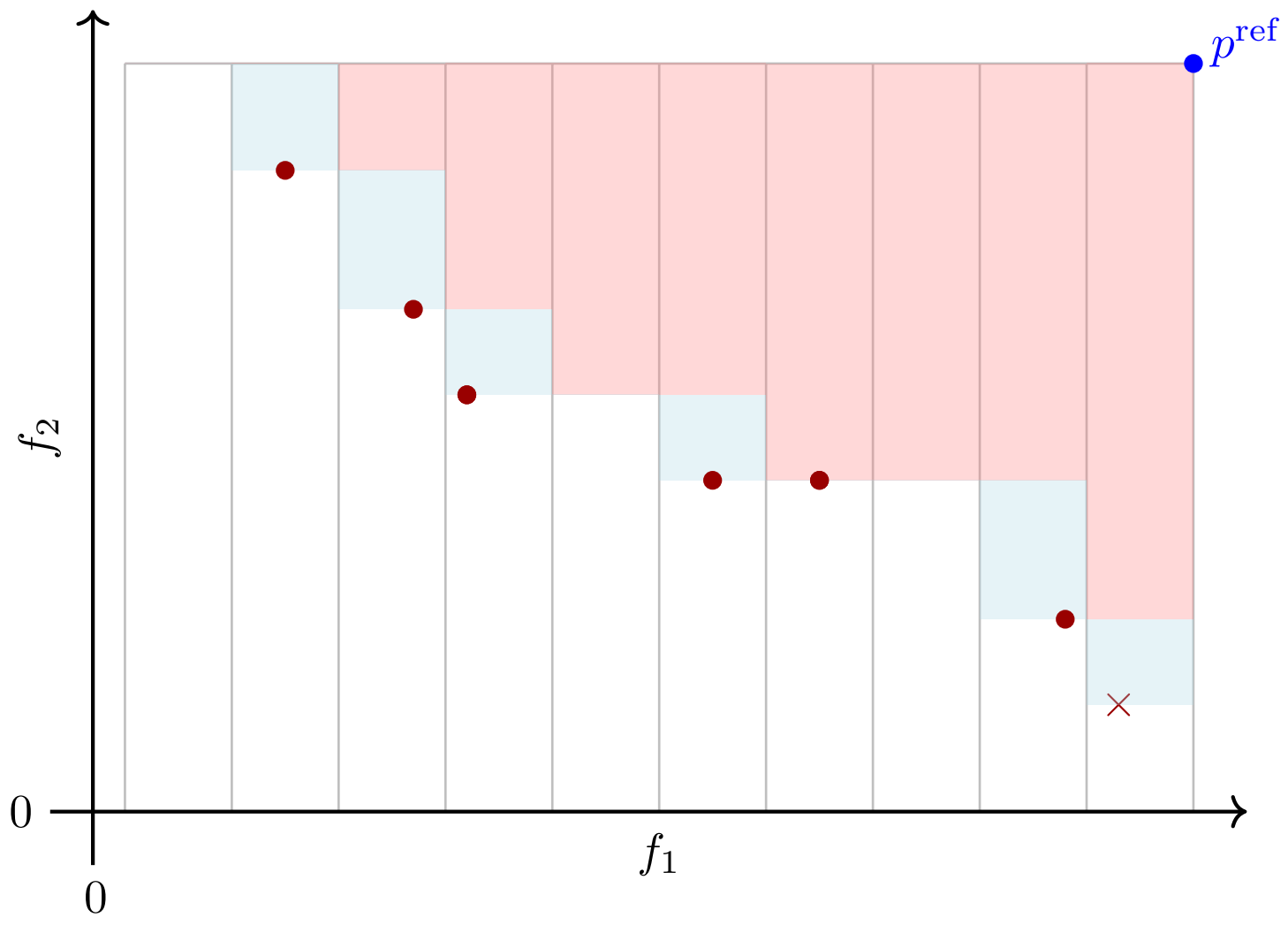}

\end{tikzpicture}
	\caption{Illustration of the activated area~$\AA^s(x)$ and the error margin~$\EM^s(x)$ for scenario~$s$ and first-stage decision~$x$. The red region shows the activated area, which consists of the vertical portions of the slabs activated in the master problem. These portions are based on feasible second-stage decisions whose corresponding points, shown in red, lie on or to the left of the left boundary of the slab they activate. Because of the positive slab width, this stepwise activated region generally underestimates the hypervolume~$\HV^s(x)$. The blue region represents the error margin used to bound this discretization error from above. It consists of the leftward extensions of the steps of the activated area and, in the rightmost slab, the additional term~$\addsn$ obtained from Problem~\eqref{eq:addn-point}.}
	\label{fig:grid-with-error-margin}
\end{figure}

Throughout this subsection, fix a first-stage decision~$x\in\X$ and a scenario~$s\in S$, and let~$\tilde h^s(x)$ be as in the definition of the activated area in Section~\ref{sec:model}. We introduce an \emph{error margin}~$\EM^s(x)$ that provides an upper bound on the discretization error~$\HV^s(x)-\AA^s(x)$ induced by the positive slab width. 
Geometrically, the error margin is obtained by extending each step of the activated area leftward by one slab width, as illustrated by the blue region in Figure~\ref{fig:grid-with-error-margin}. If the leftmost slab is activated, it is extended only to the vertical axis.
To define~$\EM^s(x)$ formally, we extend the activation heights to all
slab indices~$i\in G$ by setting
\begin{align*}
	\bar h^{s,i}(x) \coloneqq
	\begin{cases}
		0,               & \text{if } i \notin G^s, \\
		\tilde h^{s,i}(x), & \text{if } i \in G^s.
	\end{cases}
\end{align*}

Note that the master problem does not explicitly account for points whose first coordinate lies in the rightmost slab to the right of~$t^n$, since such points cannot activate any slab. Their potential contribution to the hypervolume within the rightmost slab must therefore be bounded separately. Since an exact treatment would lead to a quadratic model, we compute an upper bound by solving
\begin{align}\label{eq:addn-point}
	\begin{split}
		\min \quad & f_2^s(x, y) \\
		\text{s.t.} \quad & t^n \le f_1^s(x,y) \le t^n + w \\[1ex]
		& f_2^s(x,y) \le \refpt_2 - \bar h^{s,n}(x) \\
		& y \in \Y^s(x)
	\end{split}
\end{align}
If Problem~\eqref{eq:addn-point} has an optimal solution $\hat y$, we set 
\[
    \addsn
    \coloneqq
    w\cdot
    \bigl(
        \refpt_2-\bar h^{s,n}(x)-f_2^s(x,\hat y)
    \bigr);
\]
otherwise, we set~$\addsn\coloneqq0$.
The error margin is then given by
\begin{equation}\label{eq:EM}
    \begin{aligned}
        \EM^s(x)
        &=
        t^1\cdot\bar h^{s,1}(x)
        +
        \sum_{i=2}^{n}
        w\cdot
        \bigl(
            \bar h^{s,i}(x)-\bar h^{s,i-1}(x)
        \bigr)
        +
        \addsn \\
        &=
        t^1\cdot\bar h^{s,1}(x)
        +
        w\cdot
        \bigl(
            \bar h^{s,n}(x)-\bar h^{s,1}(x)
        \bigr)
        +
        \addsn.
    \end{aligned}
\end{equation}

The error margin depends on the slab width and is designed to vanish as the slab width approaches zero. We now prove that the hypervolume~$\HV^s(x)$ is indeed bounded above by the sum of the activated area~$\AA^s(x)$ and the error margin~$\EM^s(x)$.

\begin{proposition}\label{prop:AVplusEM}
	Let~$s \in S$ be a scenario and let~$x \in \X$ be a first-stage decision. For~$s$ and~$x$, let~$\HV^s(x)$ denote the hypervolume (with respect to the reference point $\refpt$), let $\AA^s(x)$ be the activated area (i.e., the objective value of $\MAS_{\{s\}}(x)$), and let $\EM^s(x)$ be the error margin as defined in Equation~\eqref{eq:EM}. Then
	\begin{equation*}
		\HV^s(x) \le \AA^s(x) + \EM^s(x).
	\end{equation*}
\end{proposition}

\begin{proof}
	Fix any~$s\in S$ and any~$x\in\X$. If~$\MAS_{\{s\}}(x)$ is infeasible, then $\AA^s(x)=\HV^s(x)=0$ by the conventions and observations in Section~\ref{sec:model}. Moreover, our convention $\tilde h^{s,i}(x)=0$ for all~$i\in G^s$ implies $\bar h^{s,i}(x)=0$ for all~$i\in G$, and Problem~\eqref{eq:addn-point} is infeasible, since any feasible solution would yield a point in $[0,\refpt_1]\times[0,\refpt_2]$. Thus, $\addsn=0$ and Equation~\eqref{eq:EM} gives~$\EM^s(x)=0$, so the claim follows. Hence, assume in the following that~$\MAS_{\{s\}}(x)$ is feasible.
    For~$\tau\in[0,\refpt_1]$, define
	\begin{equation*}
		\mu(\tau)\coloneqq\min\{f_2^s(x,y):y\in\Y^s(x),\,f_1^s(x,y)\le\tau\}
	\end{equation*}
	(with~$\mu(\tau)\coloneqq+\infty$ if no such~$y$ exists), and let
	\begin{equation*}
		\phi(\tau)\coloneqq\max\{0,\,\refpt_2-\mu(\tau)\}\in[0,\refpt_2].
	\end{equation*}
	Since~$\mu$ is non-increasing, $\phi$ is non-decreasing. We use~$t^{n+1}\coloneqq\refpt_1$ and recall~$t^i=\refpt_1-w\,(n-i+1)$, so that~$t^{i+1}-t^i=w$ for~$i\in[n]$ and~$0\le t^1\le w$.

    \emph{Step 1 (hypervolume as an integral).}
    First note that, by Assumption~\ref{ass:external-stability}, the region weakly dominated by all feasible points for~$s$ and~$x$ and bounded by~$\refpt$ coincides with~$\mathcal R(\PND^s(x))$. Hence, in the following, considering all feasible points for~$s$ and~$x$ is equivalent to considering only the nondominated points for~$s$ and~$x$ when determining the hypervolume.
    
    A point~$(\tau,\eta)\in[0,\refpt_1]\times[0,\refpt_2]$ is weakly
    dominated by some feasible point~$f^s(x,y)$, $y\in\Y^s(x)$, if and
    only if there exists~$y\in\Y^s(x)$ with
    $f_1^s(x,y)\le\tau$ and~$f_2^s(x,y)\le\eta$, i.e., if and only if
    $\eta\ge\mu(\tau)$. Hence, for each fixed~$\tau$, the vertical cross-section of the dominated region at first coordinate~$\tau$, i.e., the set of values of~$\eta$ for which~$(\tau,\eta)$ belongs to the dominated region, is~$[\mu(\tau),\refpt_2]$ if $\mu(\tau)\le\refpt_2$ and is empty otherwise. Its length is therefore exactly~$\phi(\tau)$.
    Integrating the lengths of these cross-sections over~$\tau$ yields
    \begin{equation*}
        \HV^s(x)
        =\operatorname{area}\bigl(\mathcal R(\PND^s(x))\bigr)
        =\int_0^{\refpt_1}\phi(\tau)\,d\tau.
    \end{equation*}

	\emph{Step 2 (activated area as a left Riemann sum).}
    For every~$i\in G^s$, any feasible activation height satisfies
    \begin{equation*}
        h^{s,i}
        \le
        \max\bigl\{0,\,\refpt_2-\mu(t^i)\bigr\}
        =
        \phi(t^i).
    \end{equation*}
    Indeed, if~$h^{s,i}>0$, Constraint~\eqref{eq:master-4} implies $a^{s,i}=1$, and Constraints~\eqref{eq:master-5} and~\eqref{eq:master-6} yield $f_1^s(x,y^{s,i})\le t^i$ and $f_2^s(x,y^{s,i})\le\refpt_2-h^{s,i}$, respectively, so that $h^{s,i}\le\refpt_2-\mu(t^i)$. If~$h^{s,i}=0$, the inequality holds trivially.
    
    Conversely, this bound can be attained simultaneously for all $i\in G^s$. If~$\phi(t^i)>0$, we choose a minimizer~$y$ in the definition of~$\mu(t^i)$ and set~$y^{s,i}=y$, $a^{s,i}=1$, and $h^{s,i}=\phi(t^i)$. If~$\phi(t^i)=0$, we set~$h^{s,i}=0$ and $a^{s,i}=0$, and the remaining slab-specific requirements on~$y^{s,i}$ then reduce to
    \[
        y^{s,i}\in\Y^s(x),\qquad
        f_1^s(x,y^{s,i})\le\refpt_1,\qquad
        f_2^s(x,y^{s,i})\le\refpt_2.
    \]
    A decision~$y^{s,i}$ satisfying these requirements exists because $\MAS_{\{s\}}(x)$ is feasible, and the same decision can be assigned to~$y^{s,j}$ for all indices~$j\in G^s$ with $\phi(t^j)=0$. Moreover, if the monotonicity constraints~\eqref{constr:monotonicity} are included, they are satisfied because~$\phi$ is non-decreasing. Since~$w>0$ and, for the single scenario~$s$, maximizing~$z$ amounts to maximizing the sum of the activation areas, optimality therefore yields
	\begin{equation}\label{eq:height-char}
		\tilde h^{s,i}(x)=\max\bigl\{0,\,\refpt_2-\mu(t^i)\bigr\}=\phi(t^i) \qquad \text{for every } i\in G^s.
	\end{equation}
    For~$i\in G\setminus G^s$, the definition of~$G^s$ implies that there is no feasible pair~$(x',y)$ with $f_1^s(x',y)\le t^i$ and~$f_2^s(x',y)\le\refpt_2$. Hence, for the fixed first-stage decision~$x$, either no $y\in\Y^s(x)$ satisfies~$f_1^s(x,y)\le t^i$, or every such~$y$ satisfies~$f_2^s(x,y)>\refpt_2$. In either case, $\phi(t^i)=0=\bar h^{s,i}(x)$.
    Combining this with~\eqref{eq:height-char}, we obtain $\bar h^{s,i}(x)=\phi(t^i)$ for all~$i\in[n]$, and therefore
	\begin{equation*}
		\AA^s(x)=\sum_{i\in G^s}w\cdot\tilde h^{s,i}(x)=\sum_{i=1}^{n}w\cdot\phi(t^i),
	\end{equation*}
	which is the left-endpoint Riemann sum of~$\int\phi$ over the slabs~$[t^i,t^{i+1}]$, $i\in[n]$.

	\emph{Step 3 (bounding the error).} Splitting the integral in Step~1 into the leftmost strip~$[0,t^1]$ and the~$n$ slabs, and subtracting the sum from Step~2,
	\begin{equation*}
		\HV^s(x)-\AA^s(x)=\int_0^{t^1}\!\!\phi(\tau)\,d\tau+\sum_{i=1}^{n}\left(\int_{t^i}^{t^{i+1}}\!\!\phi(\tau)\,d\tau-w\cdot\phi(t^i)\right).
	\end{equation*}
	Because~$\phi$ is non-decreasing and~$0\le\phi\le\refpt_2$, we bound the three types of terms separately:
	\begin{itemize}
		\item \emph{Strip:} We have $\int_0^{t^1}\phi(\tau)\,d\tau\le t^1\cdot\phi(t^1)$.
        \item \emph{Slabs~$i\in[n]\setminus\{n\}$:} $\phi(\tau)\le\phi(t^{i+1})$ for~$\tau\le t^{i+1}$, hence $\int_{t^i}^{t^{i+1}}\phi(\tau)\,d\tau-w\cdot\phi(t^i)\le w\cdot\bigl(\phi(t^{i+1})-\phi(t^i)\bigr)$. Summing over~$i\in[n]\setminus\{n\}$ telescopes to~$w\cdot\bigl(\phi(t^n)-\phi(t^1)\bigr)$.
        \item \emph{Last slab~$i=n$:} Since~$\phi$ is non-decreasing,
        \begin{align*}
            \int_{t^n}^{\refpt_1}\phi(\tau)\,d\tau
            -w\cdot\phi(t^n)
            =
            \int_{t^n}^{\refpt_1}
            \bigl(\phi(\tau)-\phi(t^n)\bigr)\,d\tau 
            \le
            w\cdot\bigl(\phi(\refpt_1)-\phi(t^n)\bigr)
            =\addsn.
        \end{align*}
        Indeed, if~$\phi(\refpt_1)>\phi(t^n)$, a minimizer~$y$ in the definition of~$\mu(\refpt_1)$ must satisfy $f_1^s(x,y)>t^n$, since otherwise it would also be feasible in the definition of~$\mu(t^n)$. Moreover, using $\bar h^{s,n}(x)=\phi(t^n)$, we have $f_2^s(x,y)=\mu(\refpt_1)<\refpt_2-\bar h^{s,n}(x)$. Hence,~$y$ is feasible for Problem~\eqref{eq:addn-point}. Conversely, every feasible solution of Problem~\eqref{eq:addn-point} is feasible in the definition of~$\mu(\refpt_1)$. Thus, an optimal solution~$\hat y$ of Problem~\eqref{eq:addn-point} satisfies $f_2^s(x,\hat y)=\mu(\refpt_1)$, and the last equality follows.
        If~$\phi(\refpt_1)=\phi(t^n)$, then, since~$\phi$ is non-decreasing, $\phi(\tau)=\phi(t^n)$ for all $\tau\in[t^n,\refpt_1]$. Hence,
        \[
            \int_{t^n}^{\refpt_1}
            \bigl(\phi(\tau)-\phi(t^n)\bigr)\,d\tau=0.
        \]
        Moreover, Problem~\eqref{eq:addn-point} is either infeasible, in which case~$\addsn=0$ by definition, or its optimal value equals $\refpt_2-\phi(t^n)$, since any strictly smaller value would imply $\phi(\refpt_1)>\phi(t^n)$. Thus, $\addsn=0$ also in this case.
	\end{itemize}
	Adding the three bounds and using~$\bar h^{s,i}(x)=\phi(t^i)$,
	\begin{align*}
		\HV^s(x)-\AA^s(x)
		&\le t^1\cdot\phi(t^1)+w\cdot\bigl(\phi(t^n)-\phi(t^1)\bigr)+\addsn\\
		&=(t^1-w)\cdot\bar h^{s,1}(x)+w\cdot\bar h^{s,n}(x)+\addsn=\EM^s(x),
	\end{align*}
	which is the claimed inequality.
\end{proof}

\subsubsection{Optimality gap and slab refinement}\label{sec:gapeval}

We now use the error-margin bound to derive computable bounds on the
optimal objective value of the overall problem and the resulting
optimality gap. We subsequently use these bounds to determine the
factor by which the slab width is reduced in a refinement step. To
state these bounds, let~$\OPT$ denote the optimal objective value of
the overall problem, i.e.,
\begin{equation*}
	\OPT \coloneqq \max_{x\in\X}\min_{s\in S} \HV^s(x).
\end{equation*}

For a current master solution~$(\tilde x,\tilde z)$, the smallest hypervolume over all scenarios for the first-stage decision~$\tilde x$ provides a lower bound on~$\OPT$. An upper bound will be obtained from the current master objective value~$\tilde z$ plus a uniform upper bound on the error margin over all scenarios~$s\in S$ and all first-stage decisions~$x\in\X$. We first derive such a uniform upper bound, which is independent of the current master solution and can therefore be computed in advance.

\begin{lemma}\label{lemma:EM_bound}
    For each~$s\in S$ and each~$x\in\X$, we have
    \begin{equation*}
        \EM^s(x)
        \le
        \max\left(
            \{0\}
            \cup
            \left\{
                w\cdot\bigl(\refpt_2-f_2^s(x',y)\bigr):
                (x',y)\in\X\times\Y^s,\,
                f_1^s(x',y)\le\refpt_1
            \right\}
        \right).
    \end{equation*}
\end{lemma}

\begin{proof}
    Let~$M^s$ denote the right-hand side of the claimed inequality.
    By the definition of~$\EM^s(x)$  in~\eqref{eq:EM} and~$t^1\le w$, we have
    \begin{equation*}
        \EM^s(x)
        =
        (t^1-w)\cdot\bar h^{s,1}(x)
        +w\cdot\bar h^{s,n}(x)
        +\addsn
        \le
        w\cdot\bar h^{s,n}(x)+\addsn.
    \end{equation*}

    We distinguish two cases based on whether Problem~\eqref{eq:addn-point} has an optimal solution. If an optimal solution~$\hat y$ exists, then, by the definition of~$\addsn$,
    \begin{align*}
        \EM^s(x)
        &\le
        w\cdot\bar h^{s,n}(x)
        +w\cdot
        \bigl(
            \refpt_2-\bar h^{s,n}(x)-f_2^s(x,\hat y)
        \bigr)\\
        &=
        w\cdot\bigl(\refpt_2-f_2^s(x,\hat y)\bigr)
        \le M^s,
    \end{align*}
    where the last inequality follows because $(x,\hat y)\in\X\times\Y^s$ and, by the constraints of Problem~\eqref{eq:addn-point} and the definition of~$t^n$, $f_1^s(x,\hat y)\le t^n+w=\refpt_1$.

    If Problem~\eqref{eq:addn-point} has no optimal solution, then
    $\addsn=0$, so
    $\EM^s(x)\le w\cdot\bar h^{s,n}(x)$.
    If~$\bar h^{s,n}(x)=0$, then
    $\EM^s(x)\le0\le M^s$.
    Otherwise, $\bar h^{s,n}(x)>0$ implies~$n\in G^s$. Moreover,
    $\MAS_{\{s\}}(x)$ must be feasible, since otherwise $\tilde h^{s,i}(x)=0$ for all~$i\in G^s$ by convention. Hence, Equation~\eqref{eq:height-char} from the proof of Proposition~\ref{prop:AVplusEM} applies. With~$\mu$ and~$\phi$ defined as in that proof, we thus obtain
    \[
        0<\bar h^{s,n}(x)=\tilde h^{s,n}(x)=\phi(t^n).
    \]
    By the definition of~$\phi$, this implies that $\phi(t^n)=\refpt_2-\mu(t^n)$. By the definition of~$\mu(t^n)$, there therefore exists a $y\in\Y^s(x)$ with
    \[
        f_1^s(x,y)\le t^n\le\refpt_1
        \qquad\text{and}\qquad
        f_2^s(x,y)=\mu(t^n)
        =\refpt_2-\bar h^{s,n}(x).
    \]
    Hence,
    \begin{equation*}
        \EM^s(x)
        \le
        w\cdot\bar h^{s,n}(x)
        =
        w\cdot\bigl(\refpt_2-f_2^s(x,y)\bigr)
        \le M^s.
    \end{equation*}
\end{proof}

\noindent
We define~$\EMmax$ as a uniform upper bound on the error margin over all scenarios and all first-stage decisions:
\begin{equation}\label{eq:EMmax}
    \begin{aligned}
    \EMmax \coloneqq \max_{s\in S} \max\bigl(
            \{0\}
            \cup
            \bigl\{
                w\cdot\bigl(\refpt_2-f_2^s(x,y)\bigr):\;
                &(x,y)\in\X\times\Y^s,\,\\
                &f_1^s(x,y)\le\refpt_1
            \bigr\}
        \bigr)
    \end{aligned}
\end{equation}
This quantity does not depend on the current master solution and can therefore be computed at the beginning of the algorithm by solving, for each scenario~$s\in S$, a relatively easy linear MIP that minimizes~$f_2^s(x,y)$ over~$(x,y)\in\X\times\Y^s$ subject to $f_1^s(x,y)\le\refpt_1$. If this problem is infeasible, the corresponding scenario-specific bound is zero.

If~$\EMmax=0$, then no scenario~$s\in S$ admits a solution~$(x,y)\in\X\times\Y^s$ satisfying $f_1^s(x,y)\le\refpt_1$ and $f_2^s(x,y)<\refpt_2$. In particular, for every~$s\in S$, the set~$f^s(\X,\Y^s)$ does not intersect the box~$[0,\refpt_1)\times[0,\refpt_2)$, so that $\HV^s(x)=0$ for every~$s\in S$ and~$x\in\X$. Hence, the optimal objective value of Problem~\eqref{eq:hv-maxmin-problem} is zero and every first-stage decision is optimal (see Proposition~\ref{prop:finite-termination}~(a)). In this case, the algorithm terminates and reports this (see Line~\ref{algline:emmax-zero} in Algorithm~\ref{alg:slabdecomp}).

\smallskip

In the following theorem, we use scenario-wise hypervolume lower bounds for the current first-stage decision to obtain a lower bound on~$\OPT$. For an upper bound, we combine Proposition~\ref{prop:AVplusEM}, the characterization of the master-solution objective value in Lemma~\ref{lemma:master-decomposition}, and the uniform error-margin bound from Lemma~\ref{lemma:EM_bound}.

\begin{theorem}\label{thm:bounds}
    Let~$(\tilde x,\tilde z)$ be a master solution of~$\MAS_D$
    for some nonempty scenario subset~$\emptyset\neq D\subseteq S$,
    and let~$\OPT$ be the optimal objective value of the
    overall problem. For each scenario~$s\in S$, let
    $\HV_{\mathrm{LB}}^s(\tilde x)$ be a valid lower bound on the
    hypervolume for~$s$ and~$\tilde x$, and let~$\EMmax$ be defined
    as in Equation~\eqref{eq:EMmax}. Then
    \[
        \min_{s\in S}\HV_{\mathrm{LB}}^s(\tilde x)
        \;\le\;
        \OPT
        \;\le\;
        \tilde z+\EMmax.
    \]
\end{theorem}
\begin{proof}
    The lower bound follows directly from the definition of~$\OPT$:
	\begin{equation*}
		\min_{s\in S} \HV_{\text{LB}}^s(\tilde x) \le \min_{s\in S} \HV^s(\tilde x) \le \max_{x\in\X}\min_{s\in S} \HV^s(x) = \OPT.
	\end{equation*}
    The upper bound is shown as follows:
	\begin{align*}
		\OPT \quad
		 & = \quad\max_{x\in\X}\min_{s\in S} \HV^s(x)                                        \\
		 & \stackrel{\mathmakebox[\widthof{$\le$}]{\text{Prop. \ref{prop:AVplusEM}}}}{\le} \quad
		\max_{x\in\X}\min_{s\in S} \left(\AA^s(x)+\EM^s(x)\right)                            \\
		 & \stackrel{\mathmakebox[\widthof{$\le$}]{}}{\le} \quad
		\max_{x\in\X}\left[\min_{s\in S}\AA^s(x) + \max_{s\in S}\EM^s(x)\right]              \\
		 & \stackrel{\mathmakebox[\widthof{$\le$}]{}}{\le} \quad
		\max_{x\in\X}\min_{s\in S}\AA^s(x) + \max_{x\in\X}\max_{s\in S}\EM^s(x)              \\
		 & \stackrel{\mathmakebox[\widthof{$\le$}]{}}{\le} \quad
		\max_{x\in\X}\min_{s\in D}\AA^s(x) + \max_{x\in\X}\max_{s\in S}\EM^s(x)              \\
		 & \stackrel{\mathmakebox[\widthof{$\le$}]{\text{Lem. \ref{lemma:master-decomposition}}}}{=} \quad
		\tilde z + \max_{x\in\X}\max_{s\in S}\EM^s(x)                                        \\
		 & \stackrel{\mathmakebox[\widthof{$\le$}]{\text{Lem. \ref{lemma:EM_bound}}}}{\le} \quad
		\tilde z + \EMmax.
	\end{align*}
\end{proof}

For a maximization problem with a positive lower bound~$\LB$ and an upper bound~$\UB$ on the optimal objective value, the corresponding relative optimality gap is $\mathrm{gap}\coloneqq \frac{\UB}{\LB}-1$ (see, e.g., \cite{MIP-gap-gurobi}). Using the bounds from Theorem~\ref{thm:bounds}, we therefore define the computed optimality gap of a master solution~$(\tilde x,\tilde z)$ as
\begin{equation}\label{eq:gap}
    \mathrm{gap}
    \coloneqq
    \frac{\tilde z+\EMmax}
         {\min_{s\in S}\HV_{\mathrm{LB}}^s(\tilde x)}
    -1.
\end{equation}
If~$\min_{s\in S}\HV_{\mathrm{LB}}^s(\tilde x)=0$, we use the convention~$\mathrm{gap}\coloneqq+\infty$.

If the computed gap exceeds the target gap~$\TG$, reducing the slab width may decrease the gap. However, it cannot in general be predicted a priori how~$\tilde z$ and $\min_{s\in S}\HV_{\mathrm{LB}}^s(\tilde x)$ change when the slab width is modified; in particular, either quantity may increase or decrease. By contrast,~$\EMmax$ scales linearly with the slab width: multiplying~$w$ by a factor~$\alpha$ also multiplies~$\EMmax$ by $\alpha$. We therefore estimate the refinement factor~$\alpha$ by assuming that~$\tilde z$ and $\min_{s\in S}\HV_{\mathrm{LB}}^s(\tilde x)$ remain unchanged and choosing~$\alpha$ such that the resulting estimated gap equals~$\TG$. Provided that $\min_{s\in S}\HV_{\mathrm{LB}}^s(\tilde x)>0$ and~$\EMmax>0$, this gives
\begin{equation}\label{eq:alpha-estimation}
    \TG
    =
    \frac{\tilde z+\alpha\cdot\EMmax}
         {\min_{s\in S}\HV_{\mathrm{LB}}^s(\tilde x)}
    -1
    \quad\Longrightarrow\quad
    \alpha
    =
    \frac{
        (1+\TG)\min_{s\in S}\HV_{\mathrm{LB}}^s(\tilde x)-\tilde z
    }{\EMmax}.
\end{equation}

Since~$\tilde z$ and $\min_{s\in S}\HV_{\mathrm{LB}}^s(\tilde x)$ may change after the slab width is refined and the master problem is re-solved, the estimate in Equation~\eqref{eq:alpha-estimation} need not yield the refinement factor that is actually required. If~$\alpha$ is chosen too large, the refinement may be insufficient, so that additional iterations are needed. If, on the other hand,~$\alpha$ is chosen smaller than necessary, the resulting slab decomposition may be unnecessarily fine, increasing the time required to solve the master problem. To balance this trade-off, we introduce a \emph{safety factor}~$\sigma\in(0,1)$ and multiply the estimate in Equation~\eqref{eq:alpha-estimation} by~$\sigma$. This biases the refinement toward a slightly smaller slab width and thereby reduces the risk of an insufficient refinement. We exclude~$\sigma=1$ in order to establish finite termination of the algorithm in Section~\ref{sec:correctness}.
This safety factor~$\sigma$ is part of the input of Algorithm~\ref{alg:slabdecomp} and stays fixed throughout the run. We investigate the influence of~$\sigma$ experimentally in Section~\ref{sec:comp-alpha}.

The estimate in Equation~\eqref{eq:alpha-estimation} can become arbitrarily small or even nonpositive, since its numerator
\[
    (1+\TG)\min_{s\in S}\HV_{\mathrm{LB}}^s(\tilde x)-\tilde z
\]
can have either property. This also includes the case $\min_{s\in S}\HV_{\mathrm{LB}}^s(\tilde x)=0$, for which the computed gap is~$+\infty$ and Equation~\eqref{eq:alpha-estimation} was not derived directly. To ensure that the slab width remains positive, we introduce a fixed lower bound $\alpha_{\min}\in(0,\sigma)$ on the refinement factor, which is also part of the input of Algorithm~\ref{alg:slabdecomp}. Thus, we use
\begin{equation}\label{eq:alpha}
    \alpha
    =
    \max\left\{
        \sigma\cdot
        \frac{
            (1+\TG)\min_{s\in S}\HV_{\mathrm{LB}}^s(\tilde x)-\tilde z
        }{\EMmax},
        \alpha_{\min}
    \right\}.
\end{equation}
Note that this is well defined, since~$\EMmax>0$ holds in every iteration of Algorithm~\ref{alg:slabdecomp}: the case~$\EMmax=0$ is detected in Line~\ref{algline:emmax-zero} before the loop is entered, and each refinement multiplies~$\EMmax$ by a factor~$\alpha\ge\alpha_{\min}>0$.

Whenever the algorithm refines the slab width, the target gap has not been reached in the current iteration. If $\min_{s\in S}\HV_{\mathrm{LB}}^s(\tilde x)>0$, then
$\mathrm{gap}<+\infty$, and hence
\begin{align*}
    \TG < \mathrm{gap}
    &\enspace\Rightarrow\enspace
    \TG
    <
    \frac{\tilde z+\EMmax}
         {\min_{s\in S}\HV_{\mathrm{LB}}^s(\tilde x)}
    -1
    \\
    &\enspace\Rightarrow\enspace
    \frac{
        (1+\TG)\min_{s\in S}\HV_{\mathrm{LB}}^s(\tilde x)-\tilde z
    }{\EMmax}
    <1
    \\
    &\enspace\overset{\eqref{eq:alpha}}{\Rightarrow}\enspace
    \alpha<\sigma,
\end{align*}
where the last implication also uses
$\alpha_{\min}<\sigma$. If $\min_{s\in S}\HV_{\mathrm{LB}}^s(\tilde x)=0$, the first term in the maximum in~\eqref{eq:alpha} is nonpositive, and therefore $\alpha=\alpha_{\min}<\sigma$. Consequently, whenever the slab width is refined, the refinement factor satisfies
\[
    \alpha\in[\alpha_{\min},\sigma).
\]
Since~$\EMmax$ scales linearly with the slab width, it is multiplied by the same factor. In particular, after each refinement, $\EMmax$ is less than~$\sigma$ times its previous value.

\subsection{Correctness and finite termination}\label{sec:correctness}

Recall from the outline in Section~\ref{sec:outline} that, in each iteration in which the master problem is feasible and the target gap is not yet reached, the algorithm either adds a scenario or refines the slab decomposition. The following two lemmas justify this branching. The first one shows that, for the scenarios that are already contained in~$D$, the activated area for the current first-stage decision~$\tilde x$ is at least the master objective
value~$\tilde z$. Hence, the check in
Line~\ref{algline:branchcondition} of
Algorithm~\ref{alg:slabdecomp} only needs to consider scenarios in
$S\setminus D$.

\begin{lemma}\label{lemma:AA-master-feas}
	Let~$\emptyset\neq D\subseteq S$, and let~$(\tilde x,\tilde z)$
    be a master solution of~$\MAS_D$. Then $\AA^s(\tilde x) \ge \tilde z$ for every scenario~$s\in D$.
\end{lemma}
\begin{proof}
	By Lemma~\ref{lemma:master-decomposition}, we have $\tilde z=\min_{s\in D}\AA^s(\tilde x)$. Hence, $\AA^s(\tilde x) \ge \tilde z$ for every~$s\in D$.
\end{proof}

In Line~\ref{algline:branchcondition}, the algorithm checks whether $\AA_{\LB}^s(\tilde x)\ge\tilde z$ holds for every scenario $s\in S\setminus D$. Since the epsilon-constraint method returns a subset of the nondominated set, the resulting activated-area estimate is a lower bound, i.e., $\AA_{\LB}^s(\tilde x)\le\AA^s(\tilde x)$ (see Appendix~\ref{sec-app:AAcomb}). Hence, a successful check implies $\AA^s(\tilde x)\ge\tilde z$ for all~$s\in S\setminus D$ and, by Lemma~\ref{lemma:AA-master-feas}, for all~$s\in S$.

If~$\tilde z>0$, then $\AA^s(\tilde x)\ge\tilde z>0$ implies that $\MAS_{\{s\}}(\tilde x)$ is feasible for every~$s\in S\setminus D$, since we use the convention $\AA^s(\tilde x)=0$ whenever this problem is infeasible; see Section~\ref{sec:model}. If~$\tilde z=0$, this feasibility need not follow from the successful check. In particular, if~$\MAS_{\{s\}}(\tilde x)$ is infeasible for some $s\in S\setminus D$, then $\HV^s(\tilde x)=0$ and, consequently, $\HV_{\LB}^s(\tilde x)=0$, so the computed gap is~$+\infty$ and the algorithm cannot terminate in the current iteration. The following lemma shows that, whenever the single-scenario problems are feasible, a successful check guarantees that adding further scenarios to~$D$ does not change the optimal objective value of the master problem and leaves~$\tilde x$ optimal.

\begin{lemma}\label{lemma:scenario-addition}
    Let~$\emptyset\neq D\subseteq S$, let~$(\tilde x,\tilde z)$ be a master solution of~$\MAS_D$, and assume that, for every scenario $s\in S\setminus D$, the problem~$\MAS_{\{s\}}(\tilde x)$ is feasible and $\AA^s(\tilde x)\ge\tilde z$.
    Then, for every scenario subset~$D'$ with $D\subseteq D'\subseteq S$, the problem~$\MAS_{D'}$ is feasible, its optimal objective value equals~$\tilde z$, and~$\tilde x$ is the first-stage decision of a master solution of~$\MAS_{D'}$.
\end{lemma}
\begin{proof}
	Let~$D\subseteq D' \subseteq S$. The problem~$\MAS_{\{s\}}(\tilde x)$ is feasible for every~$s\in D'$: for~$s\in D$, this follows from the master solution of~$\MAS_D$ (see the proof of Lemma~\ref{lemma:master-decomposition}), and for~$s\in D'\setminus D$, it holds by assumption.
    Combining optimal solutions of the problems~$\MAS_{\{s\}}(\tilde x)$, $s\in D'$, and setting the master objective variable to $\min_{s\in D'}\AA^s(\tilde x)$ therefore yields a feasible solution of~$\MAS_{D'}$ with first-stage decision~$\tilde x$.
    In particular, $\MAS_{D'}$ is feasible; let~$\tilde z'$ denote its optimal objective value.

	By Lemma~\ref{lemma:master-decomposition} applied to~$D'$ and to~$D$, and since~$D\subseteq D'$,
	\begin{equation*}
		\tilde z' \;=\; \max_{x\in\X}\min_{s\in D'}\AA^s(x) \;\le\; \max_{x\in\X}\min_{s\in D}\AA^s(x) \;=\; \tilde z.
	\end{equation*}
	Conversely, the assumption together with Lemma~\ref{lemma:AA-master-feas} yields $\AA^s(\tilde x)\ge\tilde z$ for all~$s\in D'$, so the feasible solution of~$\MAS_{D'}$ constructed above has objective value $\min_{s\in D'}\AA^s(\tilde x)\ge\tilde z$. Hence, $\tilde z' = \tilde z$, and the constructed solution is optimal. Thus,~$\tilde x$ is the first-stage decision of a master solution of~$\MAS_{D'}$.
\end{proof}

We can now establish the correctness and finite termination of Algorithm~\ref{alg:slabdecomp}. One quantity plays a central role here: since the epsilon-constraint method is run with a fixed step size~$\delta>0$, the lower bounds~$\HV_{\mathrm{LB}}^s(\tilde x)$ computed in Line~\ref{algline:findHV} of Algorithm~\ref{alg:slabdecomp} may underestimate the true hypervolumes. In contrast to the discretization error of the slab decomposition, this error is not reduced during the course of the algorithm. It is, however, uniformly bounded: in the following, let~$\Delta\ge 0$ be a constant such that
\begin{equation}\label{eq:Delta}
	\HV^s(x) - \HV_{\mathrm{LB}}^s(x) \;\le\; \Delta \qquad \text{for all } s\in S \text{ and } x\in\X.
\end{equation}
By Lemma~\ref{lemma:hv-eps-constr-error}, the choice~$\Delta \coloneqq \delta\cdot\refpt_2$ always satisfies~\eqref{eq:Delta} (note that removing the points that do not dominate the reference point when evaluating Line~\ref{algline:findHV}, as described in Section~\ref{sec:outline}, does not change the hypervolume of the computed subset, since such points contribute no area in Definition~\ref{def:HV}).
If the epsilon-constraint method returns the complete nondominated set, then~$\Delta=0$ is a valid choice. By Lemma~\ref{lemma:hv-eps-constr-error}, this is the case whenever the nondominated sets~$\PND^s(x)$ are finite and~$\delta$ is at most the smallest first-component distance between two distinct nondominated points---in particular, whenever the~$f_1^s$ attain only discrete values and~$\delta$ is chosen accordingly (for integer-valued~$f_1^s$, $\delta\le1$).

The error bound~$\Delta$ limits the optimality gaps that the algorithm can certify. Consequently, our finite-termination guarantee requires the target gap~$\TG$ to be compatible with~$\Delta$, as made precise by Condition~\eqref{eq:termination-condition} in the following proposition.

\begin{proposition}\label{prop:finite-termination}
	Let the scenario set~$S$ be finite and let~$\Delta$ be as in~\eqref{eq:Delta}.
	\begin{enumerate}
		\item[(a)] If Algorithm~\ref{alg:slabdecomp} terminates with a first-stage decision~$\tilde x$, then
		      \begin{equation*}
			      \min_{s\in S}\HV^s(\tilde x) \;\ge\; \frac{\OPT}{1+\TG}.
		      \end{equation*}
		      Thus, $\tilde x$ solves Problem~\eqref{eq:hv-maxmin-problem} with an optimality gap of at most the target gap~$\TG$. 
              If, instead, Algorithm~\ref{alg:slabdecomp} terminates in Line~\ref{algline:emmax-zero} or in Line~\ref{algline:mas-infeasible}, then~$\OPT=0$, and every first-stage decision~$x\in\X$ is optimal for Problem~\eqref{eq:hv-maxmin-problem}. 
		\item[(b)] If
		      \begin{equation}\label{eq:termination-condition}
			      \Delta \;<\; \frac{\TG}{1+\TG}\cdot\OPT,
		      \end{equation}
		      then Algorithm~\ref{alg:slabdecomp} terminates after finitely many iterations and returns a first-stage decision satisfying~(a). (Note that Condition~\eqref{eq:termination-condition} implicitly requires~$\OPT>0$ and~$\TG>0$.)
	\end{enumerate}
\end{proposition}
\begin{proof}
	We first show~(a). The algorithm returns~$\tilde x$ only if $\mathrm{gap}\le\TG$ holds in Line~\ref{algline:generalterminationcrit}, in which case $\min_{s\in S}\HV_{\mathrm{LB}}^s(\tilde x)>0$ and, by Equation~\eqref{eq:gap},
	\begin{equation*}
		\tilde z + \EMmax \;\le\; (1+\TG)\cdot \min_{s\in S}\HV_{\mathrm{LB}}^s(\tilde x).
	\end{equation*}
	Together with Theorem~\ref{thm:bounds}, this yields
	\begin{equation*}
		\min_{s\in S}\HV^s(\tilde x) \;\ge\; \min_{s\in S}\HV_{\mathrm{LB}}^s(\tilde x) \;\ge\; \frac{\tilde z + \EMmax}{1+\TG} \;\ge\; \frac{\OPT}{1+\TG}.
	\end{equation*}

    Now suppose that the algorithm terminates in Line~\ref{algline:emmax-zero} or in Line~\ref{algline:mas-infeasible}. In the first case, we have~$\EMmax=0$, so Equation~\eqref{eq:EMmax} and~$w>0$ imply that no scenario~$s\in S$ admits a solution~$(x,y)\in\X\times\Y^s$ satisfying $f_1^s(x,y)\le\refpt_1$ and $f_2^s(x,y)<\refpt_2$. Hence, for every~$s\in S$, the set~$f^s(\X,\Y^s)$ does not intersect the box~$[0,\refpt_1)\times[0,\refpt_2)$, which, by Definition~\ref{def:HV}, yields $\HV^s(x)=0$ for every~$s\in S$ and~$x\in\X$.
	In the second case, $\MAS_D$ is infeasible for the current scenario subset~$D$. Suppose, for the sake of a contradiction, that there exists~$\bar x\in\X$ with $\HV^s(\bar x)>0$ for every scenario~$s\in D$. A positive hypervolume requires, for each~$s\in D$, a second-stage decision~$\bar y^s\in\Y^s(\bar x)$ with $f_1^s(\bar x,\bar y^s)<\refpt_1$ and $f_2^s(\bar x,\bar y^s)<\refpt_2$. Setting $x=\bar x$, $y^{s,i}=\bar y^s$, $a^{s,i}=0$, and~$h^{s,i}=0$ for all~$s\in D$ and~$i\in G^s$, as well as~$z=0$, then yields a feasible solution of~$\MAS_D$ --- a contradiction. Hence, for every~$x\in\X$, there is a scenario~$s\in D$ with~$\HV^s(x)=0$.
	In both cases, $\min_{s\in S}\HV^s(x)=0$ for every~$x\in\X$. Since~$\X\neq\emptyset$, this yields~$\OPT=0$ and, in particular, that every first-stage decision~$x\in\X$ is optimal for Problem~\eqref{eq:hv-maxmin-problem}.

    \medskip

	We now show~(b), so assume that Condition~\eqref{eq:termination-condition} holds, which implies~$\TG>0$ and~$\OPT>\Delta\ge0$. Suppose, for the sake of contradiction, that the algorithm does not terminate. We first collect three observations.

    \smallskip
    First, since~$\OPT>0$, there exists~$\bar x\in\X$ such that, for every scenario~$s\in S$, there is a second-stage decision $\bar y^s\in\Y^s(\bar x)$ with $f_1^s(\bar x,\bar y^s)<\refpt_1$ and $f_2^s(\bar x,\bar y^s)<\refpt_2$; indeed, a positive hypervolume
    requires the existence of a feasible point in the box~$[0,\refpt_1)\times[0,\refpt_2)$.
    Hence, the master problem~$\MAS_D$ is feasible in every iteration. Therefore, every iteration ends by either adding a scenario in Line~\ref{algline:addscenario} or refining the slab width in Line~\ref{algline:refine}. Whenever the check in Line~\ref{algline:branchcondition} fails, there is at least one scenario in~$S\setminus D$ that can be added. Moreover, the same solutions~$(\bar x,\bar y^s)$ imply~$\EMmax>0$ in every iteration by Equation~\eqref{eq:EMmax}. Since refinements multiply~$\EMmax$ by a positive factor, the factor~$\alpha$ in Equation~\eqref{eq:alpha} is therefore always well-defined.

    \smallskip
	Second, every scenario-addition step strictly enlarges~$D\subseteq S$. Since~$S$ is finite and the algorithm does not terminate, there is hence an iteration~$i_0$ such that the slab width is refined in every iteration~$i\ge i_0$. Denoting by~$\EMmax^i$ the value of~$\EMmax$ in iteration~$i$, we obtain $\EMmax^i \le \sigma^{i-i_0}\cdot\EMmax^{i_0}$ for all~$i\ge i_0$, since each refinement multiplies~$\EMmax$ by a factor of less than~$\sigma<1$ (see Section~\ref{sec:gapeval}). In particular, $\EMmax^i\to0$ for~$i\to\infty$.

    \smallskip
    \noindent
	Third, we claim that, in every iteration~$i\ge i_0$ with master solution~$(\tilde x,\tilde z)$,
	\begin{equation}\tag{$\ast$}\label{eq:LB-claim}
		\min_{s\in S}\HV_{\mathrm{LB}}^s(\tilde x) \;\ge\; \tilde z - \Delta.
	\end{equation}
	Indeed, since iteration~$i$ refines the slab width, the check in Line~\ref{algline:branchcondition} was successful, i.e., $\AA_{\text{LB}}^s(\tilde x)\ge\tilde z$ holds for all~$s\in S\setminus D$. For~$s\in S\setminus D$, we thus have
	\begin{equation*}
		\HV_{\mathrm{LB}}^s(\tilde x) \;\ge\; \AA_{\text{LB}}^s(\tilde x) \;\ge\; \tilde z \;\ge\; \tilde z-\Delta,
	\end{equation*}
	where the first inequality follows by the same geometric argument as in the proof of Lemma~\ref{lemma:AVleHV}: both quantities are computed from the same subset of the nondominated set, and every activated vertical slab portion is contained in the region dominated by the point that activates it (intersected with~$[0,\refpt_1]\times[0,\refpt_2]$). For~$s\in D$, Equation~\eqref{eq:Delta}, Lemma~\ref{lemma:AVleHV}, and Lemma~\ref{lemma:AA-master-feas} yield
	\begin{equation*}
		\HV_{\mathrm{LB}}^s(\tilde x) \;\ge\; \HV^s(\tilde x)-\Delta \;\ge\; \AA^s(\tilde x)-\Delta \;\ge\; \tilde z-\Delta.
	\end{equation*}
	This establishes~\eqref{eq:LB-claim} and thus proves the third observation.

    \smallskip\noindent
	Now, by Condition~\eqref{eq:termination-condition},
	\begin{equation*}
		\varepsilon^* \coloneqq \frac{\TG}{1+\TG}\cdot\OPT - \Delta \;>\; 0,
	\end{equation*}
	and, by the second observation, there is an iteration~$i\ge i_0$ with~$\EMmax^i\le\varepsilon^*$. Let~$(\tilde x,\tilde z)$ be the master solution of this iteration. Theorem~\ref{thm:bounds} yields
	\begin{equation*}
		\tilde z \;\ge\; \OPT - \EMmax^i \;\ge\; \OPT-\varepsilon^* \;=\; \frac{\OPT}{1+\TG} + \Delta
	\end{equation*}
	and, in particular, $\tilde z - \Delta \ge \frac{\OPT}{1+\TG} > 0$. Together with~\eqref{eq:LB-claim}, we obtain
	\begin{equation*}
		\mathrm{gap} \;=\; \frac{\tilde z+\EMmax^i}{\min_{s\in S}\HV_{\mathrm{LB}}^s(\tilde x)} - 1
		\;\le\; \frac{\tilde z+\EMmax^i}{\tilde z-\Delta} - 1
		\;=\; \frac{\EMmax^i+\Delta}{\tilde z-\Delta}
		\;\le\; \frac{\frac{\TG}{1+\TG}\cdot\OPT}{\frac{\OPT}{1+\TG}}
		\;=\;\TG.
	\end{equation*}
	Hence, the algorithm terminates in Line~\ref{algline:generalterminationcrit} of iteration~$i$---a contradiction. Therefore, the algorithm terminates after finitely many iterations and, since~$\MAS_D$ is always feasible, it returns a first-stage decision satisfying~(a).
\end{proof}

\noindent
We conclude this section with some remarks on the finite-termination
result and its assumptions.

First, Condition~\eqref{eq:termination-condition} cannot in general be dispensed with. Suppose that the underestimation of the epsilon-constraint method is fully realized in every iteration, i.e.,
\[
    \min_{s\in S}\HV_{\mathrm{LB}}^s(\tilde x)
    =
    \min_{s\in S}\HV^s(\tilde x)-\Delta.
\]
If $\min_{s\in S}\HV_{\mathrm{LB}}^s(\tilde x)>0$, then the assumed equality and $\min_{s\in S}\HV^s(\tilde x)\le\OPT$ imply~$\OPT>\Delta$. Hence, using $\tilde z+\EMmax\ge\OPT$ (Theorem~\ref{thm:bounds}), the computed gap satisfies
\[
    \mathrm{gap}
    \ge
    \frac{\OPT}{\OPT-\Delta}-1
    =
    \frac{\Delta}{\OPT-\Delta}.
\]
This lower bound is strictly larger than~$\TG$ whenever $\Delta>\frac{\TG}{1+\TG}\cdot\OPT$, that is, whenever the strict opposite of Condition~\eqref{eq:termination-condition} holds. If instead $\min_{s\in S}\HV_{\mathrm{LB}}^s(\tilde x)=0$, the computed gap is~$+\infty$, irrespective of~$\Delta$. Thus, whenever the strict opposite of Condition~\eqref{eq:termination-condition} holds, the target gap cannot be reached in either case, no matter how far the slab decomposition is refined.

Second, if~$\OPT=0$, then $\min_{s\in S}\HV_{\mathrm{LB}}^s(\tilde x) \le \min_{s\in S}\HV^s(\tilde x) \le \OPT = 0$ in every iteration. Hence, the denominator in Equation~\eqref{eq:gap} is zero and the gap is~$+\infty$ throughout. Consequently, the gap-based termination criterion cannot be satisfied for any target gap~$\TG\in[0,1)$. This reflects the fact that a relative optimality gap is meaningless if the optimal objective value is zero. The checks in Lines~\ref{algline:emmax-zero} and~\ref{algline:mas-infeasible} of Algorithm~\ref{alg:slabdecomp} detect two special cases in which~$\OPT=0$ can be certified without evaluating a gap, but they do not cover all such cases. In the remaining cases, monitoring the upper bound~$\tilde z + \EMmax$, which converges to zero, can provide numerical evidence for~$\OPT=0$.

Third, since~$\OPT$ is not known a priori, Condition~\eqref{eq:termination-condition} cannot be verified directly before running the algorithm. It can, however, be certified during the run: by Theorem~\ref{thm:bounds}, $\OPT \ge \min_{s\in S}\HV_{\mathrm{LB}}^s(\tilde x)$, so it suffices that $\Delta < \frac{\TG}{1+\TG}\cdot\min_{s\in S}\HV_{\mathrm{LB}}^s(\tilde x)$ holds in some iteration. Note that $\min_{s\in S}\HV_{\mathrm{LB}}^s(\tilde x)$ is evaluated in each iteration of the algorithm in Line~\ref{algline:findHV}. In practice, this suggests choosing the step size~$\delta$ of the epsilon-constraint method small compared to the magnitude of the hypervolumes involved (or, in cases where exact evaluation can be guaranteed, choosing~$\delta$
accordingly). Condition~\eqref{eq:termination-condition} also provides the following theoretical insight: whenever $\OPT>0$, any target gap~$\TG>0$ can be paired with a sufficiently small step size~$\delta$ so that~$\Delta$ is small enough for Condition~\eqref{eq:termination-condition} to hold. The algorithm is then guaranteed to terminate after finitely many iterations. 
In some cases, finite termination can even be verified before running the algorithm by deriving a positive lower bound on~$\OPT$.
Suppose there is one first-stage decision that, in every scenario, admits a feasible second-stage decision whose corresponding point is strictly smaller than the reference point in both components. Since~$S$ is finite, there is a positive minimum distance to the reference point in each dimension. Hence, the rectangles bounded by these points and the reference point have a positive minimum area over all scenarios, which guarantees a positive worst-case hypervolume. 
If the evaluation step computes the hypervolume exactly, for instance because the epsilon-constraint method returns the complete nondominated set, this is already sufficient for finite termination for every positive target gap. If the evaluation is inexact, its error must additionally be sufficiently small relative to this guaranteed hypervolume.

\begin{corollary}\label{cor:finite-termination-simple}
    Let~$\TG\in(0,1)$ and suppose there exists a first-stage decision~$\bar x\in\X$ such that, for every scenario~$s\in S$, there is a second-stage decision $\bar y^s\in\Y^s(\bar x)$ with
    \[
        f_1^s(\bar x,\bar y^s)<\refpt_1
        \qquad\text{and}\qquad
        f_2^s(\bar x,\bar y^s)<\refpt_2.
    \]
    Then the following statements hold:
    \begin{enumerate}
        \item[(a)] If the hypervolume lower bounds computed in Line~\ref{algline:findHV} of Algorithm~\ref{alg:slabdecomp} are exact throughout the run of the algorithm, then Algorithm~\ref{alg:slabdecomp} terminates after finitely many iterations and returns a first-stage decision with an optimality gap of at most~$\TG$. By Lemma~\ref{lemma:hv-eps-constr-error}, this is the case whenever every nondominated set encountered in the evaluation step is finite and~$\delta$ is at most its smallest first-component distance between two distinct nondominated points. In particular, this holds if there exists~$d>0$ such that, for every~$s\in S$, distinct values attained by~$f_1^s$ differ by at least~$d$ and~$\delta\le d$; if all~$f_1^s$, $s\in S$, are integer-valued, it suffices that $\delta\le1$.

        \item[(b)] More generally, define $\eta_\ell\coloneqq\min_{s\in S}\bigl\{\refpt_\ell-f_\ell^s(\bar x,\bar y^s)\bigr\}>0$ for $\ell\in\{1,2\}$.
        If
        \[
            \delta\cdot\refpt_2
            <
            \frac{\TG}{1+\TG}\cdot\eta_1\eta_2,
        \]
        then the same conclusion holds.
    \end{enumerate}
\end{corollary}

\begin{proof}
    The assumption of the corollary implies~$\OPT>0$. Indeed, by
    external stability, for every~$s\in S$ the feasible point
    $f^s(\bar x,\bar y^s)$ is weakly dominated by a nondominated point.
    Hence, $\HV^s(\bar x)>0$ for every~$s\in S$, and, since~$S$ is
    finite, $\OPT\ge\min_{s\in S}\HV^s(\bar x)>0$.

    For~(a), suppose that Algorithm~\ref{alg:slabdecomp} does not
    terminate. The arguments establishing the first two observations
    in the proof of Proposition~\ref{prop:finite-termination}~(b)
    apply: after finitely many scenario additions, every subsequent
    iteration is a refinement and~$\EMmax$ converges to zero. At every
    such iteration, the argument establishing the third observation in
    that proof, using exactness of the hypervolume lower bounds at the
    current~$\tilde x$ in place of Equation~\eqref{eq:Delta}, yields
    $\min_{s\in S}\HV_{\mathrm{LB}}^s(\tilde x)\ge\tilde z$.
    Hence, eventually
    $\EMmax\le\frac{\TG}{1+\TG}\OPT$. By
    Theorem~\ref{thm:bounds}, we then have
    $\tilde z\ge\frac{\OPT}{1+\TG}>0$, and Equation~\eqref{eq:gap}
    gives $\mathrm{gap}\le\EMmax/\tilde z\le\TG$, a contradiction.
    The optimality guarantee follows from
    Proposition~\ref{prop:finite-termination}~(a).

    The sufficient conditions for exactness stated in~(a) follow directly from Lemma~\ref{lemma:hv-eps-constr-error}.

    For~(b), the same domination argument gives
    $\HV^s(\bar x)\ge\eta_1\eta_2$ for every~$s\in S$, and thus
    $\OPT\ge\eta_1\eta_2$. Moreover,
    Lemma~\ref{lemma:hv-eps-constr-error} shows that
    $\Delta=\delta\cdot\refpt_2$ is a valid choice in
    Equation~\eqref{eq:Delta}. The assumed inequality therefore implies
    $\Delta<\frac{\TG}{1+\TG}\OPT$, so
    Condition~\eqref{eq:termination-condition} holds. The claim follows
    from Proposition~\ref{prop:finite-termination}~(b).
\end{proof}

\section{Computational experiments}\label{sec:experiments}\label{sec:comp}

In this section, we evaluate the slab decomposition algorithm (Algorithm~\ref{alg:slabdecomp}) on two applications, namely a two-stage robust selection problem introduced in Section~\ref{sec:sel-problem} and a two-stage robust assignment problem in Section~\ref{sec:ass-problem}. Most of our experiments and analyses are carried out on the selection problem (Section~\ref{sec:comp-sel}), and the assignment problem (Section~\ref{sec:comp-ass-half}) then serves to confirm the main effects on a structurally different and computationally harder problem. Our experiments address five main questions: (i) how the running time and the number of iterations scale with the number~$|S|$ of scenarios, the instance size, the target gap~$\TG$, and the difference between the first component of the reference point and the smallest attainable value of the first objective (Section~\ref{sec:comp-runtime}); (ii) how the starting slab width~$w$ should be chosen and whether a good choice can be found a priori (Section~\ref{sec:comp-slabwidth}); (iii) how the orientations of the slab decomposition (horizontal versus vertical) compare (Section~\ref{sec:comp-orientation}); 
(iv) how the choice of the reference point influences the running time and the returned first-stage decision (Section~\ref{sec:comp-refpoint}); and (v) how often the refinement factor~$\alpha$ from Equation~\eqref{eq:alpha} avoids repeated refinements, and how the safety factor~$\sigma<1$ affects the number of iterations and the running time (Section~\ref{sec:comp-alpha}).

\subsection{Application problems and instances}\label{sec:applications}

In the following subsections, we describe the applications on which we test the algorithm, and the instance generator shared by both applications.

\subsubsection{Selection problem}\label{sec:sel-problem}

For~$n\in\N$ and~$p\in\{0,\dots,n\}$, the two-stage robust selection problem~\citep{kasperski2017robust,goerigk2021flexible} is to select exactly~$p$ out of~$n$ items over two stages. An item can be selected in the first stage, selected in the second stage, or not selected at
all. For each item, its first-stage costs with respect to the two objectives are given by the vectors~$C,D\in\R_{\ge0}^n$. If an item is selected in the second stage, its costs depend on the realized scenario~$s\in S$ and are given by the vectors~$c^s,d^s\in\R_{\ge0}^n$ for the first and second objective, respectively.

The numbers of items selected in the first and second stage must sum to~$p$, and under each objective, the total cost is the sum of the corresponding first- and second-stage costs of the selected items. The problem is known to be NP-hard~\citep{kasperski2017robust}. Using binary vectors~$x$ and~$y$ to indicate the items selected in the first and second stage, respectively, we obtain the following feasible sets and objective functions: 
\begin{align*}
    \X
        &= \{0,1\}^n, \\
    \Y^s(x)
        &= \left\{y\in\{0,1\}^n:
            \sum_{i\in[n]}(x_i+y_i)=p,\;
            x_i+y_i\le 1 \quad \forall i\in[n]
        \right\}, \\
    f_1^s(x,y)
        &= \sum_{i\in[n]}\left(C_i x_i+c_i^s y_i\right), \\
    f_2^s(x,y)
        &= \sum_{i\in[n]}\left(D_i x_i+d_i^s y_i\right).
\end{align*}

For the initialization of the scenario subset in Algorithm~\ref{alg:slabdecomp}, we use a problem-specific rule to select a presumably worst scenario based on the instance data (see Section~\ref{sec:outline}). To do so, $\Init()$ is chosen such that it returns a scenario~$s$ that maximizes the sum of the~$p$ smallest entries of~$c^s$ plus the sum of the~$p$ smallest entries of~$d^s$.

\subsubsection{Assignment problem}\label{sec:ass-problem}
For~$n\in\N$ and~$p\in\{0,\dots,n\}$, we consider a bi-objective two-stage robust variant of the assignment problem, combining the bi-objective setting of~\citet{przybylski2008two} with a two-stage decision structure similar to that of the selection problem in Section~\ref{sec:sel-problem} and conceptually related to recoverable robust assignment~\citep{fischer2020investigation}. The problem is defined on a complete bipartite graph with~$n$ vertices on each side and edge set~$E$, so that~$|E|=n^2$.
A matching is a set of edges no two of which share a common endpoint,
and a perfect matching is a matching that matches every vertex. The goal is to construct a perfect matching over two stages. In the first stage, at least~$p$ edges are selected and must form a matching; after the scenario~$s\in S$ is realized, the remaining vertices are matched in the second stage. For each edge, its first-stage costs with respect to the two objectives are given by the vectors~$C,D\in\R_{\ge0}^E$. If an edge is selected in the second stage, its costs depend on the realized scenario~$s\in S$ and are given by the vectors~$c^s,d^s\in\R_{\ge0}^E$ for the first and second objective, respectively.

The union of the edges selected in the first and second stages must form a perfect matching. Under each objective, the total cost is the sum of the first-stage costs of the edges selected in the first stage and the scenario-dependent second-stage costs of those selected in the second stage. For any binary vector~$z\in\{0,1\}^E$, let $E(z)=\{e\in E:z_e=1\}$. Accordingly, using binary vectors~$x$ and~$y$ to indicate the edges selected in the first and second stage, respectively, the feasible sets and objective functions are given as follows:
\begin{align*}
    \X
        &= \left\{x\in\{0,1\}^E:
            |E(x)|\ge p,\;
            E(x)\text{ is a matching}
        \right\}, \\
    \Y^s(x)
        &= \left\{y\in\{0,1\}^E:
            E(x)\cap E(y)=\emptyset,\;
            E(x)\cup E(y)\text{ is a perfect matching}
        \right\}, \\
    f_1^s(x,y)
        &= \sum_{e\in E}\left(C_e x_e+c_e^s y_e\right), \\
    f_2^s(x,y)
        &= \sum_{e\in E}\left(D_e x_e+d_e^s y_e\right).
\end{align*}

To ensure that~$E(x)$ is a matching, we simply require, for each
vertex, that it is incident to at most one selected edge. For a perfect
matching, each vertex must be incident to exactly one selected edge.

For this problem, $\Init()$ returns a scenario~$s$ that maximizes the total second-stage weight $\sum_{e\in E}(c_e^s+d_e^s)$ or, equivalently, the one with the largest average second-stage edge weight.

\subsubsection{Instance generation}\label{sec:sel-instances}

Both application problems use the same instance generator. Each instance is built from two integer weight vectors~$C$ and~$D$, which contain the first-stage weights of the first and the second objective, respectively. They have~$n$ entries for the selection problem and~$n^2$ entries for the assignment problem. Their entries are drawn independently from a normal distribution with mean~$60$ and standard deviation~$20$, truncated to the interval~$[20,100]$ and rounded to the nearest integer.

We index the scenarios by~$s=1,\dots,|S|$.
The scenario-dependent second-stage weight vectors~$c^s$ and~$d^s$ are obtained by perturbing the first-stage weights: each is scaled by a scenario-wide multiplier~$m_s$ and overlaid with independent Gaussian noise. More precisely, the $i$-th entries of~$c^s$ and~$d^s$ are
\[
    c^s_i=\max\bigl\{0,\; \round{m_s\cdot C_i+20\,\xi}\bigr\},
    \qquad
    d^s_i=\max\bigl\{0,\; \round{m_s\cdot D_i+20\,\xi'}\bigr\},
\]
where~$\round{\cdot}$ denotes rounding to the nearest integer and the perturbations~$\xi,\xi'\sim\mathcal{N}(0,1)$ are independent across entries, objectives, and scenarios. The multipliers~$m_s$ are drawn independently from a normal distribution with mean~$1$ and standard deviation~$0.1$, truncated to~$[0.8,1.2]$, so that a scenario is systematically up to $20\%$ more or less costly than the nominal first stage. Deriving the second stage from the first in this way makes the two stages correlated (an item or edge that is expensive to select now tends to be expensive later as well), while the noise still lets the two stages disagree on their individual entries.

The random numbers are generated using independent streams derived from the instance seed: one for the first-stage weights, sampled entry by entry; one for the scenario multipliers, sampled scenario by scenario; and a separate stream for the Gaussian perturbations of each scenario, again sampled entry by entry. Consequently, for a fixed seed, the first- and second-stage weight vectors of an instance with a smaller value of~$n$ are prefixes of the corresponding vectors of an instance with a larger value of~$n$. Similarly, increasing~$|S|$ leaves all previously generated scenarios unchanged. A specific instance is thus determined by~$n$, $p$, the number~$|S|$ of scenarios, and the instance seed used for the random number generators.

\subsection{Experimental setup}\label{sec:comp-setup}

All experiments were performed on a Linux system with an AMD EPYC~7542 processor ($32$~cores, $64$~threads, base frequency $2.9$\,GHz). The algorithm is implemented in Python~3.13.9, using Gurobi~13.0.2
as the MIP solver.
Gurobi is restricted to a single thread in all experiments to ensure consistent performance and comparable resource usage. 

For both application problems, the epsilon-constraint method uses a step size of~$\delta=10^{-2}$. Since the objective values are integer-valued in the unscaled instances and the scaling factors used in Section~\ref{sec:comp-orientation} are all at least~$1$, this step size is sufficiently small for the epsilon-constraint method to return the complete nondominated set in all experiments. 
The lower bound on the refinement factor is fixed at~$\alpha_{\min}=0.1$, and each run is subject to a time limit of three hours. The implementation of all algorithms and models, together with the instance generator used in the computational experiments, is available on~\zenodolink.

Unless stated otherwise, the reference point is computed using the combinatorial routines described in Appendices~\ref{sec-app:refpoint-selection} and~\ref{sec-app:refpoint-assignment}.

For each instance and for every scenario~$s$, we additionally record the distance
\[
    \Rs\coloneqq\refpt_1-\underline f_1^s
\]
between the first component of the reference point and the smallest value~$\underline f_1^s$ of the first objective, as computed for the construction of the slab set~$G^s$ in Section~\ref{sec:model}. By this construction, scenario~$s$ contributes exactly $|G^s|=\lfloor\Rs/w\rfloor$ slabs of width~$w$ to every master problem that contains it (and none if the problem defining~$\underline f_1^s$ is infeasible), so that~$\Rs$ determines the number of slabs, and with it the size of the master problem, for every slab width. Computing~$\Rs$ requires solving only a single MIP, which takes a few milliseconds per scenario on our instances. In addition, this MIP can be reused when constructing the slab sets~$G^s$ for the master problems (for horizontal slabs, the second objective takes the role of the first).

For both application problems, we additionally verified that the hypothesis of Corollary~\ref{cor:finite-termination-simple} holds on every instance in our test set. To do so, we solved a small auxiliary feasibility MIP that searches for a first-stage decision~$x\in\X$ and second-stage decisions~$y^s\in\Y^s(x)$ for the scenarios~$s\in S$ such that
\[
    f_1^s(x,y^s)\le\refpt_1-\tfrac12
    \qquad\text{and}\qquad
    f_2^s(x,y^s)\le\refpt_2-\tfrac12
    \qquad\text{for all }s\in S.
\]
These constraints imply that every point~$f^s(x,y^s)$ lies strictly below the reference point in both components, as required by Corollary~\ref{cor:finite-termination-simple}. The auxiliary MIP was feasible for every 
instance of both problems, each time with respect to the reference point used in the corresponding runs.
Since the epsilon-constraint method returns the complete nondominated set in our experiments, the hypervolume lower bounds used in the evaluation step are exact. Hence, part~(a) of Corollary~\ref{cor:finite-termination-simple} guarantees finite termination a priori on the entire test set.

In addition, we checked whether finite termination could also have been certified during the runs themselves, without solving the auxiliary MIP. For the same reason as above, we may take~$\Delta=0$. Since all target gaps considered are positive, Proposition~\ref{prop:finite-termination} then shows that~$\OPT>0$ is sufficient for finite termination. By Theorem~\ref{thm:bounds}, the positivity of~$\OPT$ can be certified during a run as soon as some iteration satisfies
\[
    \min_{s\in S}\HV_{\mathrm{LB}}^s(\tilde x)>0.
\]
Evaluating this condition on the log files, which record the per-scenario hypervolume lower bounds of every evaluation step, we found that finite termination was certified in this way in $4264$ of the $4320$ selection runs and in $2298$ of the $2430$ assignment runs reported in Section~\ref{sec:comp-ass-half}---and in each of them already at the \emph{first} evaluation step. The remaining $56$ selection and $132$ assignment runs are timeouts that never finished their first master problem and hence never evaluated a first-stage decision at all.

\subsection{Results for the selection problem}\label{sec:comp-sel}
% Figures and tables of this subsection are produced by
% python-code/auswertung-server-SEL.ipynb from the server run
% auswertung-server/SEL/run_20260913_231058_new_default (run_experiments_SEL.py).

In this section, we report the computational results for the selection problem. Each of Sections~\ref{sec:comp-runtime}--\ref{sec:comp-alpha} concludes with an answer to the corresponding question posed at the beginning of Section~\ref{sec:comp}.

\FloatBarrier
\subsubsection{Experimental design}\label{sec:comp-design}

We consider the full cross product of the following instance and algorithm parameters: target gaps $\TG\in\{5\%,10\%\}$, instance sizes $n\in\{10,20,40,60\}$ with three values of~$p$ per~$n$ ($p\in\{2,4,5\}$ for $n=10$, $p\in\{2,6,10\}$ for $n=20$, $p\in\{2,11,20\}$ for $n=40$, and $p\in\{2,16,30\}$ for $n=60$), scenario-set sizes $|S|\in\{4,16,128\}$, five random seeds for scenario generation, starting slab widths $w\in\{10,40,100\}$, and safety factors $\sigma\in\{0.5,0.8,0.95,0.99\}$ that scale the refinement factor~$\alpha$ from Equation~\eqref{eq:alpha-estimation} before each refinement (see Section~\ref{sec:comp-alpha}). This results in $4320$ runs in total.

All but $404$ of the $4320$ runs finished within the time limit of three hours. Because all preprocessing and setup steps are very fast compared to the total algorithm time (see Section~\ref{sec:comp-runtime}), \emph{all} timeouts occurred inside the loop of the slab decomposition algorithm (see Line~\ref{alg-line:general-while-loop} in Algorithm~\ref{alg:slabdecomp}). They are confined to five instance sizes---$208$ at $(n,p)=(60,30)$, $97$ at $(40,20)$, $80$ at $(60,16)$, $14$ at $(40,11)$, and $5$ at $(20,10)$. Out of these, $56$ runs (all at $n=60$) never finished the first master problem and therefore never completed an evaluation, $148$ timed out in the second master problem after one evaluation, and the remaining $200$ after two or more evaluations.

\FloatBarrier
\subsubsection{Running time and iteration counts}\label{sec:comp-runtime}

\begin{figure}[!htbp]
	\centering
	\includegraphics[width=\textwidth]{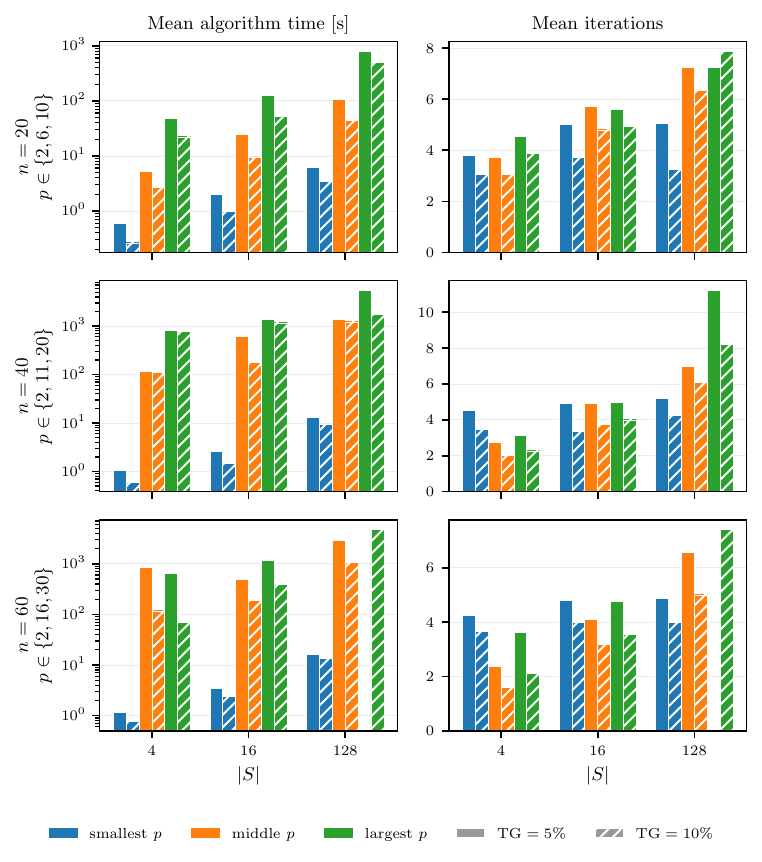}
	\caption{Mean algorithm time (left column, logarithmic scale) and mean number of iterations (right column) for the selection problem at the standard safety factor~$\sigma=0.95$, with one row each for $n=20$, $40$, and~$60$.
	Within each panel the bars are grouped by the number~$|S|$ of scenarios; within each group, the bar color distinguishes the smallest, middle, and largest value of~$p$ given in the row label, and the hatching distinguishes the target gap (solid: $5\%$, hatched: $10\%$). Each bar averages the solved runs over the five scenario seeds and the three starting slab widths. At $n=40$ and $n=60$, the bars for the two larger values of~$p$ omit the timed-out runs and are therefore subject to survivor bias.}
	\label{fig:sel-runtime}
\end{figure}

\begin{figure}[!htbp]
	\centering
	\includegraphics[width=\textwidth]{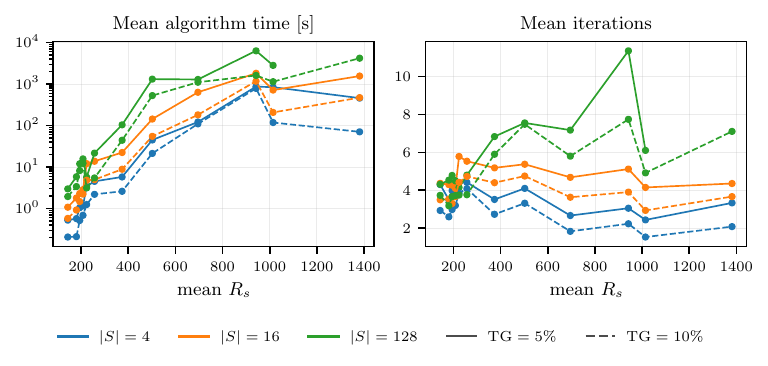}
	\caption{Mean algorithm time (left panel, logarithmic scale) and mean number of iterations (right panel) for the twelve tested instance sizes, plotted against the mean of~$\Rs$ over the scenarios of the instance (see Section~\ref{sec:comp-setup}). Colors distinguish the number~$|S|$ of scenarios and line styles the target gap; each point averages the solved runs over the scenario seeds, the starting slab widths, and the safety factors. On the four largest instance sizes the timed-out runs are missing from these averages.}
	\label{fig:sel-frontsize}
\end{figure}

Table~\ref{tab:sel-main} reports, for the default starting slab width $w=40$ (see Section~\ref{sec:comp-slabwidth}), the average time spent in the slab decomposition loop and the average number of iterations (i.e., the number of times the master problem is solved), both separately for the three safety factors $\sigma\in\{0.8,0.95,0.99\}$. Throughout this section, \emph{algorithm time} denotes everything except the computation of the reference point, i.e., the slab decomposition loop together with the computation of the initial value of~$\EMmax$ and the initial scenario subset~$D$. 
The time to compute the reference point is kept separate, since it is part of the user-specified input. The time to compute it with the methods described in Appendix~\ref{sec-app:refpoint-selection} is reported in a separate column. As can be seen, the computation of the reference point is very fast relative to the total running time and takes on average only about $2.2$\,ms and never more than $13.5$\,ms.
Algorithm times are averaged over all runs, with each timed-out run assigned a running time equal to the three-hour time limit,
while iteration counts are averaged over the runs solved within the time limit only. Figure~\ref{fig:sel-runtime} visualizes the scaling of the algorithm time and of the iteration count with~$|S|$ for the three larger values of~$n$, at the safety factor~$\sigma=0.95$, and Figure~\ref{fig:sel-frontsize} shows the same quantities as a function of the mean of~$\Rs$ over the scenarios of the instance.

\begin{sidewaystable}
\centering\footnotesize
\begin{tabular}{lll|r|rrr|rrr}
\hline
 & & & & \multicolumn{3}{c|}{$\TG=5\%$} & \multicolumn{3}{c}{$\TG=10\%$} \\
$n$ & $p$ & $|S|$ & RP [ms] & Solved & Algorithm [s] & Iterations & Solved & Algorithm [s] & Iterations \\
\hline
10 & 2 & 4 & 0.15 & 15/15 & 0.5\,/\,0.9\,/\,0.6 & 4.0\,/\,6.2\,/\,5.4 & 15/15 & 0.3\,/\,0.3\,/\,0.4 & 3.4\,/\,4.6\,/\,4.4 \\
 &  & 16 & 0.34 & 15/15 & 0.9\,/\,1.1\,/\,0.9 & 3.8\,/\,4.6\,/\,5.4 & 15/15 & 0.4\,/\,0.7\,/\,0.4 & 3.2\,/\,4.6\,/\,3.4 \\
 &  & 128 & 2.04 & 15/15 & 3.5\,/\,2.6\,/\,3.1 & 4.8\,/\,4.8\,/\,5.2 & 15/15 & 1.6\,/\,1.8\,/\,2.0 & 3.6\,/\,4.0\,/\,4.4 \\
 & 4 & 4 & 0.16 & 15/15 & 3.9\,/\,3.2\,/\,3.9 & 4.4\,/\,5.0\,/\,5.2 & 15/15 & 1.2\,/\,1.8\,/\,1.4 & 4.0\,/\,4.8\,/\,5.0 \\
 &  & 16 & 0.35 & 15/15 & 16.3\,/\,8.5\,/\,13.6 & 5.8\,/\,7.0\,/\,6.8 & 15/15 & 5.4\,/\,3.1\,/\,4.6 & 5.6\,/\,5.6\,/\,6.4 \\
 &  & 128 & 2.33 & 15/15 & 5.7\,/\,8.4\,/\,10.8 & 4.8\,/\,4.6\,/\,5.6 & 15/15 & 3.1\,/\,5.4\,/\,4.9 & 3.8\,/\,5.4\,/\,5.8 \\
 & 5 & 4 & 0.17 & 15/15 & 5.2\,/\,4.3\,/\,4.2 & 4.6\,/\,4.8\,/\,5.2 & 15/15 & 1.5\,/\,2.0\,/\,2.0 & 3.8\,/\,4.8\,/\,5.4 \\
 &  & 16 & 0.34 & 15/15 & 16.7\,/\,13.7\,/\,17.2 & 5.8\,/\,6.6\,/\,7.0 & 15/15 & 5.0\,/\,4.5\,/\,7.3 & 5.4\,/\,6.0\,/\,6.6 \\
 &  & 128 & 2.52 & 15/15 & 12.7\,/\,15.2\,/\,15.2 & 4.8\,/\,5.8\,/\,5.4 & 15/15 & 5.6\,/\,6.9\,/\,4.7 & 4.4\,/\,5.2\,/\,4.8 \\
\hline
20 & 2 & 4 & 0.18 & 15/15 & 0.6\,/\,0.6\,/\,1.0 & 3.8\,/\,4.2\,/\,4.8 & 15/15 & 0.3\,/\,0.4\,/\,0.3 & 3.2\,/\,4.0\,/\,3.8 \\
 &  & 16 & 0.47 & 15/15 & 1.5\,/\,2.0\,/\,2.3 & 4.2\,/\,5.6\,/\,5.8 & 15/15 & 1.1\,/\,1.2\,/\,1.4 & 4.0\,/\,4.6\,/\,5.8 \\
 &  & 128 & 3.07 & 15/15 & 7.4\,/\,7.8\,/\,6.1 & 5.2\,/\,5.4\,/\,5.2 & 15/15 & 3.9\,/\,4.2\,/\,6.2 & 3.6\,/\,4.0\,/\,5.8 \\
 & 6 & 4 & 0.19 & 15/15 & 4.0\,/\,2.7\,/\,4.8 & 3.2\,/\,3.2\,/\,4.8 & 15/15 & 1.5\,/\,1.5\,/\,2.2 & 2.6\,/\,3.0\,/\,3.4 \\
 &  & 16 & 0.51 & 15/15 & 21.3\,/\,22.1\,/\,29.2 & 5.4\,/\,6.2\,/\,6.8 & 15/15 & 6.3\,/\,9.2\,/\,7.3 & 4.2\,/\,5.2\,/\,5.2 \\
 &  & 128 & 3.57 & 15/15 & 147\,/\,99.6\,/\,158 & 8.0\,/\,8.2\,/\,9.0 & 15/15 & 36.9\,/\,46.2\,/\,33.5 & 6.0\,/\,7.2\,/\,7.0 \\
 & 10 & 4 & 0.24 & 15/15 & 19.6\,/\,39.1\,/\,27.6 & 3.8\,/\,5.2\,/\,5.6 & 15/15 & 5.6\,/\,5.5\,/\,6.7 & 3.0\,/\,3.6\,/\,3.6 \\
 &  & 16 & 0.64 & 15/15 & 101\,/\,122\,/\,174 & 5.2\,/\,6.2\,/\,7.6 & 15/15 & 25.3\,/\,30.9\,/\,37.0 & 4.4\,/\,5.2\,/\,5.4 \\
 &  & 128 & 4.36 & 13/15 & 2112\,/\,2824\,/\,2893 & 9.8\,/\,10.5\,/\,10.8 & 15/15 & 182\,/\,218\,/\,198 & 8.6\,/\,9.4\,/\,9.6 \\
\hline
40 & 2 & 4 & 0.26 & 15/15 & 1.1\,/\,1.3\,/\,1.2 & 4.2\,/\,5.8\,/\,5.0 & 15/15 & 0.4\,/\,0.7\,/\,0.7 & 2.8\,/\,3.8\,/\,4.4 \\
 &  & 16 & 0.74 & 15/15 & 2.5\,/\,2.9\,/\,2.6 & 4.4\,/\,5.4\,/\,5.0 & 15/15 & 1.6\,/\,1.8\,/\,2.7 & 3.4\,/\,3.8\,/\,5.4 \\
 &  & 128 & 4.93 & 15/15 & 14.9\,/\,15.3\,/\,20.6 & 5.6\,/\,6.0\,/\,8.0 & 15/15 & 10.7\,/\,12.7\,/\,13.8 & 4.8\,/\,5.6\,/\,6.2 \\
 & 11 & 4 & 0.31 & 15/15 & 9.6\,/\,9.4\,/\,14.6 & 2.2\,/\,2.6\,/\,3.0 & 15/15 & 2.8\,/\,2.8\,/\,2.5 & 1.2\,/\,1.2\,/\,1.2 \\
 &  & 16 & 0.92 & 15/15 & 152\,/\,297\,/\,396 & 4.6\,/\,5.2\,/\,6.0 & 15/15 & 33.4\,/\,35.5\,/\,33.8 & 2.8\,/\,2.8\,/\,2.8 \\
 &  & 128 & 7.20 & 12/15 & 3705\,/\,2703\,/\,2970 & 7.0\,/\,7.8\,/\,9.0 & 15/15 & 130\,/\,129\,/\,136 & 4.0\,/\,4.0\,/\,4.2 \\
 & 20 & 4 & 0.32 & 15/15 & 102\,/\,105\,/\,121 & 3.0\,/\,3.4\,/\,3.4 & 15/15 & 28.4\,/\,29.7\,/\,30.1 & 1.8\,/\,1.8\,/\,1.8 \\
 &  & 16 & 1.18 & 14/15 & 2881\,/\,600\,/\,1057 & 4.8\,/\,4.8\,/\,5.4 & 15/15 & 104\,/\,102\,/\,105 & 2.6\,/\,2.6\,/\,2.6 \\
 &  & 128 & 8.35 & 7/15 & 10536\,/\,7297\,/\,8053 & 11.0\,/\,10.0\,/\,10.3 & 12/15 & 2728\,/\,2733\,/\,2729 & 5.2\,/\,5.2\,/\,5.2 \\
\hline
60 & 2 & 4 & 0.35 & 15/15 & 1.0\,/\,1.3\,/\,1.6 & 3.6\,/\,4.6\,/\,5.6 & 15/15 & 0.6\,/\,0.8\,/\,1.0 & 2.8\,/\,3.8\,/\,4.4 \\
 &  & 16 & 0.98 & 15/15 & 2.9\,/\,4.1\,/\,4.1 & 4.0\,/\,5.6\,/\,6.2 & 15/15 & 2.5\,/\,2.6\,/\,3.0 & 4.0\,/\,4.2\,/\,5.0 \\
 &  & 128 & 7.39 & 15/15 & 18.7\,/\,21.2\,/\,25.1 & 5.2\,/\,6.2\,/\,7.4 & 15/15 & 16.6\,/\,18.3\,/\,20.1 & 5.2\,/\,5.4\,/\,6.2 \\
 & 16 & 4 & 0.42 & 15/15 & 24.8\,/\,25.0\,/\,24.2 & 1.6\,/\,1.6\,/\,1.6 & 15/15 & 22.5\,/\,23.1\,/\,22.3 & 1.4\,/\,1.4\,/\,1.4 \\
 &  & 16 & 1.25 & 14/15 & 2586\,/\,251\,/\,321 & 3.2\,/\,3.6\,/\,4.0 & 15/15 & 52.4\,/\,53.7\,/\,53.4 & 2.0\,/\,2.0\,/\,2.0 \\
 &  & 128 & 9.45 & 10/15 & 6479\,/\,4360\,/\,4863 & 5.3\,/\,6.2\,/\,5.7 & 15/15 & 225\,/\,223\,/\,223 & 2.4\,/\,2.4\,/\,2.4 \\
 & 30 & 4 & 0.51 & 12/15 & 2476\,/\,2428\,/\,2421 & 2.0\,/\,2.0\,/\,2.2 & 12/15 & 2235\,/\,2233\,/\,2230 & 1.5\,/\,1.5\,/\,1.5 \\
 &  & 16 & 1.44 & 12/15 & 3086\,/\,3078\,/\,3083 & 3.2\,/\,3.2\,/\,3.2 & 12/15 & 2564\,/\,2577\,/\,2577 & 2.5\,/\,2.5\,/\,2.5 \\
 &  & 128 & 12.37 & 0/15 & 10800\,/\,10800\,/\,10800 & -- & 9/15 & 5972\,/\,5975\,/\,5953 & 5.0\,/\,5.0\,/\,5.0 \\
\hline
\multicolumn{4}{l|}{Average $n=10$} & 100\% & 7.2\,/\,6.4\,/\,7.7 & 4.8\,/\,5.5\,/\,5.7 & 100\% & 2.7\,/\,3.0\,/\,3.1 & 4.1\,/\,5.0\,/\,5.1 \\
\multicolumn{4}{l|}{Average $n=20$} & 99\% & 268\,/\,347\,/\,366 & 5.4\,/\,6.0\,/\,6.6 & 100\% & 29.3\,/\,35.3\,/\,32.6 & 4.4\,/\,5.1\,/\,5.5 \\
\multicolumn{4}{l|}{Average $n=40$} & 91\% & 1934\,/\,1226\,/\,1404 & 4.6\,/\,5.4\,/\,5.9 & 98\% & 338\,/\,339\,/\,339 & 3.1\,/\,3.4\,/\,3.7 \\
\multicolumn{4}{l|}{Average $n=60$} & 80\% & 2830\,/\,2330\,/\,2394 & 3.5\,/\,4.2\,/\,4.5 & 91\% & 1232\,/\,1234\,/\,1231 & 2.9\,/\,3.1\,/\,3.4 \\
\hline
\multicolumn{4}{l|}{Overall average} & 92\% & 1260\,/\,977\,/\,1043 & 4.6\,/\,5.3\,/\,5.7 & 97\% & 400\,/\,403\,/\,402 & 3.7\,/\,4.2\,/\,4.5 \\
\hline
\end{tabular}
\caption{Results for the selection problem with starting slab width $w=40$. Each cell of the ``Algorithm'' and ``Iterations'' columns contains three values, one per safety factor $\sigma=0.8\,/\,0.95\,/\,0.99$; each value averages the 5 scenario seeds. ``Solved'' is the number of runs (out of $5$~seeds $\times$ $3$~values of~$\sigma$) that finished within the time limit of three hours, ``RP'' the time of the combinatorial reference-point computation, ``Algorithm'' the mean algorithm time (i.e.\ everything except the computation of the reference point) over the $5$~seeds. ``Iterations'' is the number of master solves, averaged over the solved runs only.}
\label{tab:sel-main}
\end{sidewaystable}

We make several observations. First, the algorithm benefits strongly from the scenario-addition scheme: even for $|S|=128$, on average only $3.8$ scenarios (at most~$14$) are ever added to the subset~$D$ of scenarios (including the initial one). The number of iterations grows only mildly with~$|S|$ (right column of Figure~\ref{fig:sel-runtime}), whereas the running time grows faster than the iteration count, since later iterations solve master problems over a larger set~$D$ and/or a finer slab decomposition.
The evaluation step, which applies the epsilon-constraint method to each scenario in~$S$, is cheap in comparison: over all solved runs, about $94\%$ of the algorithm time is spent solving the master problems; a detailed cost decomposition follows below.

Second, difficulty is governed by the size of the master problems, i.e., by the number of slabs the objective space is decomposed into. 
 We use~$\Rs$ in the following analysis, since it determines the number of slabs exactly. 
For a given slab width~$w$, each scenario contributes $\lfloor\Rs/w\rfloor$ slabs, together with their binary activation variables and second-stage copies, to every master problem in which it appears. Along the mean of~$\Rs$, the mean algorithm time spans almost four orders of magnitude, as can be seen in Figure~\ref{fig:sel-frontsize}. The growth comes from harder master problems with more slabs. 
Note, however, that~$\Rs$ determines the number of slabs included in the master problem for scenario~$s$, based on feasibility for \emph{some} first-stage decision. Since the master problem requires a common first-stage decision for all scenarios in~$D$, not all of these slabs can necessarily be activated with positive height for the chosen decision. The number of iterations, on the other hand, does not grow with~$\Rs$ in the same way: for a fixed slab width, a larger~$\Rs$ means that the algorithm starts with more slabs and, thus, a finer decomposition relative to the covered first-objective range, which may reduce the number of refinements needed to certify the target gap.

Third, we decompose the total time into its components: the reference-point computation, the slab decomposition setup (computation of~$\EMmax$ and of the initial set~$D$), the solves of the master problems, and the evaluation step (the epsilon-constraint computation of the nondominated sets of all scenarios together with the resulting hypervolume and activated-area values).
Over all solved runs, these account for roughly $0.01\%$, $0.1\%$, $92\%$, and $6\%$ of the total time, respectively (the remaining $1.9\%$ is instance generation and bookkeeping)---but the split depends strongly on the instance. For the easy $p=2$ instances, the evaluation step dominates (up to $48\%$ of the total time), since every iteration recomputes the nondominated sets of \emph{all} $|S|$ scenarios while the master problems stay tiny; for the harder instances, the master solves become dominant (up to $97\%$). 

Summing up, the running time is governed by the size of the master problems and therefore by the mean of~$\Rs$, along which it spans almost four orders of magnitude. The number of iterations, in contrast, depends much less on~$\Rs$ and grows only mildly with the number~$|S|$ of scenarios, since the scenario-addition scheme keeps the number of selected scenarios small. A looser target gap reduces both the running time and the number of iterations, mainly because fewer refinements are required.

\FloatBarrier
\subsubsection{Choice of the starting slab width}\label{sec:comp-slabwidth}

The starting slab width~$w$ controls a trade-off: a fine initial decomposition (small~$w$) makes refinements rare but yields master problems with many slabs from the first iteration on, whereas a coarse start keeps the early master problems small at the price of additional refinement iterations. For each run, we normalize its algorithm time by the time of the fastest starting slab width on the same instance and configuration, and report the geometric mean of this slowdown factor.

\begin{table}[!htbp]
	\centering
	\begin{tabular}{c|ccc}
		\hline
		& \multicolumn{3}{c}{Geometric mean slowdown} \\
		$w$ & $\TG=5\%$ & $\TG=10\%$ & all \\
		\hline
		$10$ & $1.75$ & $2.34$ & $2.03$ \\
		$40$ & $1.79$ & $1.85$ & $1.82$ \\
		$100$ & $1.99$ & $2.24$ & $2.12$ \\
		\hline
	\end{tabular}
	\caption{Influence of the starting slab width~$w$ for the selection problem: geometric mean, over the $1190$ configurations in which all three starting widths finished within the time limit, of the algorithm time relative to the fastest starting width for the same instance and configuration.}
	\label{tab:sel-slabwidth}
\end{table}

Table~\ref{tab:sel-slabwidth} summarizes the mean slowdowns and shows that, overall, $w=40$ is the best of the three. This makes $w=40$ our default in Table~\ref{tab:sel-main}, the more so as the marginal lead of $w=10$ at the small target gap is additionally inflated by survivorship (the hardest instances, on which $w=10$ times out, drop out of the paired comparison).
The very coarse start $w=100$ pays for its cheap early master problems with many additional refinement rounds ($2.81$ refinements on average at $\TG=5\%$, versus $2.00$ for $w=40$ and only $0.62$ for $w=10$), and the fine start $w=10$ pays for its almost refinement-free run with master problems that are large from the first iteration on.

\begin{figure}[!htbp]
	\centering
	\includegraphics[width=\textwidth]{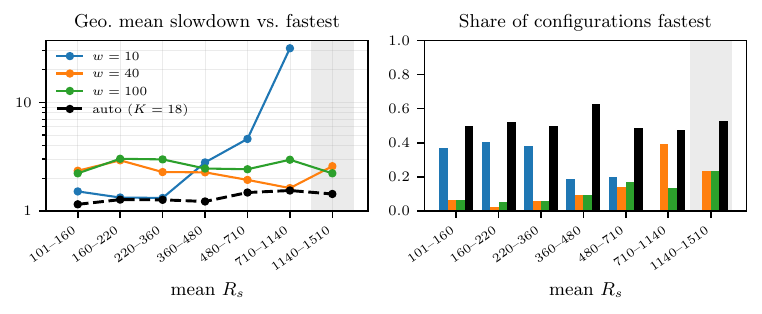}
	\caption{Starting slab width compared to~$\Rs$ for the selection problem ($\sigma=0.95$). Configurations are binned by the mean of~$\Rs$ over the scenarios~$s\in S$ of the instance. Besides the three fixed starting slab widths, the automatic width $w=\max\{\delta,\Rs[s_0]/18\}$ (``auto'', black) is shown. Left: geometric mean slowdown relative to the fastest of the four settings for the same configuration (logarithmic axis), over the configurations solved by all four. Right: share of configurations for which each setting is the fastest choice. In the shaded rightmost bin, no run with $w=10$ finished within the time limit; there, the other three settings are compared on the configurations solved by all of them.}
	\label{fig:sel-extent}
\end{figure}

A natural question is whether the best starting slab width can be predicted a priori from a property of the instance. Since~$\Rs$ determines the number of slabs, it is the natural candidate, and binning the configurations by the mean of~$\Rs$ indeed reveals a clear pattern (Figure~\ref{fig:sel-extent}, for the runs with $\sigma=0.95$): among the three fixed widths, the fine start $w=10$ is the best choice up to a mean~$\Rs$ of about~$360$ but deteriorates rapidly beyond, the moderate width~$w=40$ is the best fixed width from there up to about~$1140$, and the coarse width~$w=100$ only in the widest bin.

Measured in slabs, the three fixed widths point to a common target: over the configurations solved by all three widths, the width that turns out fastest produces a median of about $18$ slabs per scenario. This motivates determining the starting slab width a priori for each instance: with the initial scenario~$s_0$, we set $w=\max\{\delta,\Rs[s_0]/K\}$, so that the first master problem contains exactly~$K$ slabs for~$s_0$ (the lower bound~$\delta$ is never active in our experiments). 
We ran this automatic width with $K=18$ and $\sigma=0.95$ on all $360$ configurations; it is included in Figure~\ref{fig:sel-extent} as ``auto''. Indeed, it performs best: it has the smallest slowdown in every bin, and it is the fastest setting on $51\%$ of them. On the configurations solved by both ``auto'' and~$w=40$, its geometric mean running time is $41\%$ faster than that of $w=40$ and less than $4\%$ above that obtained by selecting the best fixed width in hindsight for each configuration,
and it hits the time limit as rarely as $w=40$ ($17$ times, against $59$ for $w=10$ and $22$ for $w=100$).

Therefore, our results suggest choosing the starting slab width so that the first master problem contains roughly a constant number of slabs ($18$ in our experiments), and this choice can indeed be made a priori.

\subsubsection{Horizontal versus vertical slabs}\label{sec:comp-orientation}

As noted in Section~\ref{sec:model}, the slab decomposition can equally be built from horizontal slabs, which span the objective space from~$0$ to~$\refpt_1$ in the first coordinate and are refined in the second. While vertical slabs are the default throughout this paper, horizontal slabs can be equivalently used. The best choice largely depends on the chosen starting slab width and the number of resulting slabs, and is therefore linked to the discussion in Section~\ref{sec:comp-slabwidth}. Based on those insights, the intuitive effect is that if we ``stretch'' the objective function in the decomposed direction (e.g., stretch objective~$1$ if vertical slabs are used), this orientation will become preferable until a certain number of slabs is reached, at which point the effect will reverse and the other orientation will become better. 

We now investigate this effect quantitatively by introducing scaling factors~$\hat f_1$ and~$\hat f_2$, which are applied to all first- and second-stage weights of objectives~$1$ and~$2$, respectively.
This stretches the nondominated sets, the reference point, and hence~$\Rs$ (or its counterpart $\refpt_2-\underline f_2^s$ in the second objective) in one coordinate---and thereby changes the number of slabs that results from a fixed~$w$, while leaving the feasible set, and hence the combinatorial difficulty of the instance, untouched.
For these tests, we solved $30$ instances ($n\in\{20,40\}$ with the three values of~$p$ from Table~\ref{tab:sel-main}, $|S|=16$, and five scenario seeds each; $\TG=5\%$, $w=40$, $\sigma=0.95$) once with each orientation, for $\hat f_1\in\{1.4,2,4,10\}$ at $\hat f_2=1$, for $\hat f_2\in\{1.4,2,4,10\}$ at $\hat f_1=1$, and in the unscaled base case. Of the resulting $270$ pairs, $248$ are complete; the remaining ones are excluded from the computation because one of the two orientations hit the time limit, which happens only at the two most extreme weight factors.

\begin{table}
\centering\footnotesize
\begin{tabular}{l|cccc|c|cccc}
\hline
 & \multicolumn{4}{c|}{objective 1 stretched} & base & \multicolumn{4}{c}{objective 2 stretched} \\
$\hat f_1$ & 10 & 4 & 2 & 1.4 & 1 & 1 & 1 & 1 & 1 \\
$\hat f_2$ & 1 & 1 & 1 & 1 & 1 & 1.4 & 2 & 4 & 10 \\
\hline
Time factor & 3.06 & 1.16 & 1.12 & 0.95 & 0.94 & 1.13 & 1.24 & 1.09 & 0.77 \\
Iteration factor & 0.37 & 0.48 & 0.76 & 0.84 & 1.00 & 1.08 & 1.43 & 2.19 & 2.93 \\
\hline
\end{tabular}
\caption{Effect of asymmetric objectives on the two slab orientations (selection problem, $\TG=5\%$, $|S|=16$, starting slab width $w=40$, $\sigma=0.95$). Every weight of objective~$1$ (resp.~$2$) is multiplied by the factor $\hat f_1$ (resp.~$\hat f_2$), which stretches the fronts in that coordinate without changing the feasible set. Each column is one such scaling over $n\in\{20,40\}$ with the three values of $p$ from Table~\ref{tab:sel-main} and five scenario seeds, solved once with vertical and once with horizontal slabs. Reported is the geometric mean over the pairs that both orientations solved (24, 26, 30, 30, 30, 30, 30, 27, 21 of $30$ per column, left to right) of the ratio vertical\,/\,horizontal, for the algorithm time and for the number of iterations; a value below~$1$ means that vertical slabs are the cheaper orientation.}
\label{tab:sel-orientation-scaling}
\end{table}

The results, shown in Table~\ref{tab:sel-orientation-scaling}, confirm the two sides of the trade-off that we expected. 
Starting from a balanced time factor of $1$ (this does not lie exactly at the unscaled base case because, due to some random bias in the instances, the values of~$\Rs$ are slightly larger than their counterparts in objective~$2$), going left in the direction of larger objective-$1$ weights, the vertical orientation is preferable at first (finer decomposition) but then becomes slower (larger master problems) with larger factors~$\hat f_1$. The same effect can be seen in the other direction, with the horizontal orientation being better first, but then overtaken by the vertical one in the long run.
In terms of iterations, however, aligning the slabs with the stretched objective always helps: in the vertical case, for example, stretching objective~$1$ results in a finer decomposition from the first iteration on, and fewer slab-refinement steps have to be taken in total.

The practical conclusion is that neither slab orientation is uniformly preferable; the better choice depends largely on the resulting number of slabs.

\FloatBarrier

\subsubsection{Choice of the reference point}\label{sec:comp-refpoint}

\begin{table}
\centering\footnotesize\setlength{\tabcolsep}{4pt}
\begin{tabular}{rr cc|rr cc}
\hline
$n$ & $p$ & $\refpt_{\mathrm{ND}}$ & $\refpt_{\mathrm{dys}}$ & $n$ & $p$ & $\refpt_{\mathrm{ND}}$ & $\refpt_{\mathrm{dys}}$ \\
\hline
10 & 2 & $(203,246)$ & $(234,266)$ & 20 & 2 & $(216,246)$ & $(252,269)$ \\
 & 4 & $(370,441)$ & $(424,476)$ &  & 6 & $(542,620)$ & $(686,704)$ \\
 & 5 & $(442,520)$ & $(512,551)$ &  & 10 & $(894,963)$ & $(1069,1080)$ \\
\hline
40 & 2 & $(224,246)$ & $(277,285)$ & 60 & 2 & $(224,246)$ & $(277,285)$ \\
 & 11 & $(1002,1058)$ & $(1233,1247)$ &  & 16 & $(1385,1458)$ & $(1740,1818)$ \\
 & 20 & $(1680,1742)$ & $(1966,2038)$ &  & 30 & $(2383,2585)$ & $(2837,3039)$ \\
\hline
\end{tabular}
\caption{Reference points for each instance size $(n, p)$ of the selection problem: the smallest ND-including point~$\refpt_{\mathrm{ND}}$ (used by default) and the dystopian point~$\refpt_{\mathrm{dys}}$. We show seed~$1$ only (the values are similar across seeds) and $|S|=128$, which yields the componentwise largest reference points across the different sizes of $|S|$ used. Table~\ref{tab:sel-refpoints-full} in the appendix lists all instances.}
\label{tab:sel-refpoints}
\end{table}

Appendix~\ref{sec-app:refpoint} discusses two admissible choices of the reference point: the \emph{dystopian} point~$\refpt_{\mathrm{dys}}$, obtained by maximizing each objective separately for every scenario with no fixed variables (in particular, no fixed first-stage decision), and the smallest ND-including point~$\refpt_{\mathrm{ND}}$ of Appendix~\ref{sec-app:refpoint-selection}, the componentwise smallest point that is weakly dominated by every nondominated point for every first-stage decision and scenario. All experiments so far use~$\refpt_{\mathrm{ND}}$. Table~\ref{tab:sel-refpoints} lists both reference points for one exemplary seed and $|S|=128$. The reference point is a modeling input, and the two choices define different hypervolumes. Hence, they could lead to different optimal first-stage decisions. To analyze this effect, we repeated a paired experiment on the $45$ instances given by $n\in\{20,40,60\}$ with the three values of~$p$ from Table~\ref{tab:sel-main}, $|S|=16$, and five scenario seeds each ($\TG=5\%$, $w=40$, $\sigma=0.95$), once with each of the two corresponding reference points. We deliberately use a fixed starting slab width here instead of the automatic width of Section~\ref{sec:comp-slabwidth}: since $\Rs[s_0]=\refpt_1-\underline f_1^{s_0}$ depends on the reference point, the automatic rule would assign the dystopian point a correspondingly coarser starting width and thereby mix the effect of the reference point with that of the starting width.

Of the $45$ pairs of runs, all but one pair terminated within the time limit, leaving $44$ complete pairs. On these instances, the dystopian point lies noticeably further out: the area it spans exceeds that of the corresponding smallest ND-including point by $50\%$ on average, and by at least $5\%$ on every instance (see Table~\ref{tab:sel-refpoints-full} in Appendix~\ref{sec-app:refpoint-selection}).

To compare the returned first-stage decisions, we re-evaluate both of them against one common reference point, recomputing the nondominated sets~$\PND^s(x^*)$ of all scenarios for both of the two fixed first-stage decisions~$x^*$, and comparing $\min_{s\in S} \HV^s(x^*)$. This is done twice, once with~$\refpt_{\mathrm{dys}}$ and once with~$\refpt_{\mathrm{ND}}$ as the common yardstick used in the calculation of the hypervolumes.
Table~\ref{tab:sel-refpoint-compare} reports the outcome together with the time, the certified gap, the iteration count, and the cost of computing the reference point itself.

\begin{table}
\centering\footnotesize\setlength{\tabcolsep}{3pt}
\begin{tabular}{ll|rrr|rrr|rrr}
\hline
 & & \multicolumn{3}{c|}{$n=20$} & \multicolumn{3}{c|}{$n=40$} & \multicolumn{3}{c}{$n=60$} \\
$p$ &  & 2 & 6 & 10 & 2 & 11 & 20 & 2 & 16 & 30 \\
\hline
Time [s] & $\refpt_{\mathrm{ND}}$ & 2.3 & 22.1 & 116 & 2.9 & 268 & 578 & 3.5 & 241 & 1051 \\
 & $\refpt_{\mathrm{dys}}$ & 1.9 & 13.2 & 65.1 & 2.2 & 299 & 589 & 2.8 & 86.6 & 727 \\
Gap [\%] & $\refpt_{\mathrm{ND}}$ & 4.29 & 4.20 & 4.27 & 4.27 & 4.73 & 4.86 & 4.10 & 4.37 & 4.24 \\
 & $\refpt_{\mathrm{dys}}$ & 4.34 & 4.22 & 4.44 & 3.64 & 4.22 & 4.36 & 3.87 & 3.82 & 3.13 \\
Iterations & $\refpt_{\mathrm{ND}}$ & 5.6 & 6.2 & 6.2 & 5.4 & 5.2 & 4.8 & 5.6 & 3.6 & 3.2 \\
 & $\refpt_{\mathrm{dys}}$ & 6.4 & 5.2 & 6.0 & 5.0 & 4.2 & 3.8 & 4.4 & 2.2 & 2.8 \\
$\refpt$ time [ms] & $\refpt_{\mathrm{ND}}$ & 0.39 & 0.51 & 0.60 & 0.69 & 0.92 & 1.14 & 1.04 & 1.33 & 1.50 \\
 & $\refpt_{\mathrm{dys}}$ & 37.22 & 39.88 & 39.10 & 57.10 & 49.70 & 58.13 & 76.82 & 74.26 & 84.57 \\
\hline
$\HV^\ast$ vs.\ $\refpt_{\mathrm{dys}}$ & $\refpt_{\mathrm{ND}}$ & 34.6 & 165.4 & 256.3 & 48.1 & 627.2 & 1040.9 & 58.3 & 1403.4 & 2363.6 \\
 & $\refpt_{\mathrm{dys}}$ & 34.6 & 165.3 & 257.3 & 48.2 & 628.0 & 1040.7 & 58.3 & 1398.9 & 2365.6 \\
$\HV^\ast$ vs.\ $\refpt_{\mathrm{ND}}$ & $\refpt_{\mathrm{ND}}$ & 20.6 & 88.5 & 139.8 & 24.8 & 298.8 & 536.5 & 29.9 & 648.8 & 1157.2 \\
 & $\refpt_{\mathrm{dys}}$ & 20.6 & 88.2 & 139.6 & 24.9 & 295.6 & 534.7 & 29.9 & 641.8 & 1151.6 \\
\hline
$x^\ast$ identical &  & 5/5 & 3/5 & 2/5 & 4/5 & 1/5 & 0/5 & 5/5 & 0/5 & 0/4 \\
\hline
\end{tabular}
\caption{Smallest ND-including reference point $\refpt_{\mathrm{ND}}$ (Appendix~\ref{sec-app:refpoint-selection}) versus the dystopian reference point $\refpt_{\mathrm{dys}}$ on the selection problem.
Every column is one instance size $(n,p)$, each entry averages five scenario seeds (four at the one $(60,30)$ configuration whose ND-including run timed out), and the two rows of every block are paired, i.e.\ they solve exactly the same instances, once with each reference point. ``Time'' is the total algorithm time, ``Gap'' the optimality gap certified on termination (target $5\%$), ``Iterations'' the number of master solves, and ``$\refpt$ time'' the time spent computing the reference point itself. The two $\HV^\ast:=\min_{s\in S}\HV^s(x^\ast)$ blocks give the true worst-case hypervolume of the returned first-stage decision $x^\ast$ (in units of $10^3$), each measured against \emph{one common} reference point so that the two runs are comparable: against $\refpt_{\mathrm{dys}}$ in the first block and against $\refpt_{\mathrm{ND}}$ in the second. The last row counts, per column, the seeds (out of five) on which both runs return the same first-stage decision.}
\label{tab:sel-refpoint-compare}
\end{table}

For the combinatorially easier instances (especially the ones with~$p=2$), the choice of the reference point has almost no effect on the returned solution, while for the harder instances, the returned solution changes. In total, in $20$ of the $44$ pairs the two runs return exactly the same first-stage decision. Where they differ, in most cases each reference point yields the better decision (the one with the larger value for~$\min_{s\in S} \HV^s(x^*)$) under its own yardstick, as one would expect, but the differences are small in both directions (up to a relative difference of $2.4\%$). The certified gaps are also close for both reference points ($4.37\%$ against $4.02\%$ on average). Regarding the algorithm time, we notice that on average it is smaller when the dystopian point is used, even though this reference point itself is more expensive to compute. 
The mechanism here is that the target gap is a relative measure: a reference point further out inflates every hypervolume, so the same absolute approximation error translates into a smaller relative gap and the algorithm certifies~$\TG$ earlier. This is also visible in the iteration counts ($5.14$ against $4.48$ on average).

The takeaway is that the choice of reference point can substantially affect the running time, while the returned first-stage decisions achieve similar worst-case hypervolumes when evaluated against a common reference point. At a fixed starting slab width, the more distant dystopian point yielded shorter running times on average, since the larger hypervolumes made the relative target gap easier to certify. Nevertheless, the smallest ND-including reference point may still be preferable from a modeling perspective, as it avoids unnecessarily expanding the region used to evaluate the nondominated sets.

\subsubsection{Accuracy of the refinement factor~$\alpha$ and the safety factor~$\sigma$}\label{sec:comp-alpha}

Recall from Section~\ref{sec:gapeval} that the refinement factor~$\alpha$ is chosen such that the gap would exactly attain the target gap~$\TG$ \emph{if} $\tilde z$ and $\min_{s\in S}\HV_{\mathrm{LB}}^s(\tilde x)$ remained unchanged by the refinement. Since this is not necessarily the case, a single refinement does not always suffice, and the algorithm may have to refine several times in a row. We call a maximal sequence of consecutive refinement iterations (without a scenario addition in between) a \emph{refinement chain}. Every refinement chain of length greater than one means that the predicted~$\alpha$ was too optimistic and additional master solves were required. 
The safety factor~$\sigma<1$, which we introduced in Section~\ref{sec:gapeval}, causes the refinement to be slightly more aggressive, and we sweep over $\sigma\in\{0.5,0.8,0.95,0.99\}$ to measure how effective this heuristic remedy is.

\begin{figure}[!htbp]
	\centering
	\includegraphics[width=\textwidth]{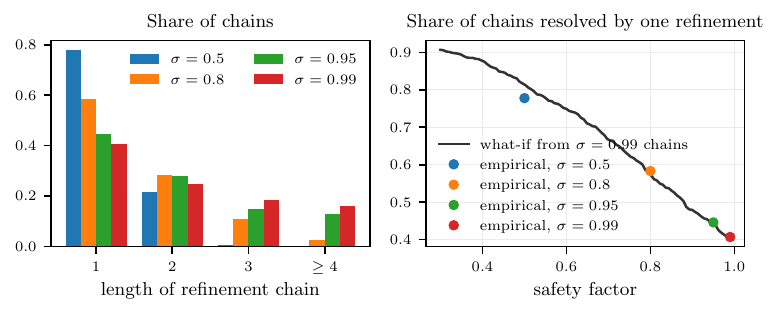}
	\caption{Accuracy of the refinement factor~$\alpha$. Left: distribution of the lengths of refinement chains per safety factor~$\sigma$. Right: retrospective estimate of the share of chains for which a single refinement with safety factor \(c\) would have reached a slab width no larger than the final width of the recorded chain (\(\sigma=0.99\)). The dots show the empirically observed shares of length-one chains for the four tested values of \(\sigma\).}
	\label{fig:sel-alpha}
\end{figure}

As can be seen in Figure~\ref{fig:sel-alpha}, at the highest value of~$\sigma$ considered ($\sigma=0.99$), $59\%$ of the chains have length greater than one, which means that the predicted refinement factor is too optimistic in more than half of the cases. Decreasing the safety factor reduces the share of chains longer than one monotonically. 
The right panel of Figure~\ref{fig:sel-alpha} provides a retrospective, width-based estimate of this trade-off. For a refinement chain recorded with safety factor~$\sigma$ and refinement factors $\alpha_1,\dots,\alpha_k$, a hypothetical smaller safety
factor~$c\le\sigma$ would have yielded the first refinement factor
\[
    \alpha_1(c)=
    \max\left\{\frac{c}{\sigma}\cdot\alpha_1,\alpha_{\min}\right\}.
\]
We count a chain as potentially resolvable by a single refinement if $\alpha_1(c)\le\prod_{i=1}^k\alpha_i$,
i.e., if the first refinement would have reduced the slab width to at most the final width reached by the recorded chain. This criterion provides only a heuristic estimate, since refining directly to a different slab width may change the master solution and the subsequent decisions of the algorithm. We derive the curve from the chains recorded with the highest safety factor considered,~$\sigma=0.99$, because the retrospective calculation is only valid for values of~$c$ that do not exceed the safety factor~$\sigma$ used to generate the chains. The empirically observed shares of length-one chains are broadly consistent with this retrospective estimate.

This behavior is reflected directly in the number of iterations, as seen in Table~\ref{tab:sel-main}. However, the algorithm time does not follow the same monotone trend, because a smaller~$\sigma$ also overshoots the required slab width in many cases and therefore produces unnecessarily large master problems. 

Overall, at $\sigma=0.99$, more than half of the refinement chains require multiple refinements. A smaller safety factor reduces the frequency of such repeated refinements but may increase running time by producing unnecessarily large master problems. Thus, a moderate value such as our default~$\sigma=0.95$ provides a good compromise in our experiments.

\FloatBarrier
\subsection{Results for the assignment problem}\label{sec:comp-ass-half}

We now turn to the two-stage robust assignment problem of Section~\ref{sec:ass-problem}. The largest values of~$\Rs$ are smaller than for the selection problem (instance means of up to $984$, against $1505$ in Section~\ref{sec:comp-sel}), so we use starting slab widths that are half as large as there, namely $w\in\{5,20,50\}$.  
The runs otherwise follow the structure of Section~\ref{sec:comp-design}, with the sizes $n\in\{6,8,10\}$ ($p\in\{0,2,3\}$ for $n=6$, $p\in\{0,2,4\}$ for $n=8$, and $p\in\{0,2,5\}$ for $n=10$), scenario-set sizes $|S|\in\{4,16,128\}$, and safety factors $\sigma\in\{0.8,0.95,0.99\}$; with five scenario seeds each, this gives $2430$ runs. Here~$p$ is the \emph{minimum} number of first-stage edges, so $p=0$ is the extreme case in which the entire matching may be deferred to the second stage. The instances are much smaller in~$n$ than for the selection problem, but the $n^2$ second-stage variables per scenario and the matching structure make the resulting master problems comparably hard or harder.

At $n=10$, the computation of the smallest ND-including reference point becomes a significant computational cost: the combinatorial routine of Appendix~\ref{sec-app:refpoint-assignment} needs $140$\,s on average (at most $326$\,s) for $|S|\le16$, and for all $15$ instances with $|S|=128$ it exceeded its time limit of $600$\,s. The runs on these instances therefore fell back to the dystopian reference point. They are marked in Tables~\ref{tab:ass-main-halfsw} and~\ref{tab:ass-refpoints-halfsw} and should be read with Section~\ref{sec:comp-refpoint} in mind: a reference point further out can make the relative target gap easier to certify, so their running times may be optimistic compared with the other rows, particularly when assessing
the scaling with~$|S|$. The reference points for the respectively smallest values of~$p$ are listed in Table~\ref{tab:ass-refpoints-halfsw}.

\begin{table}
\centering\footnotesize\setlength{\tabcolsep}{6pt}
\begin{tabular}{r|ccc}
\hline
$n$ & $|S|=4$ & $|S|=16$ & $|S|=128$ \\
\hline
$6$ & $(545, 561)$ & $(565, 563)$ & $(618, 640)$ \\
$8$ & $(705, 771)$ & $(738, 793)$ & $(828, 831)$ \\
$10$ & $(882, 933)$ & $(947, 1001)$ & $(1109, 1169)^\dagger$ \\
\hline
\end{tabular}
\caption{Reference point~$\refpt$ used by the algorithm for the assignment problem---the smallest ND-including point of Appendix~\ref{sec-app:refpoint-assignment}---for each instance size~$n$ and scenario-set size~$|S|$. The point depends only mildly on~$p$ (shown here for seed~$1$ and the smallest~$p$; across the values of~$p$ it varies by at most $2.6\%$ there), and, apart from the entry marked~$^\dagger$, none of its components grows by more than $17\%$ with~$|S|$. Entries marked~$^\dagger$ are the dystopian point instead: the combinatorial routine exceeded its time limit of $600$\,s there.}
\label{tab:ass-refpoints-halfsw}
\end{table}

\begin{sidewaystable}
\centering\footnotesize
\begin{tabular}{lll|r|rrr|rrr}
\hline
 & & & & \multicolumn{3}{c|}{$\TG=5\%$} & \multicolumn{3}{c}{$\TG=10\%$} \\
$n$ & $p$ & $|S|$ & RP [s] & Solved & Algorithm [s] & Iterations & Solved & Algorithm [s] & Iterations \\
\hline
6 & 0 & 4 & 0.29 & 15/15 & 4.8\,/\,4.8\,/\,5.5 & 2.6\,/\,2.8\,/\,3.2 & 15/15 & 2.6\,/\,2.4\,/\,2.8 & 2.0\,/\,2.0\,/\,2.0 \\
 &  & 16 & 1.49 & 15/15 & 20.4\,/\,17.0\,/\,18.7 & 4.0\,/\,4.2\,/\,4.8 & 15/15 & 4.8\,/\,5.2\,/\,5.7 & 2.2\,/\,2.2\,/\,2.2 \\
 &  & 128 & 9.41 & 15/15 & 79.6\,/\,295\,/\,278 & 4.0\,/\,4.2\,/\,5.0 & 15/15 & 16.4\,/\,13.1\,/\,15.7 & 2.6\,/\,2.6\,/\,2.6 \\
 & 2 & 4 & 0.39 & 15/15 & 18.5\,/\,15.6\,/\,18.5 & 3.6\,/\,3.6\,/\,3.6 & 15/15 & 11.0\,/\,10.6\,/\,10.6 & 2.6\,/\,2.6\,/\,2.6 \\
 &  & 16 & 1.24 & 15/15 & 205\,/\,192\,/\,125 & 4.4\,/\,4.8\,/\,5.0 & 15/15 & 41.7\,/\,38.1\,/\,38.7 & 3.8\,/\,3.8\,/\,3.8 \\
 &  & 128 & 11.04 & 15/15 & 107\,/\,166\,/\,145 & 4.2\,/\,5.0\,/\,5.0 & 15/15 & 18.1\,/\,17.7\,/\,15.0 & 2.4\,/\,2.4\,/\,2.4 \\
 & 3 & 4 & 0.22 & 15/15 & 13.1\,/\,16.4\,/\,18.6 & 3.4\,/\,3.8\,/\,4.2 & 15/15 & 8.5\,/\,7.8\,/\,7.8 & 3.0\,/\,3.0\,/\,3.0 \\
 &  & 16 & 1.01 & 15/15 & 131\,/\,86.2\,/\,158 & 5.4\,/\,5.8\,/\,6.4 & 15/15 & 27.0\,/\,22.9\,/\,34.9 & 3.8\,/\,3.6\,/\,4.2 \\
 &  & 128 & 8.29 & 15/15 & 268\,/\,213\,/\,150 & 5.0\,/\,5.4\,/\,5.2 & 15/15 & 20.5\,/\,21.3\,/\,21.4 & 3.2\,/\,3.2\,/\,3.2 \\
\hline
8 & 0 & 4 & 5.09 & 15/15 & 20.7\,/\,20.5\,/\,19.5 & 1.6\,/\,1.6\,/\,1.6 & 15/15 & 11.9\,/\,12.4\,/\,10.8 & 1.4\,/\,1.4\,/\,1.4 \\
 &  & 16 & 18.31 & 15/15 & 363\,/\,378\,/\,403 & 2.6\,/\,2.6\,/\,2.6 & 15/15 & 34.3\,/\,39.4\,/\,34.1 & 2.2\,/\,2.2\,/\,2.2 \\
 &  & 128 & 141 & 15/15 & 417\,/\,377\,/\,373 & 4.0\,/\,4.2\,/\,4.0 & 15/15 & 207\,/\,190\,/\,186 & 3.0\,/\,3.0\,/\,3.0 \\
 & 2 & 4 & 4.56 & 15/15 & 48.0\,/\,55.4\,/\,49.6 & 2.6\,/\,2.6\,/\,2.6 & 15/15 & 36.1\,/\,43.0\,/\,45.1 & 2.2\,/\,2.2\,/\,2.2 \\
 &  & 16 & 18.79 & 15/15 & 488\,/\,485\,/\,488 & 4.0\,/\,4.0\,/\,4.0 & 15/15 & 151\,/\,146\,/\,145 & 3.0\,/\,3.0\,/\,3.0 \\
 &  & 128 & 137 & 15/15 & 2446\,/\,1952\,/\,1738 & 6.2\,/\,6.0\,/\,5.8 & 15/15 & 342\,/\,315\,/\,335 & 3.4\,/\,3.4\,/\,3.4 \\
 & 4 & 4 & 3.32 & 15/15 & 298\,/\,259\,/\,276 & 2.6\,/\,2.6\,/\,2.6 & 15/15 & 257\,/\,255\,/\,241 & 2.6\,/\,2.6\,/\,2.6 \\
 &  & 16 & 12.60 & 15/15 & 1165\,/\,1028\,/\,1263 & 5.2\,/\,5.2\,/\,5.4 & 15/15 & 864\,/\,854\,/\,896 & 4.8\,/\,4.8\,/\,4.8 \\
 &  & 128 & 77.90 & 12/15 & 2472\,/\,2418\,/\,2556 & 5.2\,/\,5.5\,/\,5.8 & 12/15 & 2340\,/\,2345\,/\,2321 & 4.2\,/\,4.2\,/\,4.2 \\
\hline
10 & 0 & 4 & 54.66 & 15/15 & 99.7\,/\,92.2\,/\,90.8 & 1.8\,/\,1.8\,/\,1.8 & 15/15 & 48.3\,/\,48.0\,/\,50.5 & 1.4\,/\,1.4\,/\,1.4 \\
 &  & 16 & 237 & 15/15 & 2353\,/\,2278\,/\,2371 & 2.0\,/\,2.0\,/\,2.0 & 15/15 & 445\,/\,423\,/\,440 & 1.6\,/\,1.6\,/\,1.6 \\
 &  & 128 & $>600^\dagger$ & 12/15 & 4549\,/\,4531\,/\,4573 & 3.5\,/\,3.5\,/\,3.5 & 15/15 & 1889\,/\,1811\,/\,1814 & 2.6\,/\,2.6\,/\,2.6 \\
 & 2 & 4 & 70.05 & 15/15 & 252\,/\,258\,/\,251 & 2.0\,/\,2.0\,/\,2.0 & 15/15 & 50.9\,/\,45.0\,/\,46.8 & 1.4\,/\,1.4\,/\,1.4 \\
 &  & 16 & 262 & 12/15 & 2664\,/\,2639\,/\,2682 & 3.0\,/\,3.0\,/\,3.0 & 12/15 & 2268\,/\,2258\,/\,2263 & 1.8\,/\,1.8\,/\,1.8 \\
 &  & 128 & $>600^\dagger$ & 6/15 & 7604\,/\,7544\,/\,7523 & 3.0\,/\,3.0\,/\,3.0 & 15/15 & 1340\,/\,1335\,/\,1309 & 3.0\,/\,3.0\,/\,3.0 \\
 & 5 & 4 & 41.85 & 15/15 & 243\,/\,240\,/\,235 & 2.2\,/\,2.2\,/\,2.2 & 15/15 & 189\,/\,187\,/\,171 & 2.0\,/\,2.0\,/\,2.0 \\
 &  & 16 & 175 & 15/15 & 1971\,/\,1982\,/\,1941 & 3.8\,/\,3.8\,/\,3.8 & 15/15 & 825\,/\,788\,/\,850 & 3.0\,/\,3.0\,/\,3.0 \\
 &  & 128 & $>600^\dagger$ & 3/15 & 10227\,/\,10272\,/\,10292 & 6.0\,/\,6.0\,/\,6.0 & 3/15 & 9228\,/\,9185\,/\,9196 & 4.0\,/\,4.0\,/\,4.0 \\
\hline
\multicolumn{4}{l|}{Average $n=6$} & 100\% & 94.2\,/\,112\,/\,102 & 4.1\,/\,4.4\,/\,4.7 & 100\% & 16.7\,/\,15.5\,/\,17.0 & 2.8\,/\,2.8\,/\,2.9 \\
\multicolumn{4}{l|}{Average $n=8$} & 98\% & 857\,/\,775\,/\,796 & 3.8\,/\,3.8\,/\,3.8 & 98\% & 472\,/\,467\,/\,468 & 3.0\,/\,3.0\,/\,3.0 \\
\multicolumn{4}{l|}{Average $n=10$} & 80\% & 3329\,/\,3315\,/\,3329 & 2.7\,/\,2.7\,/\,2.7 & 89\% & 1809\,/\,1787\,/\,1793 & 2.1\,/\,2.1\,/\,2.1 \\
\hline
\multicolumn{4}{l|}{Overall average} & 93\% & 1427\,/\,1401\,/\,1409 & 3.6\,/\,3.7\,/\,3.8 & 96\% & 766\,/\,756\,/\,760 & 2.7\,/\,2.7\,/\,2.7 \\
\hline
\end{tabular}
\caption{Results for the assignment problem with starting slab width $w=20$. Each cell of the ``Algorithm'' and ``Iterations'' columns contains three values, one per safety factor $\sigma=0.8\,/\,0.95\,/\,0.99$; each value averages the five scenario seeds. ``RP'' is the time of the combinatorial reference-point computation (Appendix~\ref{sec-app:refpoint-assignment}; independent of the target gap, the starting slab width and~$\sigma$, and shown also where no run finished, as it is computed before the algorithm starts), ``Solved'' the number of runs (out of 5 seeds $\times$ $3$ values of~$\sigma$) that finished within the time limit of three hours, ``Algorithm'' the wall-clock time of the decomposition loop together with the construction of the initial master problem, i.e.\ everything except the reference point; timed-out runs enter it with the limit of $10\,800$\,s. ``Iterations'' is the number of master solves, averaged over the solved runs only. Here $p$ is the minimum number of first-stage edges, so $p=0$ allows the entire matching to be deferred to the second stage. In the rows marked~$^\dagger$ the combinatorial routine exceeded its time limit of $600$\,s, and the runs used the dystopian reference point instead (Section~\ref{sec:comp-refpoint}).}
\label{tab:ass-main-halfsw}
\end{sidewaystable}

\begin{figure}[!htbp]
	\centering
	\includegraphics[width=\textwidth]{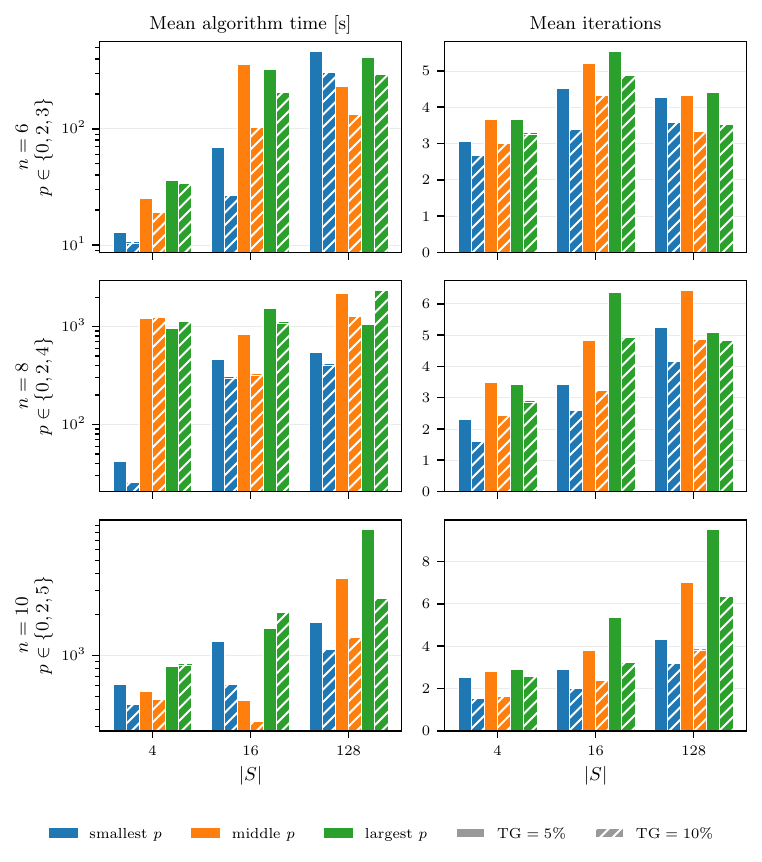}
	\caption{Assignment problem: mean algorithm time (left column, logarithmic scale) and mean number of iterations (right column) at the standard safety factor~$\sigma=0.95$, with one row each for $n=6$, $8$, and~$10$. Within each panel the bars are grouped by the number~$|S|$ of scenarios; within each group, the bar color distinguishes the smallest, middle, and largest value of~$p$ given in the row label, and the hatching the target gap (solid: $5\%$, hatched: $10\%$). Each bar averages the solved runs over the five scenario seeds and the three starting slab widths $w\in\{5,20,50\}$. At $n=8$ and especially at $n=10$, where $15\%$ and $32\%$ of the runs timed out, the averages are subject to survivor bias; at $n=10$ with $|S|=128$, the runs use the dystopian reference point.}
	\label{fig:ass-runtime-halfsw}
\end{figure}

\begin{figure}[!htbp]
	\centering
	\includegraphics[width=\textwidth]{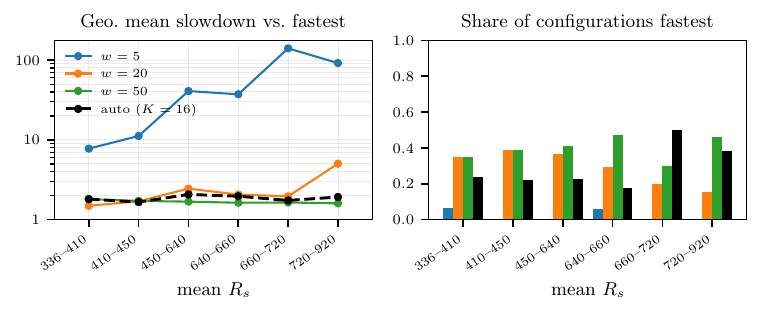}
	\caption{Starting slab width compared to~$\Rs$ for the assignment problem ($\sigma=0.95$). Configurations are binned by the mean of~$\Rs$ over the scenarios~$s\in S$ of the instance. Besides the three fixed starting slab widths, the automatic width $w=\max\{\delta,\Rs[s_0]/16\}$ (``auto'', black) is shown. Left: geometric mean slowdown relative to the fastest of the four settings for the same configuration (logarithmic axis). Right: share of configurations for which each setting is the fastest choice. Both panels use only the configurations solved by all four settings, so they understate the disadvantage of $w=5$, which times out on $104$ of the $270$ configurations.}
	\label{fig:ass-slabwidth-halfsw}
\end{figure}

Table~\ref{tab:ass-main-halfsw} reports the results in the same format as Table~\ref{tab:sel-main}, for the middle starting width $w=20$, the counterpart of the default $w=40$ of the selection problem, and Figure~\ref{fig:ass-runtime-halfsw} shows the scaling of the algorithm time and the iteration count with~$|S|$ at the standard safety factor~$\sigma=0.95$; the behavior is similar to the one observed for the selection problem (Figure~\ref{fig:sel-runtime}). The scenario-addition scheme remains as effective as for the selection problem: even for $|S|=128$, on average only $3.8$ scenarios (at most $12$) are ever added to~$D$ (including the initial one). 

We also repeated the analysis of the starting slab width effects, shown in Figure~\ref{fig:ass-slabwidth-halfsw}. The behavior is similar to that for the selection problem, with a slight tilt toward the coarser choice of starting slab width. The geometric mean slowdowns relative to the fastest choice are $18.7$, $1.73$, and $1.45$ for $w=5$, $20$, and $50$, respectively (over the $496$ configurations solved by all three widths).
The coarsest width $w=50$ is thus the best choice overall. At the tighter target gap of~$5\%$ it is almost tied with the middle width~$w=20$ ($1.50$ against $1.49$). As before, a fine start pays for its refinement-free runs ($w=5$ never refines at all, against on average $0.51$ and $1.75$ refinements for $w=20$ and $w=50$ at $\TG=5\%$) with larger master problems, and in Figure~\ref{fig:ass-slabwidth-halfsw} ($\sigma=0.95$) the slowdown of $w=5$ climbs steeply with~$\Rs$, from $7.8$ in the bin with the smallest values to between $37$ and $141$ in the three widest bins. Among the fixed widths, $w=20$ and $w=50$ are equally often fastest for small values of~$\Rs$, and $w=50$ for all larger ones. Measured in slabs, the fastest fixed width produces a median of about $16$ slabs per scenario, close to the $18$ of the selection problem. Accordingly, we also ran the automatic width of Section~\ref{sec:comp-slabwidth} with $K=16$ and $\sigma=0.95$, included in Figure~\ref{fig:ass-slabwidth-halfsw}. It is as good as the best fixed width $w=50$ overall, with the same number of timeouts ($7$, against $16$ for $w=20$ and $104$ for $w=5$) and almost the same mean algorithm time ($632$\,s against $635$\,s), and it is $24\%$ faster than $w=20$ on the configurations solved by both. Against $w=50$, it is $7\%$ slower on these configurations, being slower on the easy and faster on the hard ones.

\section{Conclusion}\label{sec:concl}

We have presented a column-and-constraint generation algorithm for bi-objective two-stage robust optimization that uses the hypervolume indicator to evaluate the nondominated sets of each scenario. A slab decomposition of the objective space turns the master problem into a linear MIP with a bounded discretization error. The algorithm adaptively decides whether to add a scenario or to refine the decomposition and ultimately certifies a user-specified optimality gap.

We evaluated the algorithm on a two-stage robust selection and a two-stage robust assignment problem. The dominant cost is the master problem, and its size is governed by the number of slabs it contains. 
In our experiments, the fastest configurations were those in which the range of the first objective was covered by about $15$ to $20$ slabs. The starting slab width should therefore be chosen accordingly or, equivalently, the desired number of starting slabs can be prescribed directly to the algorithm and the width derived from it.

We further showed in our tests that, even with $128$ scenarios, fewer than four of them were added to the scenario subset on average. This observation shows that the main computational advantage of \CCG{} carries over to our bi-objective version of two-stage robust optimization. Although every scenario now contributes an entire nondominated set instead of a single objective value, the worst case was determined by a few scenarios in all our experiments, just as it typically is for a scalar objective. This makes the first-scalarize-then-robustify paradigm practical with the hypervolume. 

Furthermore, if a nonzero final optimality gap is accepted, an exact model of the quality indicator is not needed. The algorithm balances two sources of inexactness---the restricted set of scenarios and the approximation of the hypervolume---against the target optimality gap, and it only requires a linear lower bound on the indicator whose error can be bounded and reduced by refinement. This decouples the choice of the indicator from the question of whether it admits a tractable exact formulation.

Finally, our experiments on the selection problem show that the choice of reference point can noticeably affect the running time, while the returned first-stage decisions achieve similar worst-case hypervolumes when evaluated against a common reference point (differences of at most~$2.4\%$). A more distant reference point increases the range covered by the slab decomposition and can thus lead to larger master problems for a fixed starting slab width. However, it can also make the relative target gap easier to certify. In our experiments, the more distant dystopian reference point resulted in shorter running times on average, showing that a smaller reference point is not necessarily computationally preferable. Since the reference point also determines the hypervolume objective itself, its choice remains an important modeling decision. We characterized the computation of the smallest ND-including reference point as a three-level problem and developed combinatorial procedures for both applications.

\enlargethispage{\baselineskip}

These insights also indicate where future work is most worthwhile. A tighter master problem formulation, or more spatially adaptive refinement that acts only where the error margin is large, would attack the bottleneck directly.
For problems whose exact reference point is expensive, such as larger assignment instances, a cheaper way of combinatorially computing a smallest ND-including reference point could be explored.

Beyond that, the framework is not necessarily restricted to the hypervolume: it could potentially be extended to other unary quality indicators for which a linear lower bound with a controllable error can be constructed. Combining several indicators could also be investigated in future work.

\backmatter

\section*{Statements and declarations}

\subsection*{Funding}
\noindent
This work was supported by the German Federal Ministry of Research, Technology and Space (BMFTR) [grant number~05M22WTA].

\subsection*{Author contributions}
{\parindent0pt
\textbf{Johannes Kager:} Methodology, Software, Validation, Formal analysis, Investigation, Data Curation, Writing - Original Draft, Writing - Review \& Editing, Visualization.

\textbf{Clemens Thielen:} Conceptualization, Methodology, Validation, Formal analysis, Writing - Review \& Editing, Supervision, Project Administration, Funding acquisition.
}

\subsection*{Declaration of interest}
\noindent
Declarations of interest: none

% \subsection*{Acknowledgments}
% \noindent

\subsection*{Data availability}
\noindent
The implementations of the algorithms and models as well as the instance generator used in the computational experiments are published on~\zenodolink.

\bibliography{literature}
%% if required, the content of the .bbl file can be included here once it is generated
%%\input main.bbl

\newpage

\begin{appendices}

\section{An exact quadratic model of the hypervolume}\label{sec-app:hv-opt}

As stated in Section~\ref{sec:fund-biobj}, the hypervolume for a scenario~$s\in S$ and a first-stage decision~$x\in\X$ can also be modeled exactly when the set~$\PND^s(x)$ of nondominated points is not explicitly available. The corresponding model is the following optimization problem, which is a mixed-integer quadratic program with linear constraints and in which $n\geq 2$ second-stage decisions $y^{s,1},\dots,y^{s,n}$ describe (a subset of) the nondominated set:
\begin{equation}\label{eq:hv-opt}\tag{OPT-HV}
	\begin{aligned}
		\max \quad
		 & \lefteqn{\bigl(\refpt_1-f_1^s(x,y^{s,n})\bigr)\cdot \bigl(\refpt_2-f_2^s(x,y^{s,n})\bigr) \;+} && \\ 
		 & \lefteqn{\sum_{i=1}^{n-1} \bigl(f_1^s(x,y^{s,i+1})-f_1^s(x,y^{s,i})\bigr)\cdot \bigl(\refpt_2-f_2^s(x,y^{s,i})\bigr)} &&                                           \\
		\mathrm{s.t.}\quad
		 & f_1^s(x,y^{s,i}) \le f_1^s(x,y^{s,i+1}) \le \refpt_1                                                   &  & \forall i\in[n-1] \\
		 & f_2^s(x,y^{s,i+1}) \le f_2^s(x,y^{s,i}) \le\refpt_2                                                    &  & \forall i\in[n-1] \\
		 & y^{s,i} \in \Y^s(x)                                                                                    &  & \forall i\in[n].
	\end{aligned}
\end{equation}

The constraints of \eqref{eq:hv-opt} force the points $f^s(x,y^{s,1}),\dots,f^s(x,y^{s,n})$ to weakly dominate~$\refpt$ and to be ordered by non-decreasing first and non-increasing second component, so that the objective function sums the areas of the corresponding rectangles as in Equation~\eqref{eq:hv-calculation}. Consequently, if $\PND^s(x)$ is a finite set, and if~$n$ is large enough for the~$n$ second-stage decisions to describe all nondominated points weakly dominating~$\refpt$, the optimal objective value of \eqref{eq:hv-opt} equals~$\HV^s(x)$. We state Problem~\eqref{eq:hv-opt} only to illustrate that an exact model of the hypervolume is quadratic; it is not used algorithmically anywhere in this paper, which is why we omit a formal proof of its correctness.

\section{Building $\AA$ combinatorially}\label{sec-app:AAcomb}
In the following, let~$s\in S$ be a scenario and~$x\in\X$ be a first-stage decision. Let~$P = \left\{p^1, \dots, p^m\right\} \subseteq \PND^s(x)$ be a subset of the nondominated set for~$s$ and~$x$, and let~$\PSOL = \left\{(x,\tilde y^1),\dots,(x,\tilde y^m)\right\}$ be a set of corresponding efficient solutions---e.g., as returned by the epsilon-constraint method in Section~\ref{sec:fund-biobj}.
Assume the points in~$P$ are ordered such that, for increasing index, the first component increases. Further, let the index set~$G^s$, the vector of slab positions~$(t^i)_{i\in G^s}$, and the width $w$ describe a slab decomposition as in Section~\ref{sec:model}. For each~$i\in G^s$, we define
\begin{align*}
	j(i) \coloneqq
	\begin{cases*}
		-1                                                         & if $\nexists j \in [m]: p_1^j \le t^i \text{ and } p_2^j \le \refpt_2$ \\
		\max \{j\in [m]: p_1^j \le t^i \text{ and } p_2^j \le \refpt_2\} & else.
	\end{cases*}
\end{align*}

We now show that any such subset~$P$ can be used to construct a feasible solution of the master problem~$\MAS_{\{s\}}(x)$, whose objective value is therefore a lower bound on~$\AA^s(x)$---this is the lower bound~$\AA_{\text{LB}}^s(x)$ used in Algorithm~\ref{alg:slabdecomp}, where~$P$ is the (generally incomplete) subset returned by the epsilon-constraint method. If~$P$ contains \emph{all} nondominated points for~$s$ and~$x$ that weakly dominate the reference point, the constructed solution is even optimal, so that~$\AA^s(x)$ can be computed combinatorially in this case.
\begin{lemma}\label{lemma:AA-construct}
	If $G^s$ is empty, then $(x, (), 0, (), ())$ is an (optimal) master solution for~$\MAS_{\{s\}}(x)$ (with $()$ being an empty vector, since no~$y,h$ or $a$-variables exist in this case).
	Otherwise, let for each~$i\in G^s$,
	\begin{align*}
				y^{s,i} &= \begin{cases*}
					\tilde y^{j(i)} & if $j(i) \ne -1$ \\
					\text{any element } y\in\Y^s(x) \text{ with } f^s(x,y)\in [0,\refpt_1]\times[0,\refpt_2] & else
				\end{cases*}  \\
		h^{s,i}  & =
		\begin{cases*}
			\refpt_2-p^{j(i)}_2 & if $j(i) \ne -1$ \\
			0                   & else
		\end{cases*}
		\qquad &                                         \\
		a^{s,i}  & =
		\begin{cases*}
			1 & if $j(i) \ne -1$ \\
			0 & else
		\end{cases*}
		\qquad &                                         
	\end{align*}
	and~$z = \sum_{i \in G^s} w\cdot h^{s,i}$.
    Then the following holds:
	\begin{enumerate}
		\item[(a)] $(x, y^s, z, h^s, a^s)$ is a feasible solution of~$\MAS_{\{s\}}(x)$ with objective value~$z$; in particular, $z \le \AA^s(x)$.
		\item[(b)] If~$P$ contains every nondominated point for~$s$ and~$x$ that weakly dominates the reference point, i.e., if $\{p\in\PND^s(x): p_\ell\le\refpt_\ell \text{ for } \ell=1,2\}\subseteq P$, then $(x, y^s, z, h^s, a^s)$ is a master solution for~$\MAS_{\{s\}}(x)$, i.e., $z = \AA^s(x)$.
	\end{enumerate}
\end{lemma}
\begin{proof}
	The first case with empty $G^s$ is trivial.
	For the second case, we first show feasibility for \(\MAS_{\{s\}}\).

	By construction, \(h^{s,i}\ge 0\) and \(a^{s,i}\in\{0,1\}\) for all \(i\in G^s\). Moreover,
	\(y^{s,i}\in \Y^s( x)\) for all \(i\in G^s\), and \( x\in \X\). It remains to verify the constraints~\eqref{eq:master-3} to~\eqref{eq:master-6}.

	Constraints~\eqref{eq:master-3} and~\eqref{eq:master-4} are satisfied directly by construction.
	For constraints~\eqref{eq:master-5} and~\eqref{eq:master-6}, we consider two cases.
	First, assume~$j(i)=-1$. Then, $a^{s,i}=0$ and~$h^{s,i}=0$. By construction, $y^{s,i}$ satisfies $f_1^s(x,y^{s,i}) \le \refpt_1$ and $f_2^s(x,y^{s,i}) \le \refpt_2$, and thus,
	\begin{equation*}
		f_1^s(x,y^{s,i}) \le \refpt_1 = (1-a^{s,i})\refpt_1+a^{s,i} t^i,
	\end{equation*}
	and
	\begin{equation*}
		f_2^s(x,y^{s,i}) \le \refpt_2 = \refpt_2 - h^{s,i}.
	\end{equation*}

	For the case of~$j(i)\ne-1$, we have $a^{s,i}=1$ and $h^{s,i}=\refpt_2-p_2^{j(i)}$. By definition, we have~$p_1^{j(i)} \le t^i$. Thus,
	\begin{equation*}
		f_1^s(x,y^{s,i}) = f_1^s(x,\tilde y^{j(i)}) = p^{j(i)}_1 \le t^i = (1-a^{s,i})\refpt_1+a^{s,i} t^i,
	\end{equation*}
	and
	\begin{equation*}
		f_2^s(x,y^{s,i}) = f_2^s(x,\tilde y^{j(i)}) = p_2^{j(i)} = \refpt_2 - h^{s,i}.
	\end{equation*}
	So the constructed solution is feasible, and its objective value is~$z$, which proves~(a).

	For~(b), assume that every nondominated point for~$s$ and~$x$ that weakly dominates the reference point is contained in~$P$; it remains to prove optimality. Let \((x, \check y^s,\check z,\check h^s,\check a^s)\) be any feasible solution of \(\MAS_{\{s\}}(x)\) for the fixed first-stage decision \(x\). For each \(i\in G^s\), feasibility implies
	\[
		f_1^s(x,\check y^{s,i})\le (1-\check a^{s,i})\refpt_1+\check a^{s,i} t^i.
	\]
	If \(\check a^{s,i}=0\), then \(\check h^{s,i}\le \refpt_2\check a^{s,i}=0\), hence \(\check h^{s,i}=0\le h^{s,i}\).
    If \(\check a^{s,i}=1\), then
	\[
		f_1^s(x,\check y^{s,i})\le t^i.
	\]
	Let \(q\) be a nondominated point in \(\PND^s(x)\) weakly dominating  \(f^s(x,\check y^{s,i})\) ($q$ exists by Assumption~\ref{ass:external-stability}, since~$\check y^{s,i}\in\Y^s(x)$). Then \(q_1\le t^i \le \refpt_1\) and, using the constraint \(f_2^s(x,\check y^{s,i})\le \refpt_2-\check h^{s,i}\), also \(q_2 \le f_2^s(x,\check y^{s,i}) \le \refpt_2\). Hence, $q$ weakly dominates the reference point and, by the assumption of~(b), $q\in P$, say $q=p^{j_q}$. Since \(p_1^{j_q}\le t^i\) and \(p_2^{j_q}\le\refpt_2\), the definition of \(j(i)\) yields \(j_q\le j(i)\). Since the points in~$P$ are ordered increasingly in the first component (and, thus, decreasingly in the second component), this implies
	\[
		p_2^{j(i)} \le p_2^{j_q} = q_2 \le f_2^s(x,\check y^{s,i}).
	\]
	Using the constraint \(f_2^s(x,\check y^{s,i})\le \refpt_2-\check h^{s,i}\), we obtain
	\[
		\check h^{s,i} \le \refpt_2-f_2^s(x,\check y^{s,i})
		\le \refpt_2-p_2^{j(i)}=h^{s,i}.
	\]
	Thus, \(\check h^{s,i}\le h^{s,i}\) for every \(i\in G^s\). Summing over all slabs yields
	\[
		\check z \le \sum_{i\in G^s} w\,\check h^{s,i}
		\le \sum_{i\in G^s} w\,h^{s,i} = z.
	\]
	Therefore, no feasible solution attains an objective value larger than \(z\), and the constructed solution is optimal for \(\MAS_{\{s\}}(x)\), which proves~(b).
\end{proof}

\section{Reference point}\label{sec-app:refpoint}

Throughout this paper, we evaluate sets of points using the hypervolume indicator with respect to a fixed reference point~$\refpt$. 
Note that neither the definition of the hypervolume nor the algorithms and models require the reference point to be dominated by all points in the nondominated sets (for each scenario in~$S$ and each first-stage decision in~$\X$).
By definition, our algorithm only considers or finds points that (weakly) dominate the reference point. Therefore, the validity of the algorithm is independent of the concrete choice of~$\refpt$.

Regardless of the general validity, it is naturally desirable to use a reference point that is weakly dominated by each point in the nondominated set~$\PND^s(x)$ for each combination of a scenario~$s\in S$ and a first-stage decision~$x\in\X$. 
We call this property of the reference point \emph{ND-including}. When trying to find such a reference point, it is not advisable to choose unnecessarily large coordinates as this could lead to two undesired effects. First, there is a direct relation between the number of slabs and the horizontal coordinate of the reference point (if using vertical slabs), and a larger coordinate may lead to a longer running time because of the additional (essentially useless) variables. Second, if one of the coordinates of the reference point is much larger than the other, one could deliberately influence the choice of the first-stage decision. For example, if the horizontal coordinate is chosen very large compared to the vertical coordinate, a negative shift of the nondominated points in the vertical direction would make the hypervolume grow more than the same shift in the horizontal direction, and therefore a first-stage decision~$x$ whose nondominated sets~$\PND^s(x)$ (for~$s\in S$) include points with smaller second component is preferred.

The most naive way to find an ND-including reference point is to choose its coordinates as componentwise upper bounds on the objective values attained by any decision in any scenario:
\begin{align*}
	\refpt_\ell \coloneqq \max_{s \in S} \max_{\substack{x^s\in\X \\ y^s \in \Y^s(x^s)}} f_\ell^s(x^s,y^s),\quad \textrm{for } \ell=1,2.
\end{align*}
We call this the \emph{dystopian reference point} or simply \emph{dystopian point}. It can be easily computed by solving the inner maximization problem separately for each scenario.

Finding an ND-including reference point with smallest possible coordinates is not as straightforward. This leads to a three-level optimization problem (see, e.g., \cite{beck2021gentle} for an introduction to multi-level optimization in general). For instance, the optimization problem of finding the horizontal (first) coordinate of the reference point can be described as follows:
\begin{equation}\label{eq:reflexmin}\tag{Ref}
	\begin{aligned}
		\refpt_1 = \max_{s \in S} \quad\max_{x^s \in \X} \quad & f_1^s(x^s,y^s)                                                              \\
		\textrm{s.t.}\quad                                     & y^s \in \underset{y'\in\Y^s(x^s)}{\arg\,\lex\min}\left(f_2^s(x^s,y'),f_1^s(x^s,y')\right).
	\end{aligned}
\end{equation}
This finds, for each scenario and for each first-stage decision (out of possibly infinitely many), the horizontal coordinate of the point corresponding to the lexicographic minimum of~$(f^s_2,f^s_1)$. Out of these values for each scenario and each first-stage decision, the largest one is used. Depending on the problem, there are ways to tackle this three-level problem using appropriate reformulation techniques. We provide a short guide to these in Appendix~\ref{sec-app:ref-analytic-guide}. We call such an ND-including reference point with componentwise smallest possible coordinates the \emph{smallest ND-including reference point} or simply \emph{smallest ND-including point}.

For the problems considered in this paper, Problem~\eqref{eq:reflexmin} can be solved exactly by combinatorial routines: for the selection problem, it is even solvable in polynomial time, and for the assignment problem, an (accelerated) complete enumeration is fast enough for all tested instances except the largest ones ($n=10$ with $|S|=128$), where it exceeded its time limit and the dystopian point was used instead (see Section~\ref{sec:comp-ass-half}). We present the routines in Appendix~\ref{sec-app:ref-enumeration}.

We emphasize again that the selection of the reference point is ultimately a modeling choice and should reflect the decision maker's notion of what constitutes a ``worst relevant'' outcome. Different choices of~$\refpt$ may change the numerical value of the hypervolume and the relative emphasis of regions of the nondominated set, but do not affect the correctness of the proposed method. In applications, the user is therefore free to specify~$\refpt$ according to problem-specific considerations.

\subsection{Analytic determination}\label{sec-app:ref-analytic-guide}

As shown, searching for the smallest ND-including reference point leads to a bi-level problem (Problem~\eqref{eq:reflexmin}) where the inner problem is a lexicographic minimization problem, which can be easily reformulated as a bi-level optimization problem itself, making it three levels in total. We restate the problem here with the three levels and omit the superscript~$s$ on the $y$-variables for better readability:

\begin{equation}\label{eq:threelevel}\tag{3-Level}
	\begin{aligned}
		\max_{s \in S} \quad \max_{x,y}\quad & f^s_1(x,y) \nonumber                         \\
		\text{s.t.}\quad
		                                     & x \in \X, \nonumber                        \\[0.3em]
		                                     & y \in \arg\min_{y'} \; f^s_1(x,y') \nonumber \\
		                                     & \hspace{3.2em}
		\text{s.t.}\quad
		y' \in \arg\min_{y''} \; f^s_2(x,y'') \nonumber                                     \\
		                                     & \hspace{9.0em}
		\text{s.t.}\quad
		y'' \in \Y^s(x).
	\end{aligned}
\end{equation}

Depending on the structure of the underlying problem, there are several possible ways to tackle this problem.

The easiest way to tackle \eqref{eq:reflexmin} is to, instead of replacing the $\lex\min$ part by a third level of optimization, replace it with an adequate linear combination of the two objectives such that
\begin{equation*}
	{\arg\lex\min}_{y\in\Y^s(x)}\{f^s_2(x,y),f^s_1(x,y)\} = {\arg\min}_{y\in\Y^s(x)} \alpha f^s_2(x,y) + \beta f^s_1(x,y),
\end{equation*}
with $\alpha, \beta \in \R_{>0}$. To achieve this, the weights have to be chosen such that a smaller $f^s_2$-value is always preferred over a smaller $f^s_1$-value. Given a point~$p$ and a point~$q$ in the objective space such that $p_2 < q_2$ and $p_1 > q_1$, the model should always choose~$p$ as the one with smaller (better) objective (because of its smaller second component). We claim that this can be achieved by choosing $\alpha$ and $\beta$ such that:
\begin{equation}\label{eq:refpoint-lincomb}
	\frac\beta\alpha < \min_{\substack{y,y'\in\Y^s(x)\\f^s_1(x,y) > f^s_1(x,y')\\f^s_2(x,y') > f^s_2(x,y)}} \frac{f^s_2(x,y')-f^s_2(x,y)}{f^s_1(x,y)-f^s_1(x,y')}.
\end{equation}
Indeed, we have:
\begin{align*}
	\frac\beta\alpha                       & < \frac{q_2-p_2}{p_1-q_1} \\
	\Rightarrow\quad\beta (p_1-q_1)        & < \alpha (q_2-p_2)        \\
	\Rightarrow\quad\alpha p_2 + \beta p_1 & < \alpha q_2 + \beta q_1,
\end{align*}
which is the desired result.

The magnitude of~$\alpha$ and~$\beta$ can be chosen such that the weights and the objective values fall into reasonable and stable numeric ranges. Note that, without additional assumptions, doing this $\lex\min$-reformulation is only possible when the feasible set of second-stage decisions is finite (and even in these cases it can be hard to evaluate the expression on the right-hand side in~\eqref{eq:refpoint-lincomb}).

In case the $\lex\min$-reformulation can be done, we are left with a bi-level (possibly mixed-integer) linear problem. If no integrality constraints are given, the problem can be solved, e.g., with the Kth-best algorithm, as explained in Section~6.1 of \cite{beck2021gentle}.
If there are integrality constraints on some of the variables, this algorithm does not work (unless they can be provably relaxed). For this case, \cite{moore1990mixed} published a branch-and-bound method that involves solving further reformulations of the bi-level problem by using relaxations and bounding on the integer variables. These problems themselves have to be solved using one of the single-level reformulation techniques (usually using KKT reformulations and linearization). We refer the reader to Section~7.3 of \cite{beck2021gentle} for more detail.

In case the $\lex\min$-reformulation cannot be used, other approaches can be tried under certain conditions. For example, if the third (innermost) level of Problem~\eqref{eq:threelevel} has no integrality constraints, we can then get rid of the third level by using reformulations that are otherwise used to transform bi-level into single-level problems, such as the KKT reformulation and then linearizing the quadratic complementarity conditions by adding additional binary variables (see Sections~4.2 and~4.4 in~\cite{beck2021gentle}). This then leads to a bi-level mixed-integer linear problem, and the branch-and-bound procedure by~\cite{moore1990mixed} can be applied as mentioned for the case above. Note that non-integrality of the third level can also be achieved when relaxing the integrality constraints leads to the same problem (e.g., when the coefficient matrix is totally unimodular and the right-hand sides are integral, \cite{schrijver1998theory}).

\subsection{Combinatorial determination}\label{sec-app:ref-enumeration}

As noted at the beginning of this appendix, Problem~\eqref{eq:reflexmin} can be solved exactly by combinatorial routines for the two problems considered in this paper. For the selection problem, we derive a polynomial-time routine that exploits the greedy structure of the second stage. For the assignment problem, we use an enumeration of the first-stage decisions, which is viable for the instance sizes considered, and additionally describe how it can be accelerated by solving assignment subproblems instead of enumerating matchings.

\subsubsection{Selection problem}\label{sec-app:refpoint-selection}

For the selection problem (see Section~\ref{sec:sel-problem}), no enumeration of first-stage decisions is needed at all: Problem~\eqref{eq:reflexmin} can be solved exactly in time~$\mathcal{O}(|S|\,n\log n)$. The reason is the greedy structure of the second stage. For a fixed first-stage decision, the lexicographically optimal second-stage completion simply selects the remaining items with the smallest second-stage weights, so the only way the first stage can influence the second stage is by \emph{blocking} items, i.e., removing them from the pool available to the second stage. This interaction can be resolved by a threshold argument, as made precise in the following proposition. We state it for the first coordinate~$\refpt_1$ of the reference point; the statement for~$\refpt_2$ is obtained by interchanging the roles of~$C$ and~$c^s$ with~$D$ and~$d^s$, respectively.

\begin{proposition}\label{prop:sel-refpoint}
	Let~$s\in S$ be fixed and assume, by renumbering the items, that the pairs~$(d^s_i, c^s_i)$ are lexicographically non-decreasing in~$i$ (see Definition~\ref{def:lex}). For~$t\in\{0,1,\dots,p\}$, define
	\begin{align*}
		V(t)\coloneqq\sum_{i=1}^{t} c^s_i \;+\; \sum_{i=1}^{t-1} \bigl(C_i-c^s_i\bigr)^+ \;+\; \Top_{p-t}\bigl(C_{t+1},\dots,C_n\bigr),
	\end{align*}
	where~$(a)^+\coloneqq\max\{0,a\}$, and where $\Top_m(a_1,\dots,a_k)$ denotes the sum of the $m$ largest values among $a_1,\dots,a_k$, with $\Top_0\coloneqq0$. Then the optimal value of the inner maximization of Problem~\eqref{eq:reflexmin} for scenario~$s$ equals~$\max_{t\in\{0,\dots,p\}}V(t)$.
\end{proposition}
\begin{proof}
	We identify a first-stage decision~$x$ with the set~$X=\{i\in[n]:x_i=1\}$ and write~$k_x\coloneqq|X|$. Since~$\Y^s(x)=\varnothing$ if and only if~$k_x>p$, exactly the sets~$X$ with~$k_x\le p$ are feasible for the inner maximization, and the second stage then selects a set~$Y\subseteq U\coloneqq[n]\setminus X$ of~$k_y\coloneqq p-k_x$ of the remaining items. As the contribution of~$X$ to both objectives is fixed, a second-stage selection is optimal for the inner lexicographic minimization if and only if it lexicographically minimizes~$\bigl(\sum_{i\in Y}d^s_i,\,\sum_{i\in Y}c^s_i\bigr)$ over all~$Y\subseteq U$ with~$|Y|=k_y$.

	\emph{Step 1 (greedy second stage).} We claim that the \emph{greedy selection}~$\hat Y(X)$, consisting of the~$k_y$ items of~$U$ with the smallest indices (recall the sorting of the items), is optimal, and that every optimal selection attains the same pair of sums. To see this, note first that, for any feasible~$Y$, the~$j$-th smallest~$d^s$-value in~$Y$ is greater than or equal to the~$j$-th smallest~$d^s$-value in~$U$ for every~$j\in[k_y]$. Hence, $\sum_{i\in Y}d^s_i$ is minimal if and only if the~$d^s$-values of~$Y$ are exactly the~$k_y$ smallest~$d^s$-values in~$U$, counted with multiplicity. Letting~$\tau$ denote the~$k_y$-th smallest~$d^s$-value in~$U$, this means that every minimizer contains all items~$i\in U$ with~$d^s_i<\tau$ and is completed by items with~$d^s_i=\tau$. The secondary sum~$\sum_{i\in Y}c^s_i$ is then minimized exactly by completing with tied items having the smallest~$c^s$-values, which is precisely what the greedy selection does (by the sorting of the items). Since both components of a lexicographic minimum are unique, $f_1^s(x,y)$ has the same value for every optimal second-stage selection~$y$, so we may evaluate the inner maximization using~$\hat Y(X)$.

	\emph{Step 2 (upper bound).} Let~$X$ with~$k_x\le p$ be given. If~$k_y=0$, then $f_1^s(x,y)=\sum_{i\in X}C_i\le\Top_p(C_1,\dots,C_n)=V(0)$. So let~$k_y\ge1$ and let~$t$ be the largest index of an item in~$\hat Y(X)$. By the greedy rule, every item~$i\le t$ belongs to either~$X$ or~$\hat Y(X)$, and item~$t$ belongs to~$\hat Y(X)$. Hence, $|X\cap[t]|=t-k_y$ and, consequently, $|X\setminus[t]|=k_x-(t-k_y)=p-t$, which in particular yields~$t\le p$. We obtain
	\begin{align*}
		f_1^s(x,\hat y)&=\sum_{i\in[t]\setminus X}c^s_i+\sum_{i\in X\cap[t]}C_i+\sum_{i\in X\setminus[t]}C_i
		=\sum_{i=1}^{t}c^s_i+\sum_{i\in X\cap[t]}\bigl(C_i-c^s_i\bigr)+\sum_{i\in X\setminus[t]}C_i\\
		&\le\sum_{i=1}^{t}c^s_i+\sum_{i=1}^{t-1}\bigl(C_i-c^s_i\bigr)^++\Top_{p-t}\bigl(C_{t+1},\dots,C_n\bigr)=V(t),
	\end{align*}
	where the inequality uses~$t\notin X$ and~$|X\setminus[t]|=p-t$.

	\emph{Step 3 (achievability).} Let~$t\in\{0,\dots,p\}$. For~$t=0$, choosing~$X$ as a set of~$p$ items with the largest values~$C_i$ yields the value~$V(0)$. For~$t\ge1$, set~$A\coloneqq\{i\in[t-1]:C_i>c^s_i\}$, let~$B$ consist of~$p-t$ items~$i>t$ with the largest values~$C_i$, and put~$X\coloneqq A\cup B$. Then~$k_x=|A|+p-t\le p-1$, so~$X$ is feasible with~$k_y=t-|A|\ge1$, and the greedy selection is~$\hat Y(X)=[t]\setminus A$, as these are the~$k_y$ items with the smallest indices outside of~$X$. By Step~1, $\hat Y(X)$ is an optimal second-stage selection for~$X$ (and every optimal selection yields the same value of~$f_1^s$), so the value attained by~$X$ in the inner maximization is~$f_1^s(x,\hat y)$. Using~$\sum_{i\in A}C_i=\sum_{i\in A}\bigl(C_i-c^s_i\bigr)+\sum_{i\in A}c^s_i$, this gives
	\begin{align*}
		f_1^s(x,\hat y)=\sum_{i\in A}C_i+\sum_{i\in B}C_i+\sum_{i\in[t]\setminus A}c^s_i
		=\sum_{i=1}^{t}c^s_i+\sum_{i\in A}\bigl(C_i-c^s_i\bigr)+\sum_{i\in B}C_i=V(t). &\qedhere
	\end{align*}
\end{proof}

Algorithm~\ref{alg:sel-refpoint} directly implements Proposition~\ref{prop:sel-refpoint}. Regarding its running time, note that, after sorting, the first two sums in the definition of~$V(t)$ can be maintained incrementally in constant time per step when iterating over~$t=0,1,\dots,p$. The tail sums~$\Top_{p-t}(C_{t+1},\dots,C_n)$ for all~$t\in\{0,\dots,p\}$ can be computed in a single backward pass by maintaining a max-heap of the not-yet-counted values together with a running sum. We initialize the heap with~$C_{p+1},\dots,C_n$ and the running sum with zero, which is the value~$\Top_0(C_{p+1},\dots,C_n)$ for~$t=p$, and then process~$t=p-1,p-2,\dots,0$: in each step, we insert~$C_{t+1}$ into the heap, extract the maximum, and add it to the running sum, which then equals~$\Top_{p-t}(C_{t+1},\dots,C_n)$. In total, this yields a running time of~$\mathcal{O}(n\log n)$ per scenario and coordinate, i.e., $\mathcal{O}(|S|\,n\log n)$ overall.

\begin{algorithm}
	\caption{Routine to find the smallest ND-including reference point for the selection problem.}\label{alg:sel-refpoint}
	\begin{algorithmic}[1]
		\Statex \textbf{Input:} $n\in\N$ and $p\in\{0,\dots,n\}$, scenario set $S$, first-stage weights $C,D\in\R^{n}$, second-stage weights $c^s,d^s\in\R^{n}$ for all $s\in S$
		\Statex \textbf{Output:} Reference point $\refpt=(\refpt_1,\refpt_2)$

		\State $(\refpt_1,\refpt_2)\gets(-\infty,-\infty)$
		\For{$s\in S$ and $\ell\in\{1,2\}$}
		\If{$\ell=1$}
		\State $(F,f,g)\gets(C,c^s,d^s)$
		\Else
		\State $(F,f,g)\gets(D,d^s,c^s)$
		\EndIf
		\State renumber the items such that the pairs $(g_i,f_i)$ are lexicographically non-decreasing in~$i$
		\For{$t=0,1,\dots,p$}
		\State $V(t)\gets\sum_{i=1}^{t} f_i+\sum_{i=1}^{t-1}\bigl(F_i-f_i\bigr)^++\Top_{p-t}\bigl(F_{t+1},\dots,F_n\bigr)$
		\EndFor
		\State $\refpt_\ell\gets\max\bigl\{\refpt_\ell,\ \max_{t\in\{0,\dots,p\}}V(t)\bigr\}$
		\EndFor
		\State \Return $(\refpt_1,\refpt_2)$
	\end{algorithmic}
\end{algorithm}

The routine is cheap in practice: it never exceeds $13.5$\,ms and stays below $0.25\%$ of the total algorithm time on every single run. Table~\ref{tab:sel-refpoints-full} lists the resulting reference point for every instance of Section~\ref{sec:comp-design}, next to the dystopian point of Appendix~\ref{sec-app:refpoint} for comparison. It is precisely this that makes the largest instances of our grid tractable at all; a naive full enumeration of all first-stage selections, whose effort grows like $|S|\cdot\sum_{k\le p}\binom{n}{k}$, would visit about $6.2\cdot 10^{11}$ selections at $(n,p)=(40,20)$ and about $6.4\cdot 10^{17}$ at $(60,30)$, and be entirely out of reach.

\begingroup\scriptsize\setlength{\tabcolsep}{1.5pt}
\begin{longtable}{lll|cc|cc|cc}
\hline
 & & & \multicolumn{2}{c|}{$|S|=4$} & \multicolumn{2}{c|}{$|S|=16$} & \multicolumn{2}{c}{$|S|=128$} \\
$n$ & $p$ & Seed & $\mathrm{ND}$ & $\mathrm{dys}$ & $\mathrm{ND}$ & $\mathrm{dys}$ & $\mathrm{ND}$ & $\mathrm{dys}$ \\
\hline
\endfirsthead
\multicolumn{9}{l}{\footnotesize Table~\ref{tab:sel-refpoints-full}, continued.} \\
\hline
 & & & \multicolumn{2}{c|}{$|S|=4$} & \multicolumn{2}{c|}{$|S|=16$} & \multicolumn{2}{c}{$|S|=128$} \\
$n$ & $p$ & Seed & $\mathrm{ND}$ & $\mathrm{dys}$ & $\mathrm{ND}$ & $\mathrm{dys}$ & $\mathrm{ND}$ & $\mathrm{dys}$ \\
\hline
\endhead
\hline
\endfoot
\hline
\caption{Reference points for each instance $(n, p, \text{seed}, |S|)$ of the selection problem: the smallest ND-including point~$\refpt_{\mathrm{ND}}$ (column ``ND'') computed by the combinatorial routine of Appendix~\ref{sec-app:refpoint-selection}, which the algorithm uses throughout, and the dystopian point~$\refpt_{\mathrm{dys}}$ (column ``dys'') obtained by maximizing each objective separately in every scenario with both stages free. Both depend only on the instance, not on the algorithm parameters. By construction $\refpt_{\mathrm{ND}}\le\refpt_{\mathrm{dys}}$ componentwise
.}\label{tab:sel-refpoints-full} \\
\endlastfoot
10 & 2 & 1 & $(185,202)$ & $(200,230)$ & $(185,219)$ & $(208,230)$ & $(203,246)$ & $(234,266)$ \\
 &  & 2 & $(171,193)$ & $(221,238)$ & $(189,193)$ & $(229,256)$ & $(217,209)$ & $(253,268)$ \\
 &  & 3 & $(155,198)$ & $(193,199)$ & $(170,208)$ & $(227,234)$ & $(204,221)$ & $(253,259)$ \\
 &  & 4 & $(190,214)$ & $(241,216)$ & $(195,214)$ & $(250,282)$ & $(238,225)$ & $(269,282)$ \\
 &  & 5 & $(206,215)$ & $(222,236)$ & $(206,219)$ & $(222,274)$ & $(231,277)$ & $(244,293)$ \\
 & 4 & 1 & $(324,361)$ & $(344,412)$ & $(333,394)$ & $(385,421)$ & $(370,441)$ & $(424,476)$ \\
 &  & 2 & $(337,367)$ & $(397,430)$ & $(337,367)$ & $(397,463)$ & $(411,367)$ & $(444,463)$ \\
 &  & 3 & $(286,340)$ & $(339,372)$ & $(340,380)$ & $(393,435)$ & $(363,420)$ & $(397,481)$ \\
 &  & 4 & $(363,376)$ & $(422,406)$ & $(371,376)$ & $(426,489)$ & $(418,400)$ & $(483,489)$ \\
 &  & 5 & $(358,389)$ & $(390,408)$ & $(358,389)$ & $(390,476)$ & $(390,445)$ & $(455,494)$ \\
 & 5 & 1 & $(385,435)$ & $(406,497)$ & $(405,469)$ & $(451,506)$ & $(442,520)$ & $(512,551)$ \\
 &  & 2 & $(459,446)$ & $(463,515)$ & $(459,446)$ & $(463,549)$ & $(502,489)$ & $(523,549)$ \\
 &  & 3 & $(343,413)$ & $(401,454)$ & $(408,453)$ & $(457,529)$ & $(431,503)$ & $(464,569)$ \\
 &  & 4 & $(444,450)$ & $(498,488)$ & $(452,450)$ & $(512,572)$ & $(507,482)$ & $(581,572)$ \\
 &  & 5 & $(431,480)$ & $(469,489)$ & $(431,480)$ & $(469,562)$ & $(483,529)$ & $(536,582)$ \\
\hline
20 & 2 & 1 & $(191,204)$ & $(224,230)$ & $(214,231)$ & $(224,232)$ & $(216,246)$ & $(252,269)$ \\
 &  & 2 & $(184,175)$ & $(258,238)$ & $(220,201)$ & $(258,257)$ & $(246,212)$ & $(264,268)$ \\
 &  & 3 & $(178,198)$ & $(225,204)$ & $(250,208)$ & $(260,234)$ & $(250,237)$ & $(271,259)$ \\
 &  & 4 & $(206,179)$ & $(257,245)$ & $(206,190)$ & $(257,282)$ & $(238,249)$ & $(286,282)$ \\
 &  & 5 & $(198,220)$ & $(239,272)$ & $(213,220)$ & $(274,274)$ & $(231,277)$ & $(274,303)$ \\
 & 6 & 1 & $(518,558)$ & $(575,619)$ & $(527,607)$ & $(599,660)$ & $(542,620)$ & $(686,704)$ \\
 &  & 2 & $(502,519)$ & $(617,629)$ & $(529,519)$ & $(617,682)$ & $(576,549)$ & $(688,689)$ \\
 &  & 3 & $(493,494)$ & $(601,547)$ & $(580,528)$ & $(693,641)$ & $(580,606)$ & $(693,672)$ \\
 &  & 4 & $(553,552)$ & $(659,655)$ & $(560,556)$ & $(659,708)$ & $(594,599)$ & $(738,708)$ \\
 &  & 5 & $(523,581)$ & $(625,626)$ & $(566,599)$ & $(659,715)$ & $(580,686)$ & $(677,735)$ \\
 & 10 & 1 & $(812,881)$ & $(895,970)$ & $(840,916)$ & $(929,1044)$ & $(894,963)$ & $(1069,1080)$ \\
 &  & 2 & $(825,828)$ & $(941,974)$ & $(825,828)$ & $(941,1018)$ & $(949,859)$ & $(1055,1044)$ \\
 &  & 3 & $(774,810)$ & $(907,855)$ & $(905,824)$ & $(1031,999)$ & $(905,936)$ & $(1034,1008)$ \\
 &  & 4 & $(855,862)$ & $(975,986)$ & $(890,862)$ & $(975,1025)$ & $(941,914)$ & $(1092,1066)$ \\
 &  & 5 & $(827,909)$ & $(954,967)$ & $(895,932)$ & $(982,1066)$ & $(895,1002)$ & $(989,1099)$ \\
\newpage
40 & 2 & 1 & $(192,185)$ & $(233,246)$ & $(193,201)$ & $(247,246)$ & $(224,246)$ & $(277,285)$ \\*
 &  & 2 & $(219,183)$ & $(271,280)$ & $(220,185)$ & $(271,280)$ & $(220,208)$ & $(278,280)$ \\
 &  & 3 & $(192,198)$ & $(250,219)$ & $(236,209)$ & $(268,241)$ & $(248,221)$ & $(292,261)$ \\
 &  & 4 & $(201,210)$ & $(273,245)$ & $(201,210)$ & $(273,285)$ & $(206,225)$ & $(286,285)$ \\
 &  & 5 & $(183,193)$ & $(248,272)$ & $(208,219)$ & $(288,274)$ & $(230,226)$ & $(288,303)$ \\
 & 11 & 1 & $(929,930)$ & $(1063,1102)$ & $(951,948)$ & $(1131,1198)$ & $(1002,1058)$ & $(1233,1247)$ \\
 &  & 2 & $(943,916)$ & $(1118,1172)$ & $(943,916)$ & $(1153,1232)$ & $(966,919)$ & $(1233,1232)$ \\
 &  & 3 & $(938,870)$ & $(1087,993)$ & $(1033,938)$ & $(1249,1115)$ & $(1073,971)$ & $(1297,1207)$ \\
 &  & 4 & $(951,956)$ & $(1247,1174)$ & $(951,956)$ & $(1247,1249)$ & $(1004,1006)$ & $(1247,1249)$ \\
 &  & 5 & $(964,943)$ & $(1100,1111)$ & $(979,964)$ & $(1177,1251)$ & $(992,1070)$ & $(1186,1270)$ \\
 & 20 & 1 & $(1565,1615)$ & $(1752,1813)$ & $(1610,1645)$ & $(1851,1943)$ & $(1680,1742)$ & $(1966,2038)$ \\
 &  & 2 & $(1523,1609)$ & $(1817,1907)$ & $(1523,1609)$ & $(1864,1972)$ & $(1637,1609)$ & $(1957,1972)$ \\
 &  & 3 & $(1559,1487)$ & $(1754,1665)$ & $(1817,1601)$ & $(1981,1804)$ & $(1817,1636)$ & $(2024,1907)$ \\
 &  & 4 & $(1604,1644)$ & $(1976,1885)$ & $(1604,1644)$ & $(1976,1980)$ & $(1679,1735)$ & $(2021,2074)$ \\
 &  & 5 & $(1607,1659)$ & $(1839,1842)$ & $(1666,1728)$ & $(1910,2012)$ & $(1730,1780)$ & $(1937,2073)$ \\
\hline
60 & 2 & 1 & $(192,191)$ & $(233,277)$ & $(193,213)$ & $(248,277)$ & $(224,246)$ & $(277,285)$ \\
 &  & 2 & $(223,183)$ & $(280,280)$ & $(250,235)$ & $(301,287)$ & $(250,235)$ & $(301,312)$ \\
 &  & 3 & $(192,203)$ & $(251,233)$ & $(218,203)$ & $(294,256)$ & $(248,221)$ & $(297,266)$ \\
 &  & 4 & $(220,210)$ & $(275,245)$ & $(220,210)$ & $(275,285)$ & $(226,210)$ & $(286,285)$ \\
 &  & 5 & $(183,197)$ & $(248,274)$ & $(193,208)$ & $(288,276)$ & $(224,234)$ & $(288,320)$ \\
 & 16 & 1 & $(1295,1368)$ & $(1480,1642)$ & $(1325,1411)$ & $(1567,1745)$ & $(1385,1458)$ & $(1740,1818)$ \\
 &  & 2 & $(1387,1314)$ & $(1759,1697)$ & $(1410,1350)$ & $(1759,1897)$ & $(1411,1350)$ & $(1779,1897)$ \\
 &  & 3 & $(1362,1301)$ & $(1566,1495)$ & $(1469,1361)$ & $(1791,1688)$ & $(1471,1400)$ & $(1851,1730)$ \\
 &  & 4 & $(1403,1361)$ & $(1915,1680)$ & $(1465,1384)$ & $(1915,1787)$ & $(1552,1468)$ & $(1915,1823)$ \\
 &  & 5 & $(1318,1399)$ & $(1539,1707)$ & $(1358,1400)$ & $(1671,1823)$ & $(1358,1471)$ & $(1737,1823)$ \\
 & 30 & 1 & $(2233,2429)$ & $(2505,2725)$ & $(2276,2491)$ & $(2643,2961)$ & $(2383,2585)$ & $(2837,3039)$ \\
 &  & 2 & $(2309,2378)$ & $(2867,2865)$ & $(2343,2378)$ & $(2867,3064)$ & $(2438,2394)$ & $(2940,3064)$ \\
 &  & 3 & $(2351,2270)$ & $(2637,2551)$ & $(2611,2425)$ & $(2916,2822)$ & $(2611,2498)$ & $(3036,2863)$ \\
 &  & 4 & $(2420,2400)$ & $(3055,2789)$ & $(2539,2408)$ & $(3055,2930)$ & $(2602,2571)$ & $(3091,3049)$ \\
 &  & 5 & $(2285,2419)$ & $(2637,2849)$ & $(2375,2552)$ & $(2807,3016)$ & $(2435,2552)$ & $(2879,3021)$ \\
\end{longtable}
\endgroup

\subsubsection{Assignment problem}\label{sec-app:refpoint-assignment}

For the assignment problem (see Section~\ref{sec:ass-problem}), Algorithm~\ref{alg:ass-refpoint-fast} is applied to find the smallest ND-including reference point. As in Algorithm~\ref{alg:sel-refpoint}, both coordinates are computed by the same routine via the roles~$(F,f,g)$. We assume that an edge~$e$ can equivalently be written via the nodes it connects as~$(i,j)$, and we let~$\texttt{comb}(A,k)$ denote the set of all $k$-element subsets of~$A$, each written as a tuple in a fixed order.

\begin{algorithm}
	\caption{Routine to find the smallest ND-including reference point for the assignment problem.}\label{alg:ass-refpoint-fast}
	\begin{algorithmic}[1]
		\Statex \textbf{Input:} $n \in \N$ and $p \in \{0,\dots,n\}$, scenario set $S$, first-stage weights $C,D \in \R^{n^2}$, second-stage weights $c^s, d^s\in\R^{n^2}$ for all $s\in S$
		\Statex \textbf{Output:} Reference point $\refpt = (\refpt_1, \refpt_2)$

		\State $(\refpt_1, \refpt_2) \gets (-\infty,-\infty)$
		\For{$s\in S$ and $\ell\in\{1,2\}$}
		\If{$\ell=1$}
		\State $(F,f,g)\gets(C,c^s,d^s)$
		\Else
		\State $(F,f,g)\gets(D,d^s,c^s)$
		\EndIf
		\For{$k \in \{p,p+1,\dots,n\}$, $X^L \in \texttt{comb}([n],k)$, and $X^R \in \texttt{comb}([n], k)$}
		\State $F_x \gets \max\bigl\{\textstyle\sum_{e\in M_x} F_e : M_x \text{ perfect matching between } X^L \text{ and } X^R\bigr\}$
		\If{$k = n$} \Comment{second stage is empty}
		\State $\refpt_\ell \gets \max\{\refpt_\ell,\ F_x\}$
		\State \textbf{continue}
		\EndIf
		\State $(g_y, f_y) \gets \lex\min\bigl\{\bigl(\textstyle\sum_{e\in M_y} g_e,\ \sum_{e\in M_y} f_e\bigr) : M_y \text{ perfect matching}$
		\Statex \hspace{7.65em}$\text{between } [n]\setminus X^L \text{ and } [n]\setminus X^R\bigr\}$
		\State $\refpt_\ell \gets \max\{\refpt_\ell,\ F_x + f_y\}$
		\EndFor
		\EndFor
		\State \Return $(\refpt_1, \refpt_2)$
	\end{algorithmic}
\end{algorithm}

Rather than enumerating all first- and second-stage matchings, the routine enumerates only the pairs~$(X^L, X^R)$ of vertex sets covered by the first-stage matching and solves two linear assignment problems for each. This is justified by two observations. First, the set of feasible second-stage completions---the perfect matchings between the vertices left uncovered by the first stage---depends only on the vertex sets~$X^L$ and~$X^R$ covered by the first-stage matching, and not on how these vertices are matched to each other. Since, in addition, both components of a lexicographic minimum are unique, the contribution of the second stage to the objective is the same for every first-stage matching covering~$X^L$ and~$X^R$. For fixed~$X^L$ and~$X^R$, the best objective value is therefore attained by a first-stage matching between~$X^L$ and~$X^R$ of maximum $F$-weight, which can be computed as a linear assignment problem in time~$\mathcal{O}(n^3)$ instead of enumerating all~$k!$ pairings. Second, finding the lexicographically minimal second-stage completion is a lexicographic assignment problem, which can likewise be solved in time~$\mathcal{O}(n^3)$: first, a perfect matching between~$[n]\setminus X^L$ and~$[n]\setminus X^R$ minimizing~$\sum_{e} g_e$ is computed. By complementary slackness, a perfect matching is optimal for this problem if and only if it uses only edges that are tight with respect to an optimal dual solution (see, e.g., \cite{schrijver1998theory}), so minimizing~$\sum_{e} f_e$ over the perfect matchings in the subgraph of tight edges---again a linear assignment problem---yields the lexicographic minimum. Alternatively, for integer weights, both steps can be combined into a single assignment problem with the aggregated weights~$W g_e + f_e$ for a sufficiently large constant~$W$.

A naive enumeration of all first- and second-stage matchings would consider
\begin{align*}
	\sum_{k=p}^{n}\binom{n}{k}\cdot\frac{n!}{(n-k)!}\cdot(n-k)! \;\le\; 2^n\, n!
\end{align*}
combinations per scenario and coordinate. In contrast, Algorithm~\ref{alg:ass-refpoint-fast} enumerates only the
\begin{align*}
	\sum_{k=p}^{n}\binom{n}{k}^2\;\le\;\binom{2n}{n}
\end{align*}
pairs of covered vertex sets and solves two assignment problems for each, reducing the overall effort per scenario and coordinate from~$\mathcal{O}(2^n\, n!\, n)$ to~$\mathcal{O}\bigl(\binom{2n}{n} n^3\bigr)$. 

Even in this accelerated form, the effort still grows exponentially in~$n$, and it does so considerably faster than for the dystopian reference point defined at the beginning of this appendix, which only requires maximizing each objective separately for every scenario and therefore scales very well with the instance size. In fact, for the largest instance sizes we considered in Section~\ref{sec:comp-ass-half}, the computation of the smallest ND-including point hit a time limit that we chose for the experiments. In these cases, we decided to fall back on the dystopian reference point, at the price of a reference point lying further out. This is possible since, as noted above, the choice of~$\refpt$ affects the numerical value of the hypervolume, but not the correctness of the method.

\end{appendices}

\end{document}